\documentstyle[pstricks,pst-node,cite]{amsart}

\newtheorem{theorem}{Theorem}[section]
\newtheorem{lemma}[theorem]{Lemma}
\newtheorem{conjecture}[theorem]{Conjecture}

\newtheorem{exercise}[theorem]{Exercise}

\newenvironment{proof}{{\bf Proof: }}{\qed}

\newcommand{\C}{{\Bbb C}}
\newcommand{\Z}{{\Bbb Z}}
\newcommand{\F}{{\Bbb F}}

\newcommand{\opname}[1]{\mathop{\fam0#1}}

\newcommand{\tensor}{\otimes}
\newcommand{\injects}{\hookrightarrow}
\newcommand{\Inv}{\mbox{\rm Inv}}
\newcommand{\End}{\opname{End}}
\newcommand{\Hom}{\opname{Hom}}
\renewcommand{\sl}{\opname{sl}}
\newcommand{\so}{\opname{so}}
\renewcommand{\sp}{\opname{sp}}

\newcommand{\Spin}{\opname{Spin}}
\newcommand{\id}{\opname{id.}}
\newcommand{\wt}{\opname{wt}}
\renewcommand{\j}{\Join}
\newcommand{\Ie}{{\em I.e.}}
\newcommand{\ie}{{\em i.e.}}
\newcommand{\df}{\em\emgreen}

\newcounter{fignum}
\ifx\blackandwhite\undefined
  \psset{linecolor=blue}\newrgbcolor{darkred}{.5 0 0}\newrgbcolor{emgreen}{0 .5 0}
\else
  \psset{linecolor=black}\newgray{darkred}{0}\newgray{emgreen}{0}
\fi

\psset{linewidth=.4pt,dash=3pt 3pt,doublesep=.05,dotsize=1pt 5}
\SpecialCoor
\newcommand{\psgoesto}{\hspace{.5cm}\pspicture[.5](0,-.1)(1,.1)
    \psline[linecolor=black]{->}(0,0)(1,0)
    \endpspicture\hspace{.5cm}}
\newcommand{\vertvert}{\pspicture[.4](-.6,-.5)(.6,.5)
\pscircle[linecolor=darkred,linestyle=dashed](0,0){.5}
\psbezier(.5;45)(.25;45)(.25;315)(.5;315)
\psbezier(.5;225)(.25;225)(.25;135)(.5;135)
\endpspicture}
\newcommand{\horizhoriz}{\pspicture[.4](-.6,-.5)(.6,.5)
\pscircle[linecolor=darkred,linestyle=dashed](0,0){.5}
\psbezier(.5;45)(.25;45)(.25;135)(.5;135)
\psbezier(.5;225)(.25;225)(.25;315)(.5;315)
\endpspicture}
\newcommand{\rcrossing}{\pspicture[.4](-.6,-.5)(.6,.5)
\qline(.5;135)(.5;315)
\psline[border=.1](.5;45)(.5;225)
\endpspicture}
\newcommand{\vv}{\pspicture[.4](-.6,-.5)(.6,.5)
\psbezier(.5;45)(.25;45)(.25;315)(.5;315)
\psbezier(.5;225)(.25;225)(.25;135)(.5;135)
\endpspicture}
\newcommand{\hh}{\pspicture[.4](-.6,-.5)(.6,.5)
\psbezier(.5;45)(.25;45)(.25;135)(.5;135)
\psbezier(.5;225)(.25;225)(.25;315)(.5;315)
\endpspicture}
\newcommand{\btwohvert}{\pspicture[.4](-.6,-.7)(.6,.7)
\qline(0, .25)(.35, .6)\qline(0, .25)(-.35, .6)
\qline(0,-.25)(.35,-.6)\qline(0,-.25)(-.35,-.6)
\psline[doubleline=true](0,-.25)(0,.25)
\endpspicture}
\newcommand{\btwohhoriz}{\pspicture[.4](-.8,-.5)(.8,.5)
\qline( .25,0)( .6,.35)\qline( .25,0)( .6,-.35)
\qline(-.25,0)(-.6,.35)\qline(-.25,0)(-.6,-.35)
\psline[doubleline=true](-.25,0)(.25,0)
\endpspicture}
\newcommand{\gtwohvert}{\pspicture[.4](-.6,-.7)(.6,.7)
\qline(0, .25)(.433, .5)\qline(0, .25)(-.433, .5)
\qline(0,-.25)(.433,-.5)\qline(0,-.25)(-.433,-.5)
\qline(0,-.25)(0,.25)
\endpspicture}
\newcommand{\gtwohhoriz}{\pspicture[.4](-.8,-.5)(.8,.5)
\qline( .25,0)( .5,.433)\qline( .25,0)( .5,-.433)
\qline(-.25,0)(-.5,.433)\qline(-.25,0)(-.5,-.433)
\qline(-.25,0)(.25,0)
\endpspicture}
\newcommand{\doublehvert}{\pspicture[.4](-.6,-.7)(.6,.7)
\psline[doubleline=true](0, .25)(.35, .6)\qline(0, .25)(-.35, .6)
\psline[doubleline=true](0,-.25)(-.35,-.6)\qline(0,-.25)(.35,-.6)
\qline(0,-.25)(0,.25)
\endpspicture}
\newcommand{\doublehhoriz}{\pspicture[.4](-.8,-.5)(.8,.5)
\psline[doubleline=true]( .25,0)( .6,.35)\qline( .25,0)( .6,-.35)
\psline[doubleline=true](-.25,0)(-.6,-.35)\qline(-.25,0)(-.6,.35)
\qline(-.25,0)(.25,0)
\endpspicture}
\newcommand{\singleloop}{\pspicture[.4](-.6,-.5)(.6,.5)
\pscircle(0,0){.4}
\endpspicture}
\newcommand{\doubleloop}{\pspicture[.4](-.6,-.5)(.6,.5)
\pscircle[doubleline=true](0,0){.4}
\endpspicture}
\newcommand{\doublesquare}{\pspicture[.4](-.9,-.9)(.9,.9)
\psline[doubleline=true](.4;45)(.8;45)
\psline[doubleline=true](.4;135)(.8;135)
\psline[doubleline=true](.4;225)(.8;225)
\psline[doubleline=true](.4;315)(.8;315)
\psline(.4;45)(.4;135)\psline(.4;135)(.4;225)
\psline(.4;225)(.4;315)\psline(.4;315)(.4;45)
\endpspicture}
\newcommand{\semidoublesquare}{\pspicture[.4](-.9,-.9)(.9,.9)
\psline[doubleline=true](.4;45)(.8;45)
\psline(.4;135)(.8;135)
\psline[doubleline=true](.4;225)(.8;225)
\psline(.4;315)(.8;315)
\psline(.4;45)(.4;135)\psline(.4;135)(.4;225)
\psline(.4;225)(.4;315)\psline(.4;315)(.4;45)
\endpspicture}
\newcommand{\pentanode}{\pnode(.4;90){a1}\pnode(.4;162){a2}\pnode(.4;234){a3}
\pnode(.4;306){a4}\pnode(.4;18){a5}
\pnode(.9;90){b1}\pnode(.9;162){b2}\pnode(.9;234){b3}
\pnode(.9;306){b4}\pnode(.9;18){b5}}
\newcommand{\btwox}{\pspicture[.4](-.6,-.5)(.6,.5)
\qline(.5;45)(.5;225)\qline(.5;135)(.5;315)
\endpspicture}
\newcommand{\middlearrow}{\lput{:U}{\pspicture(0,0)(0,0)
\psline[arrows=->,arrowscale=1.5](2.2pt,0)(2.3pt,0)\endpspicture}}
\newcommand{\littley}{
\qline(.2887;30)(0,0)\qline(.2887;150)(0,0)\qline(.2887;270)(0,0)}
\newcommand{\littlelam}{
\qline(.2887;90)(0,0)\qline(.2887;210)(0,0)\qline(.2887;330)(0,0)}

\begin{document}

\author{Greg Kuperberg}
\address{Department of Mathematics \\ Yale University \\
	New Haven, CT 06520}
\email{greg@@math.yale.edu}
\thanks{The author was supported by an NSF
	Postdoctoral Fellowship, grant \#DMS-9107908.}
\title{Spiders for rank 2 Lie algebras}
\begin{abstract}
A spider is an axiomatization of the representation theory of a group, quantum
group, Lie algebra, or other group or group-like object.  It is also known as
a spherical category, or a strict, monoidal category with a few extra
properties, or by several other names.  A recently useful point of view,
developed by other authors, of the representation theory of $\sl(2)$ has been
to present it as a spider by generators and relations.  That is, one has an
algebraic spider, defined by invariants of linear representations, and one
identifies it as isomorphic to a combinatorial spider, given by generators
and relations.  We generalize this approach to the rank 2 simple Lie
algebras, namely $A_2$, $B_2$, and $G_2$.  Our combinatorial rank 2 spiders
yield bases for invariant spaces which are probably related to Lusztig's
canonical bases, and they are useful for computing quantities such as
generalized $6j$-symbols and quantum link invariants.  Their definition
originates in definitions of the rank 2 quantum link invariants that were
discovered independently by the author and Francois Jaeger.
\end{abstract}
\maketitle

\section{Introduction}
\label{sintroduction}

One of the problems of classical invariant theory is to characterize, for all
$n$-tuples $V_1,\ldots,V_n$ of finite-dimensional, irreducible
representations over $\C$ of a compact group $G$ or simple Lie algebra $\frak
g$, the vector space of multilinear functions
$$f:V_1\times V_2\times\ldots\times V_n\to\C$$
which
are invariant under the action of $G$ or $\frak g$. In more modern
terminology, the problem is to characterize the dual vector
space of invariant tensors
$$\Inv(V_1\tensor V_2\tensor\ldots\tensor V_n),$$
or just $\Inv(V)$ if $V$ is a tensor product of irreducibles.  (Also, instead
of working over $\C$, one might work over some other field $\F$ of
characteristic 0.)  Of course, for a simple Lie algebra $\frak g$, the
dimension of such a vector space is given by Cartan-Weyl character theory.
But it is also useful to consider operations on invariant spaces such as
tensor products and contractions. Even for finite-dimensional, simple Lie
algebras over $\C$, these operations are not completely understood.

Interest in multilinear invariant theory was revived after the discovery of
quantum groups.  A quantum Lie group is a non-commutative, non-cocommutative
Hopf algebra $U_q(\frak g)$ which is a deformation of the universal
enveloping algebra $U(\frak g)$.  In fact, the reprentation theory of a
quantum group is just as important than the quantum group itself and may also
be defined using loop groups or conformal quantum field theory. Among other
uses of quantum groups, a quantum invariant of a tangle takes is a vector in
an invariant space $\Inv(V)$ where $V$ depends on the boundary of the
tangle and other data.

Because of non-cocommutativity, the switching map
$$\tau:V \tensor W \to W \tensor V$$
given by
$$x \tensor y \mapsto y \tensor x$$
is in general not an (equivariant) map between quantum group representations.
Thus, there is no natural symmetric group action on invariant spaces with $n$
tensor factors with $n>2$.  (There is often a braid group action.) However,
the following operations exist and are natural:
\begin{description} 
\item[1.] Tensor product:
$$\Inv(V) \tensor \Inv(W) \injects \Inv(V \tensor W)$$
\item[2.] Cyclic permutation:
$$\Inv(V \tensor W) \to
\Inv(W \tensor V)$$
\item[3.] Contraction:
$$\Inv(V \tensor V^* \tensor W) \to \Inv(W)$$
\end{description}
Since $V$ and $W$ may themselves be tensor produts, cyclic permutation of two
tensor factors yields cyclic permutation of $n$ tensor factors, but not
general permutations.  Also, contraction must be interpreted carefully,
because in a quantum group, one must reverse order when taking
duals:  $(V \tensor W)^* \cong W^* \tensor V^*$.

A spider is an abstraction of a representation theory with these
three operations.  It is a collection of vector spaces, or perhaps
modules or sets, to be thought of as invariant spaces, together
with abstract operations called join, rotation, and stitch,
to be thought of as tensor product, rotation, and contraction.
It is both convenient and conceptually important to depict these
operations with certain planar graphs.  These graphs are called
webs, hence the term ``spider''.

Another motivation of the spider operations is that they really describe the
entire equivariant tensor category of representations of a group, Lie
algebra, or quantum group.  In general,
$$\Hom(A,B) \cong \Inv(A^* \tensor B),$$
the tensor product of two homomorphisms can be defined in terms of tensor
product and cyclic permutation of invariants, and composition of homomorphisms can be
defined in terms of tensor product and contraction.  Contrariwise, the spider
operations can be defined in terms of tensor product and composition of morphisms.
For this reason, spiders are sometimes defined as a type of (non-symmetric)
tensor category.  To construct such a category, one must divide
the tensor factors of an invariant space
$$\Inv(V_1 \tensor \ldots \tensor V_n)$$
into smaller tensor products which serve as the domain and target of
a space of morphisms.  This involves arbitrary choices that are extraneous
to most of the arguments in this paper, so we will not usually treat
spiders as categories.

In this paper, we will define certain spiders in terms of generators and
relations, and we will show that they are isomorphic to the representation
theories of rank two Lie algebras and the quantum deformations of these
representation theories.  These results generalize a well-known construction
for $A_1 = \sl(2)$ that first arose in a paper of Rumer, Teller, and
Weyl~\cite{RTW:gottingen1932},that was developed later by Temperley and
Lieb~\cite{TL:prsl1971}, and that was greatly developed recently by Jones,
Kauffman, Lickorish, Masbaum, and Vogel
\cite{Jones:invent1983,Kauffman:top1987,Lickorish:cmh1992,MV:pjm1994}. 
Moreover, Frenkel and Khovanov have recently established that the bases of
invariant spaces that arise when one constructs the $A_1$ spider by
generators and relations are dual to the canonical bases of Lusztig
\cite{Lusztig:jams1990}.  We conjecture that a similar phenomenon holds in
the rank 2 cases.

\subsection{Acknowledgements}

The author would like to thank Sarah Witherspoon, Bruce Westbury, Igor
Frenkel, and Mikhail Khovanov for their attention to this paper and
the work it presents.  Also, the \TeX\, macro package PSTricks~\cite{PSTricks}
was essential for typesetting the equations and figures.

\section{The $A_1$ spider and $\sl(2,\F)$}
\label{sa1}

Let $V\cong\F^2$ be the defining representation of $\sl(2,\F)$, where $\F$ is
some field of characteristic 0.  (The complex numbers $\C$ are a good choice
for $\F$.)  Choose a parameter $a\in\F$.  Since the vector space $V$ is a
self-dual representation, there exists a non-degenerate, invariant
contraction operation $\sigma:V\tensor V\to\F$.  For each $n$, choose $n$
points on the boundary of a disk in the plane, and let $B_n$, the set of
{\df basis webs}, be the set of crossingless matchings of the $n$ points.  For
example, if $n$ is odd, $B_n$ is empty, while $B_6$ has the following 5 elements:
$$
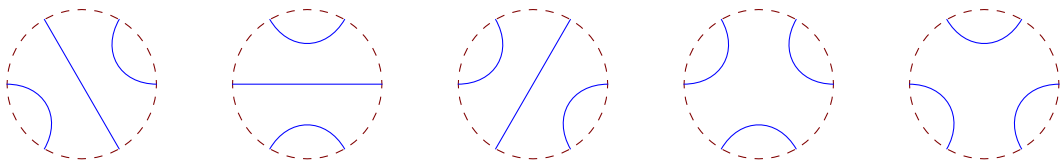
\pspicture(-1,-1)(1,1)
\pscircle[linecolor=darkred,linestyle=dashed](0,0){1}
\psbezier(1;0)(.5;0)(.5;60)(1;60)
\qline(1;120)(1;300)
\psbezier(1;180)(.5;180)(.5;240)(1;240)
\endpspicture
\hspace{1cm}
\pspicture(-1,-1)(1,1)
\pscircle[linecolor=darkred,linestyle=dashed](0,0){1}
\psbezier(1;60)(.5;60)(.5;120)(1;120)
\qline(1;180)(1;0)
\psbezier(1;240)(.5;240)(.5;300)(1;300)
\endpspicture
\hspace{1cm}
\pspicture(-1,-1)(1,1)
\pscircle[linecolor=darkred,linestyle=dashed](0,0){1}
\psbezier(1;120)(.5;120)(.5;180)(1;180)
\qline(1;240)(1;60)
\psbezier(1;300)(.5;300)(.5;0)(1;0)
\endpspicture
\hspace{1cm}
\pspicture(-1,-1)(1,1)
\pscircle[linecolor=darkred,linestyle=dashed](0,0){1}
\psbezier(1;0)(.5;0)(.5;60)(1;60)
\psbezier(1;120)(.5;120)(.5;180)(1;180)
\psbezier(1;240)(.5;240)(.5;300)(1;300)
\endpspicture
\hspace{1cm}
\pspicture(-1,-1)(1,1)
\pscircle[linecolor=darkred,linestyle=dashed](0,0){1}
\psbezier(1;60)(.5;60)(.5;120)(1;120)
\psbezier(1;180)(.5;180)(.5;240)(1;240)
\psbezier(1;300)(.5;300)(.5;0)(1;0)
\endpspicture
$$
By convention, $B_0$ has a single element, the empty disk. Let $W_n$, the {\df
web space}, be the vector space of formal linear combinations of elements of
$B_n$ with coefficients in $\F$.

\begin{theorem} There exist isomorphisms
$\phi_n:W_n\to\Inv(V^{\tensor n})$. \label{thcatalan}
\end{theorem}
Theorem~\ref{thcatalan} has few consequences in isolation; it only says that
web spaces have the same dimension as invariant spaces. By the Weyl character
formula, $\dim \Inv(V^{\tensor n})$ is the number of non-negative lattice
paths of length $n$ in one dimension, and there is a standard combinatorial
bijection between such paths and crossingless matchings.  (Idea:  Both
combinatorial sets are equivalent to balanced lists of parentheses of length
$n$, such as ``{\tt (())()}'' when $n = 6$.)

A {\df rotation} operation is the linear extension to a web space of
rotation of basis webs:
$$
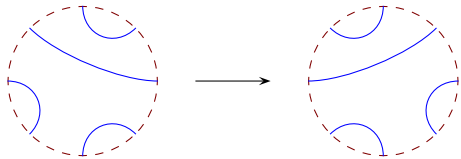
\pspicture[.5](-1,-1)(1,1)
\pscircle[linecolor=darkred,linestyle=dashed](0,0){1}
\psbezier(1;0)(.5;0)(.5;135)(1;135)
\psbezier(1;45)(.6;45)(.6;90)(1;90)
\psbezier(1;180)(.6;180)(.6;225)(1;225)
\psbezier(1;270)(.6;270)(.6;315)(1;315)
\endpspicture
\psgoesto
\pspicture[.5](-1,-1)(1,1)
\pscircle[linecolor=darkred,linestyle=dashed](0,0){1}
\psbezier(1;45)(.5;45)(.5;180)(1;180)
\psbezier(1;90)(.6;90)(.6;135)(1;135)
\psbezier(1;225)(.6;225)(.6;270)(1;270)
\psbezier(1;315)(.6;315)(.6;0)(1;0)
\endpspicture
$$
A {\df join} operation is a bilinear operation on a pair of web spaces. Given
two basis webs, their join is given by connecting their disks by a band:
$$
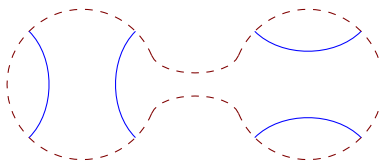
\pspicture(-1,-1)(4,1)
\psarc[linecolor=darkred,linestyle=dashed](0,0){1}{22.620}{-22.620}
\psbezier[linecolor=darkred,linestyle=dashed](!12 13 div 5 13 div)
(!12 13 div 5 39 div add 5 13 div 12 39 div sub)
(!3 12 13 div 5 39 div add sub 5 13 div 12 39 div sub)
(!3 12 13 div sub 5 13 div)
\psbezier[linecolor=darkred,linestyle=dashed](!12 13 div -5 13 div)
(!12 13 div 5 39 div add -5 13 div -12 39 div sub)
(!3 12 13 div 5 39 div add sub -5 13 div -12 39 div sub)
(!3 12 13 div sub -5 13 div)
\psarc[linecolor=darkred,linestyle=dashed](3,0){1}{202.620}{157.380}
\psbezier(1;45)(.5;45)(.5;-45)(1;-45)
\psbezier(1;135)(.5;135)(.5;-135)(1;-135)
\psbezier(!3 45 cos add 45 sin)(!3 45 cos 2 div add 45 sin 2 div)
(!3 135 cos 2 div add 135 sin 2 div)(!3 135 cos add 135 sin)
\psbezier(!3 -45 cos add -45 sin)(!3 -45 cos 2 div add -45 sin 2 div)
(!3 -135 cos 2 div add -135 sin 2 div)(!3 -135 cos add -135 sin)
\endpspicture
$$
A join operation extends bilinearly to arbitrary webs.  Finally, a {\df stitch}
operation is a linear transformation between web spaces.   Given a basis web,
its stitch at an adjacent pair of vertices is given by connecting the
vertices by an arc:
$$
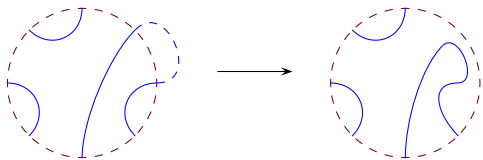
\pspicture[.5](-1,-1)(1.3,1.3)
\pscircle[linecolor=darkred,linestyle=dashed](0,0){1}
\psbezier(1;270)(.5;270)(.5;45)(1;45)
\psbezier(1;315)(.6;315)(.6;0)(1;0)
\psbezier(1;90)(.6;90)(.6;135)(1;135)
\psbezier(1;180)(.6;180)(.6;225)(1;225)
\psbezier[linestyle=dashed](1;0)(1.6;0)(1.6;45)(1;45)
\endpspicture
\psgoesto
\pspicture[.5](-1,-1)(1.3,1.3)
\pscircle[linecolor=darkred,linestyle=dashed](0,0){1}
\psbezier(1;270)(.5;270)(.3;45)(.7;45)
\psbezier(.7;45)(1;45)(1;0)(.7;0)
\psbezier(.7;0)(.3;0)(.6;315)(1;315)
\psbezier(1;90)(.6;90)(.6;135)(1;135)
\psbezier(1;180)(.6;180)(.6;225)(1;225)
\endpspicture
$$
If the result produces a closed loop, then strictly speaking, it is not
a matching of vertices on the boundary and therefore it is not a basis web.
In this case, erase the closed loop and replace it by a factor of $a$, the
element of $\F$ chosen at the beginning:
$$
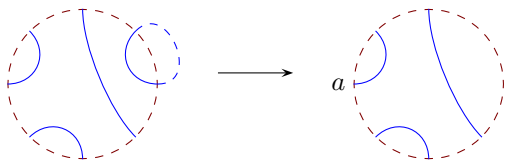
\pspicture[.5](-1,-1)(1.3,1.3)
\pscircle[linecolor=darkred,linestyle=dashed](0,0){1}
\psbezier(1;-270)(.5;-270)(.5;-45)(1;-45)
\psbezier(1;-315)(.6;-315)(.6;0)(1;0)
\psbezier(1;-90)(.6;-90)(.6;-135)(1;-135)
\psbezier(1;-180)(.6;-180)(.6;-225)(1;-225)
\psbezier[linestyle=dashed](1;0)(1.6;0)(1.6;45)(1;45)
\endpspicture
\psgoesto
\pspicture[.5](-1.3,-1)(1.3,1.3)
\rput(-1.2,0){$a$}
\pscircle[linecolor=darkred,linestyle=dashed](0,0){1}
\psbezier(1;-270)(.5;-270)(.5;-45)(1;-45)
\psbezier(1;-90)(.6;-90)(.6;-135)(1;-135)
\psbezier(1;-180)(.6;-180)(.6;-225)(1;-225)
\endpspicture
$$
The factor of $a$ is a
reminder that stitch is a linear operation whose value on a basis web might
be a non-basis web.  A stitch also extends linearly from a family of map $B_n
\to W_{n-2}$ to a map  $W_n \to W_{n-2}$.

The {\df combinatorial $A_1$ spider}, parameterized by $a$, is the list of web
spaces, together with all rotation, join, and stitch operations.

\begin{theorem}[Rumer,Teller,Weyl] If $a = -2$, then the  isomorphisms
$\phi_n:W_n \to \Inv(V^{\tensor n})$ can be uniquely chosen to send the
operations of join to tensor product, stitch to contraction, and rotation to
cyclic permutation of tensor factors composed with negation. \label{thrtw}
\end{theorem}

(Difference in sign between rotation and cyclic permutation of tensor
factors is explained in Section~\ref{sprecise}.)

The conversion of join to tensor product requires some explanation. Tensor
products depend on a linear ordering of tensor factors, but the vertices in a
web are only cyclically ordered.  To realize an explicit association between
webs and invariant tensors, it is necessary to refine the cyclic ordering to
a linear ordering. Then join or tensor product becomes a process of
concatenation:
$$
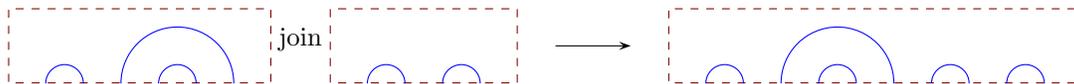
\pspicture[.5](0,0)(3.5,1)
\psframe[linecolor=darkred,linestyle=dashed](0,0)(3.5,1)
\psarc(.75,0){.25}{0}{180}
\psarc(2.25,0){.75}{0}{180}
\psarc(2.25,0){.25}{0}{180}
\endpspicture
\hspace{.1cm}\mbox{join}\hspace{.1cm}
\pspicture[.5](0,0)(2.5,1)
\psframe[linecolor=darkred,linestyle=dashed](0,0)(2.5,1)
\psarc(.75,0){.25}{0}{180}
\psarc(1.75,0){.25}{0}{180}
\endpspicture
\psgoesto
\pspicture[.5](0,0)(5.5,1)
\psframe[linecolor=darkred,linestyle=dashed](0,0)(5.5,1)
\psarc(.75,0){.25}{0}{180}
\psarc(2.25,0){.75}{0}{180}
\psarc(2.25,0){.25}{0}{180}
\psarc(3.75,0){.25}{0}{180}
\psarc(4.75,0){.25}{0}{180}
\endpspicture
$$
The general operation of join corresponds to tensor
product modified by cyclic permutations.

We explictly construct the $\phi_n$'s to demonstrate their uniqueness.
Assume some non-zero value  $r \in V \tensor V$ for $\phi_2$ of a line
segment:
$$
\phi_2\left(
\pspicture[.4](-.25,0)(1.75,.75)
\psframe[linecolor=darkred,linestyle=dashed](0,0)(1.5,.75)
\psarc(.75,0){.25}{0}{180}
\endpspicture\right) = r
$$
Any other basis web is obtained from a line segment by
repeated joins. Therefore its image under $\phi_n$ is, up to sign, a
combination of cyclic permutations and tensor products of $r$'s.  For
example, $\phi_4$ of the basis web
$$
\phi_4\left(
\pspicture[.4](-.25,0)(2.75,.75)
\psframe[linecolor=darkred,linestyle=dashed](0,0)(2.5,.75)
\psarc(.75,0){.25}{0}{180}
\psarc(1.75,0){.25}{0}{180}
\endpspicture\right) = r \tensor r,
$$
while
$$
\phi_6\left(
\pspicture[.4](-.25,0)(3.75,1)
\psframe[linecolor=darkred,linestyle=dashed](0,0)(3.5,1)
\psarc(.75,0){.25}{0}{180}
\psarc(2.25,0){.75}{0}{180}
\psarc(2.25,0){.25}{0}{180}
\endpspicture\right) = -r \tensor \rho_4(r \tensor r),$$
where $\rho_4$ is cyclic
permutation of four tensor factors. It remains only to determine $r$.  The
stitch of two arcs is one arc:
$$
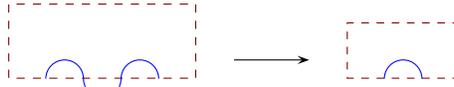
\pspicture[.5](0,-.25)(2.5,.75)
\psframe[linecolor=darkred,linestyle=dashed](0,0)(2.5,1)
\psarc(.75,0){.25}{0}{180}
\psarc(1.75,0){.25}{0}{180}
\psarc[linestyle=dashed](1.25,0){.25}{180}{0}
\endpspicture
\psgoesto
\pspicture[.5](0,-.25)(1.5,.75)
\psframe[linecolor=darkred,linestyle=dashed](0,0)(1.5,.75)
\psarc(.75,0){.25}{0}{180}
\endpspicture
$$
Applying $\phi_4$ and $\phi_2$ to the two sides yields
$$(I \tensor \sigma \tensor I)(r \tensor r) = r.$$
This equation says that $r$ is the inverse of
$\sigma$ in an indirect sense.  More explicitly, given a basis $\{e_1,e_2\}$
for $V$, one natural choice for $r$ and $\sigma$ is:
$$\sigma(e_2 \tensor e_1) = -\sigma(e_1 \tensor e_2) = 1$$
$$\sigma(e_1 \tensor e_1) = \sigma(e_2 \tensor e_2) = 0$$
$$r = e_1 \tensor e_2 - e_2 \tensor e_1$$
But in a more direct sense, $r$ is not the inverse of $\sigma$, because with
this choice of $r$ and $\sigma$, as with every other compatible pair of
choices, a closed loop, or $\sigma(r)$, is $-2$.

Checking that the $\phi_n$'s take spider operations to tensor operations
reduces to comparing definitions and simple calculations.  The least trivial
part of Theorem~\ref{thrtw} is the assertion that each $\phi_n$ is an
isomorphism.  Since the domain and target of $\phi_n$ have the same
dimension, it suffices to establish surjectivity. To prove surjectivity, we
use the isomorphism
$$\Inv(V^{\tensor 2n}) \cong \End(V^{\tensor n}).$$
By the Fundamental Theorem of Invariant Theory of Schur and Weyl,
the endomorphisms of $V^{\tensor n}$ are spanned by permutations
of tensor factors.  Such a permutation can be depicted by
a diagram of matched dots:
$$
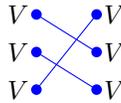
\pspicture(-.5,0)(1.5,1)
\qdisk( 0,0){2pt}\qdisk( 0,.5){2pt}\qdisk( 0,1){2pt}
\qdisk(.8,0){2pt}\qdisk(.8,.5){2pt}\qdisk(.8,1){2pt}
\qline(0,0)(.8,1)\qline(0,.5)(.8,0)\qline(0,1)(.8,.5)
\rput[r](-.1,0){$V$}\rput[r](-.1,.5){$V$}\rput[r](-.1,1){$V$}
\rput[l]( .9,0){$V$}\rput[l]( .9,.5){$V$}\rput[l]( .9,1){$V$}
\endpspicture
$$
It is a composition of many copies of the
switching map $\tau:x \tensor y \mapsto y \tensor x$
in $\End(V \tensor V)$.
A calculation demonstrates that $\tau$ lies in the image of $\phi_4$:
$$
\pspicture[.4](-.6,-.5)(.6,.5)
\pscircle[linecolor=darkred,linestyle=dashed](0,0){.5}
\qline(.5;45)(.5;225)\qline(.5;135)(.5;315)
\endpspicture
=
\pspicture[.4](-.6,-.5)(.6,.5)
\pscircle[linecolor=darkred,linestyle=dashed](0,0){.5}
\psbezier(.5;45)(.25;45)(.25;315)(.5;315)
\psbezier(.5;225)(.25;225)(.25;135)(.5;135)
\endpspicture
+
\pspicture[.4](-.6,-.5)(.6,.5)
\pscircle[linecolor=darkred,linestyle=dashed](0,0){.5}
\psbezier(.5;45)(.25;45)(.25;135)(.5;135)
\psbezier(.5;225)(.25;225)(.25;315)(.5;315)
\endpspicture
$$
It follows by multilinear expansion that every permutation is in the image
of $\phi_{2n}$.  Since the permutations span, $\phi_{2n}$
is surjective.

What if $a$ is not $-2$?  The specialization of the spider should describe
the representation theory of some object related to $\sl(2,\F)$.  The
quantum group $U_q(\sl(2,\F))$ was invented essentially for this purpose.
One parameterization of the $A_1$ spider is
$a = -q^{1/2} - q^{-1/2}$ over the field  $\C(q,q^{1/2},q^{1/3},\ldots)$.
(Note, however, that the spider is actually defined over the
ring $\C[q,q^{-1},q^{1/2},q^{1/3},\ldots]$, which allows us
to specialize to any value of $q \in \C^*$.)
Another standard notation is $a = -[2]$, where, by definition,
$$[n] = {q^{n/2} - q^{-n/2} \over q^{1/2}-q^{-1/2}}.$$

\subsection{Computations with the $A_1$ spider}

A basis web, or crossingless matching, is an equivalence class, modulo
boundary-fixing isotopies, of proper embeddings in the disk of a 1-manifold
with no circles.  Although proper embeddings of loop-free 1-manifolds suffice
to define the web space $W_n$, it is convenient to define any properly
embedded 1-manifold in the disk as a web.  Specifically, a 1-manifold with
$n$ closed loops denotes the web which is $a^n$ times the web obtained by
erasing all closed loops.

Given this meaning for all embedded 1-manifolds, basis webs
such as
$$
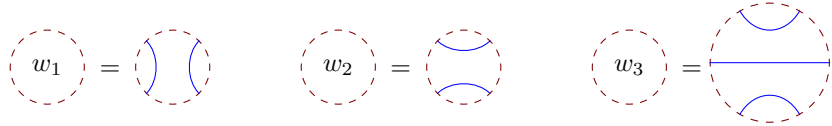
\pspicture[.4](-.6,-.5)(.6,.5)
\pscircle[linecolor=darkred,linestyle=dashed](0,0){.5}
\rput(0,0){$w_1$}
\endpspicture = \vertvert
\hspace{1cm}
\pspicture[.4](-.6,-.5)(.6,.5)
\pscircle[linecolor=darkred,linestyle=dashed](0,0){.5}
\rput(0,0){$w_2$}
\endpspicture = \horizhoriz
\hspace{1cm}
\pspicture[.4](-.6,-.5)(.6,.5)
\pscircle[linecolor=darkred,linestyle=dashed](0,0){.5}
\rput(0,0){$w_3$}
\endpspicture
=
\pspicture[.42](-.8,-1)(.8,1)
\pscircle[linecolor=darkred,linestyle=dashed](0,0){.8}
\psbezier(.8;60)(.4;60)(.4;120)(.8;120)
\qline(.8;180)(.8;0)
\psbezier(.8;240)(.4;240)(.4;300)(.8;300)
\endpspicture
$$
can be embedded in larger disks to form 1-manifolds which are
therefore other webs:
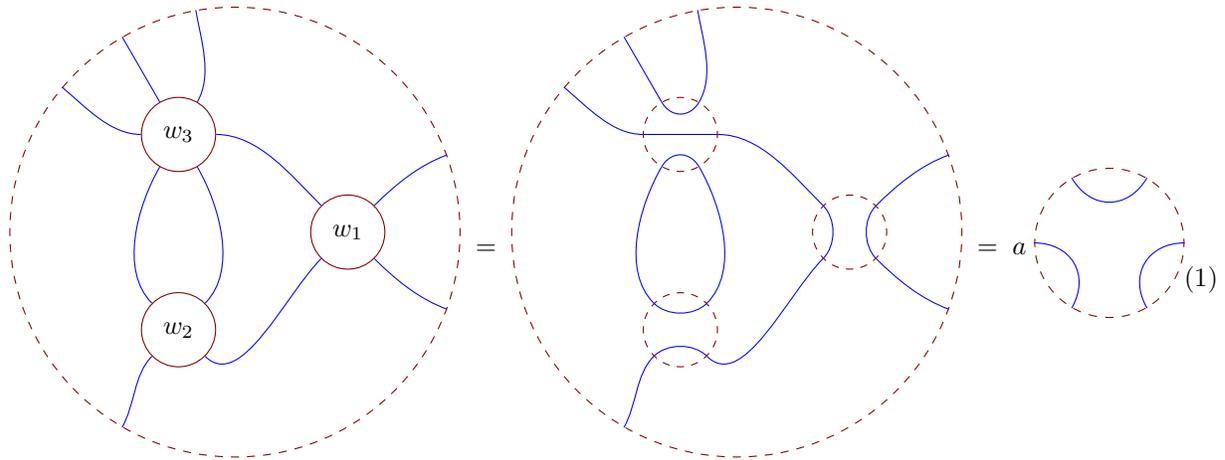
\begin{equation}
\pspicture[.45](-3.1,-3)(3.1,3)
\pscircle[linecolor=darkred,linestyle=dashed](0,0){3}
\pscircle[linecolor=darkred](1.5;0){.5}\rput(1.5;0){$w_1$}
\pscircle[linecolor=darkred](1.5;240){.5}\rput(1.5;240){$w_2$}
\pscircle[linecolor=darkred](1.5;120){.5}\rput(1.5;120){$w_3$}
\pccurve[nodesepA=.5,angleA=45,angleB=200](1.5;0)(3;20)
\pccurve[nodesepA=.5,angleA=-45,angleB=160](1.5;0)(3;-20)
\pccurve[nodesep=.5,angleA=135,angleB=0](1.5;0)(1.5;120)
\pccurve[nodesep=.5,angleA=225,angleB=-45](1.5;0)(1.5;240)
\pccurve[nodesep=.5,angleA=45,angleB=-60](1.5;240)(1.5;120)
\pccurve[nodesep=.5,angleA=135,angleB=240](1.5;240)(1.5;120)
\pccurve[nodesepA=.5,angleA=225,angleB=60](1.5;240)(3;240)
\pccurve[nodesepA=.5,angleA=60,angleB=280](1.5;120)(3;100)
\pccurve[nodesepA=.5,angleA=120,angleB=300](1.5;120)(3;120)
\pccurve[nodesepA=.5,angleA=180,angleB=320](1.5;120)(3;140)
\endpspicture
=
\pspicture[.45](-3.1,-3)(3.1,3)
\pscircle[linecolor=darkred,linestyle=dashed](0,0){3}
\rput(1.5;0){\vertvert}
\rput(1.5;240){\horizhoriz}
\rput(1.5;120){\pspicture[.42](-.5,-.5)(.5,.5)
\pscircle[linecolor=darkred,linestyle=dashed](0,0){.5}
\psbezier(.5;60)(.25;60)(.25;120)(.5;120)
\qline(.5;180)(.5;0)
\psbezier(.5;240)(.25;240)(.25;300)(.5;300)
\endpspicture}
\pccurve[nodesepA=.5,angleA=45,angleB=200](1.5;0)(3;20)
\pccurve[nodesepA=.5,angleA=-45,angleB=160](1.5;0)(3;-20)
\pccurve[nodesep=.5,angleA=135,angleB=0](1.5;0)(1.5;120)
\pccurve[nodesep=.5,angleA=225,angleB=-45](1.5;0)(1.5;240)
\pccurve[nodesep=.5,angleA=45,angleB=-60](1.5;240)(1.5;120)
\pccurve[nodesep=.5,angleA=135,angleB=240](1.5;240)(1.5;120)
\pccurve[nodesepA=.5,angleA=225,angleB=60](1.5;240)(3;240)
\pccurve[nodesepA=.5,angleA=60,angleB=280](1.5;120)(3;100)
\pccurve[nodesepA=.5,angleA=120,angleB=300](1.5;120)(3;120)
\pccurve[nodesepA=.5,angleA=180,angleB=320](1.5;120)(3;140)
\endpspicture
= \hspace{.1cm} a
\pspicture[.42](-1.1,-1)(1.1,1)
\pscircle[linecolor=darkred,linestyle=dashed](0,0){1}
\psbezier(1;60)(.5;60)(.5;120)(1;120)
\psbezier(1;180)(.5;180)(.5;240)(1;240)
\psbezier(1;300)(.5;300)(.5;0)(1;0)
\endpspicture
\label{ecomweb}
\end{equation}
If the $w_i$'s are arbitrary webs, then a diagram such as the one in
equation~(\ref{ecomweb}) denotes a multilinear expansion and is called
a compound web. For example, if
$$
\pspicture[.4](-.6,-.5)(.6,.5)
\pscircle[linecolor=darkred,linestyle=dashed](0,0){.5}
\rput(0,0){$w$}
\endpspicture
= \hspace{.1cm} 2 \vertvert + 3 \horizhoriz
$$
then
\begin{eqnarray*}
\pspicture[.42](-2.1,-1)(2.1,1)
\psellipse[linecolor=darkred,linestyle=dashed](0,0)(2,1)
\pscircle[linecolor=darkred](-1,0){.5}\rput(-1,0){$w$}
\pscircle[linecolor=darkred](1,0){.5}\rput(1,0){$w$}
\pccurve[nodesep=.5,angleA=45,angleB=135](-1,0)(1,0)
\pccurve[nodesep=.5,angleA=-45,angleB=-135](-1,0)(1,0)
\pccurve[nodesepA=.5,angleA=135,angleB=-45](-1,0)(-1.7,.527)
\pccurve[nodesepA=.5,angleA=225,angleB=45](-1,0)(-1.7,-.527)
\pccurve[nodesepA=.5,angleA=45,angleB=-135](1,0)(1.7,.527)
\pccurve[nodesepA=.5,angleA=-45,angleB=135](1,0)(1.7,-.527)
\endpspicture
& = &
4
\pspicture[.42](-2.1,-1)(2.1,1)
\psellipse[linecolor=darkred,linestyle=dashed](0,0)(2,1)
\rput(-1,0){\vertvert}
\rput(1,0){\vertvert}
\pccurve[nodesep=.5,angleA=45,angleB=135](-1,0)(1,0)
\pccurve[nodesep=.5,angleA=-45,angleB=-135](-1,0)(1,0)
\pccurve[nodesepA=.5,angleA=135,angleB=-45](-1,0)(-1.7,.527)
\pccurve[nodesepA=.5,angleA=225,angleB=45](-1,0)(-1.7,-.527)
\pccurve[nodesepA=.5,angleA=45,angleB=-135](1,0)(1.7,.527)
\pccurve[nodesepA=.5,angleA=-45,angleB=135](1,0)(1.7,-.527)
\endpspicture
+ 6
\pspicture[.42](-2.1,-1)(2.1,1)
\psellipse[linecolor=darkred,linestyle=dashed](0,0)(2,1)
\rput(-1,0){\vertvert}
\rput(1,0){\horizhoriz}
\pccurve[nodesep=.5,angleA=45,angleB=135](-1,0)(1,0)
\pccurve[nodesep=.5,angleA=-45,angleB=-135](-1,0)(1,0)
\pccurve[nodesepA=.5,angleA=135,angleB=-45](-1,0)(-1.7,.527)
\pccurve[nodesepA=.5,angleA=225,angleB=45](-1,0)(-1.7,-.527)
\pccurve[nodesepA=.5,angleA=45,angleB=-135](1,0)(1.7,.527)
\pccurve[nodesepA=.5,angleA=-45,angleB=135](1,0)(1.7,-.527)
\endpspicture \\
& + & 6
\pspicture[.42](-2.1,-1)(2.1,1)
\psellipse[linecolor=darkred,linestyle=dashed](0,0)(2,1)
\rput(-1,0){\horizhoriz}
\rput(1,0){\vertvert}
\pccurve[nodesep=.5,angleA=45,angleB=135](-1,0)(1,0)
\pccurve[nodesep=.5,angleA=-45,angleB=-135](-1,0)(1,0)
\pccurve[nodesepA=.5,angleA=135,angleB=-45](-1,0)(-1.7,.527)
\pccurve[nodesepA=.5,angleA=225,angleB=45](-1,0)(-1.7,-.527)
\pccurve[nodesepA=.5,angleA=45,angleB=-135](1,0)(1.7,.527)
\pccurve[nodesepA=.5,angleA=-45,angleB=135](1,0)(1.7,-.527)
\endpspicture
+ 9
\pspicture[.42](-2.1,-1)(2.1,1)
\psellipse[linecolor=darkred,linestyle=dashed](0,0)(2,1)
\rput(-1,0){\horizhoriz}
\rput(1,0){\horizhoriz}
\pccurve[nodesep=.5,angleA=45,angleB=135](-1,0)(1,0)
\pccurve[nodesep=.5,angleA=-45,angleB=-135](-1,0)(1,0)
\pccurve[nodesepA=.5,angleA=135,angleB=-45](-1,0)(-1.7,.527)
\pccurve[nodesepA=.5,angleA=225,angleB=45](-1,0)(-1.7,-.527)
\pccurve[nodesepA=.5,angleA=45,angleB=-135](1,0)(1.7,.527)
\pccurve[nodesepA=.5,angleA=-45,angleB=135](1,0)(1.7,-.527)
\endpspicture \\
& = & (21+4a) \vertvert + 9 \horizhoriz
\end{eqnarray*}
Another view of a compound web is that it is a sequence of joins and stitches
of the component webs.  It can be realized by many different such sequences,
but they all have the same final value.

An important class of examples of compound webs are those generated
by a web called a {\df crossing}:
$$
\pspicture[.4](-.6,-.5)(.6,.5)
\pscircle[linecolor=darkred,linestyle=dashed](0,0){.5}
\qline(.5;135)(.5;315)
\psline[border=.1](.5;45)(.5;225)
\endpspicture
=
- q^{1/4}\vertvert
- q^{-1/4}\horizhoriz
$$
Given this definition, it is easy to check the identities
$$
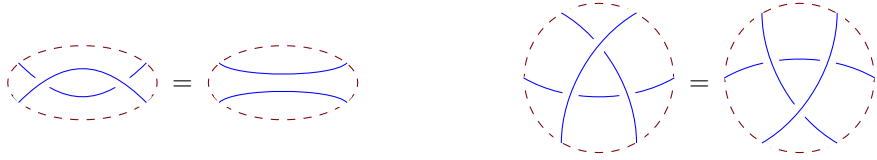
\pspicture[.4](-1.1,-.5)(1.1,.5)
\psellipse[linecolor=darkred,linestyle=dashed](0,0)(1,.5)
\pccurve[angleA=-45,angleB=-135,ncurv=1](-.85,.264)(.85,.264)
\pccurve[border=.1,angleA=45,angleB=135,ncurv=1](-.85,-.264)(.85,-.264)
\endpspicture
=
\pspicture[.4](-1.1,-.5)(1.1,.5)
\psellipse[linecolor=darkred,linestyle=dashed](0,0)(1,.5)
\pccurve[angleA=-45,angleB=-135,ncurv=.33](-.85,.264)(.85,.264)
\pccurve[border=.1,angleA=45,angleB=135,ncurv=.33](-.85,-.264)(.85,-.264)
\endpspicture
\hspace{2cm}
\pspicture[.42](-1.1,-1)(1.1,1)
\pscircle[linecolor=darkred,linestyle=dashed](0,0){1}
\pcarc[border=.1,arcangle=30](1;0)(1;180)
\pcarc[border=.1,arcangle=30](1;120)(1;300)
\pcarc[border=.1,arcangle=30](1;240)(1;60)
\endpspicture = 
\pspicture[.42](-1.1,-1)(1.1,1)
\pscircle[linecolor=darkred,linestyle=dashed](0,0){1}
\pcarc[border=.1,arcangle=30](1;180)(1;0)
\pcarc[border=.1,arcangle=30](1;300)(1;120)
\pcarc[border=.1,arcangle=30](1;60)(1;240)
\endpspicture
$$
These are known as the second and third {\df Reidemeister moves}. In
particular, a clever trick due to Kauffman \cite{Kauffman:top1987}, namely
replacing one of the crossings in the third move by its linear expansion,
reduces the third move to the second one. If we interpret a compound web made
from copies of a crossing as a tangle or link projection, the identities are
also known as the second and third Reidemeister moves.  More generally, a
projection of a tangle or link such as
$$
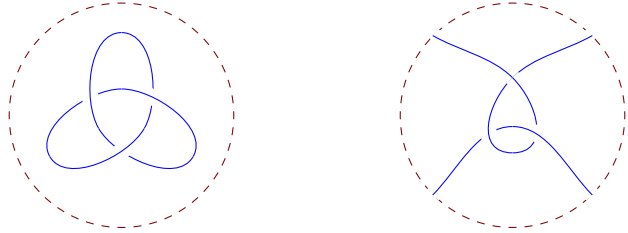
\pspicture(-1.6,-1.6)(1.6,1.6)
\pscircle[linecolor=darkred,linestyle=dashed](0,0){1.5}
\pccurve[border=.1,angleA=0,angleB=60](1.1;90)(.35;330)
\pccurve[border=.1,angleA=120,angleB=180](1.1;210)(.35;90)
\pccurve[border=.1,angleA=240,angleB=300](1.1;330)(.35;210)
\pccurve[border=.1,angleA=180,angleB=120](1.1;90)(.35;210)
\pccurve[border=.1,angleA=300,angleB=240](1.1;210)(.35;330)
\pccurve[border=.1,angleA=60,angleB=0](1.1;330)(.35;90)
\endpspicture
\hspace{2cm}
\pspicture(-1.6,-1.6)(1.6,1.6)
\pscircle[linecolor=darkred,linestyle=dashed](0,0){1.5}
\pccurve[border=.1,angleA=45,angleB=180](1.5;225)(0,-.15)
\pccurve[border=.1,angleA=210,angleB=45](1.5;45)(0,.5)
\pccurve[border=.1,angleA=225,angleB=180,ncurv=1](0,.5)(0,-.5)
\pccurve[border=.1,angleA=330,angleB=135](1.5;135)(0,.5)
\pccurve[border=.1,angleA=315,angleB=0,ncurv=1](0,.5)(0,-.5)
\pccurve[border=.1,angleA=135,angleB=0](1.5;315)(0,-.15)
\endpspicture
$$
evaluates to a vector in some web space.  The value of a link
projection lies in $W_0$, a 1-dimensional vector space with basis the empty
web.  Its single coefficient is a Laurent polynomial in $q^{1/4}$.  Given
invariance under the second and third Reidemeister moves, this function
on link projections is a regular isotopy invariant.
It is also {\df covariant} under full
isotopy; it gains a factor of $q^{\pm 3/4}$ under the first Reidemeister
move:
$$
\pspicture[.4](-.6,-.5)(.6,.5)
\pscircle[linecolor=darkred,linestyle=dashed](0,0){.5}
\psbezier(.5;120)(-.25,0)(.25,-.5)(.25,0)
\psbezier[border=.1](.5;240)(-.25,0)(.25,.5)(.25,0)
\endpspicture
=
q^{3/4}
\pspicture[.4](-.6,-.5)(.6,.5)
\pscircle[linecolor=darkred,linestyle=dashed](0,0){.5}
\psbezier(.5;240)(.25;240)(.25;120)(.5;120)
\endpspicture
$$
(Since $W_2$ is 1-dimensional, these two webs must be
proportional.) This polynomial is known as the {\df Kauffman bracket}, and up to
normalization it equals the Jones polynomial.

Another important type of compound web is a concatenation of two webs. Given a
web in $W_{a+b}$, divide its endpoints into a segment of $a$ points and a
segment of $b$ points.  Given another web in $W_{b+c}$ whose vertices are
divided into segments of length $b$ and $c$, their {\df concatenation}
consists of connecting $b$ adjacent pairs:
$$
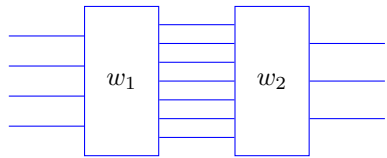
\pspicture(0,0)(5,2)
\qline(0,.4)(1,.4)\qline(0,.8)(1,.8)\qline(0,1.2)(1,1.2)\qline(0,1.6)(1,1.6)
\psframe(1,0)(2,2)\rput(1.5,1){$w_1$}
\qline(2, .25)(3, .25)\qline(2, .5 )(3, .5 )\qline(2, .75)(3, .75)
\qline(2,1   )(3,1   )\qline(2,1.25)(3,1.25)\qline(2,1.5 )(3,1.5 )
\qline(2,1.75)(3,1.75)
\psframe(3,0)(4,2)\rput(3.5,1){$w_2$}
\qline(4,.5)(5,.5)\qline(4,1)(5,1)\qline(4,1.5)(5,1.5)
\endpspicture
$$
With concatenation operations, the $A_1$ spider can be understood as a
category, isomorphic to a subcategory of the (quantum) representation
category of $\sl(2,\F)$, whose objects are segments of points and whose
arrows are webs.   For a fixed $n$, the endomorphisms of the object
consisting of $n$ points form an associative algebra, called the
{\df Temperley-Lieb} algebra.  As a unital algebra, it is generated by
$e_1,\ldots,e_{n-1}$, where $e_i$ is a basis web which is a pair of U-turns
at the $i$th and $i+1$st positions:
$$
e_i = \begin{array}{c}
\pspicture(0,0)(1,.2)\qline(0,0)(1,0)\qline(0,.2)(1,.2)\endpspicture \\
\vdots \\
\pspicture(0,0)(1,1)\qline(0,0)(1,0)\qline(0,.2)(1,.2)
\psarc(0,.6){.2}{-90}{90}\psarc(1,.6){.2}{90}{-90}\qline(0,1)(1,1)
\endpspicture \\
\vdots \\
\pspicture(0,0)(1,0)\qline(0,0)(1,0)\endpspicture\end{array}$$
It is easy to show that
$$\begin{array}{cc}
e_i e_j = e_j e_i, & i \ne j \pm 1 \\
e_i e_{i\pm1} e_i = e_i \\
e_i^2 = -[2] e_i \end{array}$$ 
is a complete set of relations for the Temperley-Lieb algebra.

\subsection{Other representations of $\sl(2,\F)$}

So far, we have only described $\Inv(V^{\tensor n})$ and not
$$\Inv(A_1 \tensor A_2 \tensor \ldots \tensor A_k)$$
for arbitrary irreducible representations $A_1,\ldots,A_k$ of $\sl(2,\F)$. 
The Lie algebra $\sl(2,\F)$ has an irreducible representation $V_n$ of
dimension $n+1$ for every non-negative $n$, and any finite-dimensional
irreducible is isomorphic to one of these.  The  representation $V_n$ can be
viewed as the $n$th symmetric power of  $V$, in which case the equivariant
contraction $\sigma$ induces an equivariant contraction
$$\sigma_n:V_n\tensor V_n \to \C.$$
It is determined by the rule that
$$\sigma_n(v^n \tensor w^n) = \sigma(v \tensor w)^n.$$
Given positive integers $n_1,\ldots,n_k$ with sum $n$, consider a circle with
$n$ distinguished points, partitioned into consecutive strings of points of
length $n_1, n_2, \ldots, n_k$. Each string of points is called an {\df
external clasp}.  The {\df clasped web space} $W(n_1,\ldots,n_k)$ is a vector
subspace of the web space $W_n$, defined as the span of those basis webs of
with no U-turns between two endpoints in the same clasp.  The set of such
basis webs is denoted $B(n_1,\ldots,n_k)$.  For example, the web
$$
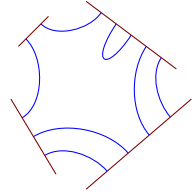
\pspicture(0,0)(2.4,2.5)
\psline[linecolor=darkred](0,1.2)(.6,.2)
\psline[linecolor=darkred](.1,1.9)(.5,2.3)
\psline[linecolor=darkred](1,2.5)(2.2,1.6)
\psline[linecolor=darkred](1,0)(2.4,1.2)
\pnode(1.28,0.24){a1}\pnode(1.56,0.48){a2}\pnode(1.84,0.72){a3}
\pnode(2.12,0.96){a4}
\pnode(2,1.75){b1}\pnode(1.8,1.9){b2}\pnode(1.6,2.05){b3}
\pnode(1.4,2.2){b4}\pnode(1.2,2.35){b5}
\pnode(.4,2.2){c1}\pnode(.2,2.0){c2}
\pnode(.15,.95){d1}\pnode(.3,.7){d2}\pnode(.45,.45){d3}
\nccurve[angleA=130.6,angleB=236.3]{a4}{b1}
\nccurve[angleA=130.6,angleB=236.3]{a3}{b2}
\nccurve[angleA=236.3,angleB=236.3,ncurv=5]{b3}{b4}
\nccurve[angleA=236.3,angleB=-45]{b5}{c1}
\nccurve[angleA=-45,angleB=31]{c2}{d1}
\nccurve[angleA=31,angleB=130.6]{d2}{a2}
\nccurve[angleA=31,angleB=130.6]{d3}{a1}
\endpspicture
$$
is not a basis web of $W(2,3,4,5)$, which instead has
basis
$$
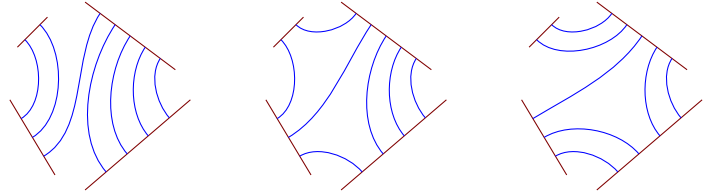
\pspicture(0,0)(2.4,2.5)
\psline[linecolor=darkred](0,1.2)(.6,.2)
\psline[linecolor=darkred](.1,1.9)(.5,2.3)
\psline[linecolor=darkred](1,2.5)(2.2,1.6)
\psline[linecolor=darkred](1,0)(2.4,1.2)
\pnode(1.28,0.24){a1}\pnode(1.56,0.48){a2}\pnode(1.84,0.72){a3}
\pnode(2.12,0.96){a4}
\pnode(2,1.75){b1}\pnode(1.8,1.9){b2}\pnode(1.6,2.05){b3}
\pnode(1.4,2.2){b4}\pnode(1.2,2.35){b5}
\pnode(.4,2.2){c1}\pnode(.2,2.0){c2}
\pnode(.15,.95){d1}\pnode(.3,.7){d2}\pnode(.45,.45){d3}
\nccurve[angleA=130.6,angleB=236.3]{a1}{b4}
\nccurve[angleA=130.6,angleB=236.3]{a2}{b3}
\nccurve[angleA=130.6,angleB=236.3]{a3}{b2}
\nccurve[angleA=130.6,angleB=236.3]{a4}{b1}
\nccurve[angleA=236.3,angleB=31]{b5}{d3}
\nccurve[angleA=-45,angleB=31]{c1}{d2}
\nccurve[angleA=-45,angleB=31]{c2}{d1}
\endpspicture
\hspace{1cm}
\pspicture(0,0)(2.4,2.5)
\psline[linecolor=darkred](0,1.2)(.6,.2)
\psline[linecolor=darkred](.1,1.9)(.5,2.3)
\psline[linecolor=darkred](1,2.5)(2.2,1.6)
\psline[linecolor=darkred](1,0)(2.4,1.2)
\pnode(1.28,0.24){a1}\pnode(1.56,0.48){a2}\pnode(1.84,0.72){a3}
\pnode(2.12,0.96){a4}
\pnode(2,1.75){b1}\pnode(1.8,1.9){b2}\pnode(1.6,2.05){b3}
\pnode(1.4,2.2){b4}\pnode(1.2,2.35){b5}
\pnode(.4,2.2){c1}\pnode(.2,2.0){c2}
\pnode(.15,.95){d1}\pnode(.3,.7){d2}\pnode(.45,.45){d3}
\nccurve[angleA=130.6,angleB=31]{a1}{d3}
\nccurve[angleA=130.6,angleB=236.3]{a2}{b3}
\nccurve[angleA=130.6,angleB=236.3]{a3}{b2}
\nccurve[angleA=130.6,angleB=236.3]{a4}{b1}
\nccurve[angleA=236.3,angleB=31]{b4}{d2}
\nccurve[angleA=236.3,angleB=-45]{b5}{c1}
\nccurve[angleA=-45,angleB=31]{c2}{d1}
\endpspicture
\hspace{1cm}
\pspicture(0,0)(2.4,2.5)
\psline[linecolor=darkred](0,1.2)(.6,.2)
\psline[linecolor=darkred](.1,1.9)(.5,2.3)
\psline[linecolor=darkred](1,2.5)(2.2,1.6)
\psline[linecolor=darkred](1,0)(2.4,1.2)
\pnode(1.28,0.24){a1}\pnode(1.56,0.48){a2}\pnode(1.84,0.72){a3}
\pnode(2.12,0.96){a4}
\pnode(2,1.75){b1}\pnode(1.8,1.9){b2}\pnode(1.6,2.05){b3}
\pnode(1.4,2.2){b4}\pnode(1.2,2.35){b5}
\pnode(.4,2.2){c1}\pnode(.2,2.0){c2}
\pnode(.15,.95){d1}\pnode(.3,.7){d2}\pnode(.45,.45){d3}
\nccurve[angleA=130.6,angleB=31]{a1}{d3}
\nccurve[angleA=130.6,angleB=31]{a2}{d2}
\nccurve[angleA=130.6,angleB=236.3]{a3}{b2}
\nccurve[angleA=130.6,angleB=236.3]{a4}{b1}
\nccurve[angleA=236.3,angleB=31]{b3}{d1}
\nccurve[angleA=236.3,angleB=-45]{b4}{c2}
\nccurve[angleA=236.3,angleB=-45]{b5}{c1}
\endpspicture
$$
Like a web with closed loops, a clasped web with U-turns has
a meaning as a web, but not a basis web.  If there are any U-turns, the
clasped web is defined to be the zero web.  Thus, $W(n_1,\ldots,n_k)$ is
also a quotient space of $W_n$.

\begin{theorem} Let $N = (n_1,\ldots,n_k)$ be a multi-index, with
$$V_N = V_{n_1} \tensor \ldots \tensor V_{n_k}.$$
When $q = 1$, there is a family of vector space isomorphisms
$$\phi_N:W(N) \to  \Inv(V_N)$$
that send join to tensor product, stitch to
contraction, and rotation to cyclic permutation of tensor factors up to sign.
\label{tha1isom}
\end{theorem}

Rotation and join of clasped webs are defined in the same way as for
unclasped webs, but stitch is more complicated.  Given a clasped web with two
adjacent external clasps of the same size, the stitch of the web is obtained
by identifying the two clasps and introducing an internal clasp, usually
depicted as a box:
$$
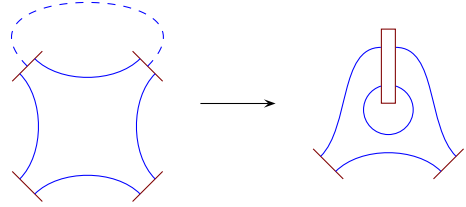
\pspicture[.5](0,0)(2,2.6)
\psline[linecolor=darkred](0,.4)(.4,0)\psline[linecolor=darkred](0,1.6)(.4,2)
\psline[linecolor=darkred](1.6,2)(2,1.6)\psline[linecolor=darkred](1.6,0)(2,.4)
\pccurve[angleA=45,angleB=135](.3,.1)(1.7,.1)
\pccurve[angleA=135,angleB=225](1.9,.3)(1.9,1.7)
\pccurve[angleA=225,angleB=315](1.7,1.9)(.3,1.9)
\pccurve[angleA=315,angleB=45](.1,1.7)(.1,.3)
\pccurve[linestyle=dashed,angleA=135,angleB=45,ncurv=2](.2,1.8)(1.8,1.8)
\endpspicture
\psgoesto
\pspicture[.5](0,0)(2,2)
\psline[linecolor=darkred](0,.4)(.4,0)
\psline[linecolor=darkred](1.6,0)(2,.4)
\psframe[linecolor=darkred](.9,1)(1.1,2)
\pccurve[angleA=45,angleB=135](.3,.1)(1.7,.1)
\pccurve[angleA=45,angleB=180](.1,.3)(.9,1.75)
\pccurve[angleA=135,angleB=0](1.9,.30)(1.1,1.75)
\pnode(1,1.25){temp}
\nccircle[angleA=180,nodesepA=.1]{temp}{.33}
\endpspicture
$$
An {\df internal clasp} of size $n$ is a particular web
in $W_{2n}$ that satisfies the axioms
$$
\pspicture[.4](-.3,0)(2.3,1)
\rput[r](-0.1,.5){$n$}\qline(0,.5)(.4,.5)
\psframe[linecolor=darkred](.4,0)(.6,1)
\qline(.6,.5)(1.4,.5)
\rput[b](1,.6){$n$}
\psframe[linecolor=darkred](1.4,0)(1.6,1)
\qline(1.6,.5)(2,.5)\rput[l](2.1,.5){$n$}
\endpspicture
=
\pspicture[.4](-.3,0)(1.3,1)
\rput[r](-0.1,.5){$n$}\qline(0,.5)(.4,.5)
\psframe[linecolor=darkred](.4,0)(.6,1)
\qline(.6,.5)(1,.5)\rput[l](1.1,.5){$n$}
\endpspicture
\hspace{2cm}
\pspicture[.4](-.3,0)(2.3,1)
\qline(0,.5)(.4,.5)\rput[r](-.1,.5){$n$}
\psframe[linecolor=darkred](.4,0)(.6,1)
\qline(.6,.833)(1,.833)\rput[l](1.1,.833){$k$}
\psarc(.6,.5){.167}{-90}{90}
\qline(.6,.167)(1,.167)\rput[l](1.1,.167){$n-k-2$}
\endpspicture = 0
\label{einternala1}
$$
Equation~(\ref{einternala1}) uses the convention that a strand labelled by $n$
denotes $n$ parallel strands.  The equation says that an internal clasp is an
idempotent of the $n$-strand Temperley-Lieb algebra and it annihilates all
basis webs other than the identity on the right. Therefore a clasp
concatenated with any web is proportional to a clasp. In particular, a clasp
is unique, if it exists, and it has the same annihilation property on the
left.  The Wenzl recursion formula \cite{Wenzl:canada1987} demonstrates that clasps do
exist,  at least for most values of $q$:
$$
\pspicture[.4](-.1,-.3)(1.1,1.3)
\rput[t](.5,-.1){$n$}\qline(.5,0)(.5,.4)
\psframe[linecolor=darkred](0,.4)(1,.6)
\qline(.5,.6)(.5,1)\rput[b](.5,1.1){$n$}
\endpspicture
=
\pspicture[.4](-.1,-.3)(1.3,1.3)
\rput[t](.5,-0.1){$n-1$}\qline(.5,0)(.5,.4)
\psframe[linecolor=darkred](0,.4)(1,.6)
\qline(.5,.6)(.5,1)\rput[b](.5,1.1){$n-1$}
\qline(1.3,0)(1.3,1)
\endpspicture
+
{[n-1] \over [n]}
\pspicture[.4](-.1,-.3)(1.35,2.3)
\rput[t](.5,-.1){$n-1$}
\qline(.5,0)(.5,.4)
\psframe[linecolor=darkred](0,.4)(1,.6)
\qline(.25,.6)(.25,1.4)\rput[l](.35,1){$n-2$}
\psframe[linecolor=darkred](0,1.4)(1,1.6)
\qline(.5,1.6)(.5,2)
\rput[b](.5,2.1){$n-1$}
\psarc(1,.6){.2}{0}{180}\qline(1.2,0)(1.2,.6)
\psarc(1,1.4){.2}{180}{0}\qline(1.2,1.4)(1.2,2)
\endpspicture
$$
Internal clasps are also called {\df magic weaving elements}
\cite{KL:recoupling}, {\df boxes} \cite{Lickorish:cmh1992}, and
{\df Jones-Wenzl idempotents} \cite{MV:pjm1994}.

An internal clasp is the concatenation of the unique basis web of
$W(n,1,1,1,\ldots,1)$ with itself:
$$
\pspicture[.5](-.3,0)(2.3,1)
\rput[r](-.1,.5){$n$}
\qline(0,.5)(.5,.5)\psline[linecolor=darkred](.5,0)(.5,1)
\psline[linecolor=darkred](1.5,0)(1.5,1)
\qline(1.5,.5)(2,.5)
\rput[l](2.1,.5){$n$}
\endpspicture
\psgoesto
\pspicture[.5](-.3,0)(1.3,1)
\rput[r](-0.1,.5){$n$}\qline(0,.5)(.4,.5)
\psframe[linecolor=darkred](.4,0)(.6,1)
\qline(.6,.5)(1,.5)
\rput[l](1.1,.5){$n$}
\endpspicture
$$
In terms of representations, this
concatenation is a composition of the form
$$V^{\tensor n} \to V_n \to V^{\tensor n}.$$
Since the composition is a non-zero idempotent, it is
the equivariant projection from $V^{\tensor n}$ to its highest-weight
irreducible summand.

\subsection{Isomorphism and equinumeration}

In this section, we prove Theorem~\ref{tha1isom}.  It is fairly easy to
construct each $\phi_N$ and to prove that is surjective.  If $n = n_1
+ \ldots + n_k$, there is a projection
$$\pi_N = \Inv(V^{\tensor n}) \to \Inv(V_N),$$
and let $\phi_N = \pi_N \circ \phi_n \circ i_N$, where $i_N$ is
the inclusion $W(N) \subset W_n$.  Suppose that the $i$th clasp of a web
$w$ has $k$ U-turns, and let $m = n - 2k$.  Then the
invariant $\pi_N \circ \phi_n(w)$ lies in the image of some map 
which has a tensor factor of the form
$V^{\tensor m} \to V_{n_i}$ , where $m < n_i$, and any such map must
be zero.

Thus, $\pi_N \circ \phi_n = \phi_N \circ j_N$, where $j_N$ is the projection
$W_n \to W(N)$.  From this, it is routine to show that the maps $phi_N$ take
spider operations to spider operations. Moreover, each $\phi_N$ is
surjective, because both $\pi_N$ and $\phi_n$ are surjective and $i_n$
complements the kernel of $\pi_N \circ \phi_n$.

We complete the proof that $\phi_N$ is an isomorphism by demonstrating that
its domain and target have the same dimension, thereby generalizing
Theorem~\ref{thcatalan} to the clasped case:

\begin{theorem} If $N$ is a multi-index, then
$$\dim W(N) = \dim \Inv(V_N).$$
\end{theorem}

\begin{proof}
The proof is by induction on $|N|$, the number of indicies in $N$.  The
relation is straightforward for $|N| \le 2$, so we first assume that $|N| = 3$.

The Clebsch-Gordan theorem states that
$$\dim \Inv(V_i \tensor V_j \tensor V_k) = 1$$
if each of $n$, $m$, and $l$ is less than or equal to the
sum of the other two and if the sum of all three is
even, and
$$\dim \Inv(V_i \tensor V_j \tensor V_k) = 0$$
otherwise.  On the other hand, $B(i,j,k)$ has at most one element:
$$
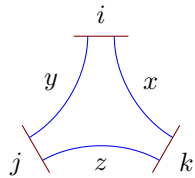
\pspicture(-1.5,-1.5)(1.5,1.5)
\pcline[linecolor=darkred](1.064;110)(1.064;70)\Aput{$i$}
\pcline[linecolor=darkred](1.064;230)(1.064;190)\Aput{$j$}
\pcline[linecolor=darkred](1.064;350)(1.064;310)\Aput{$k$}
\pccurve[angleA=270,angleB=30](1.015;100)(1.015;200)\Bput{$y$}
\pccurve[angleA=30,angleB=150](1.015;220)(1.015;320)\Bput{$z$}
\pccurve[angleA=150,angleB=270](1.015;340)(1.015;80)\Bput{$x$}
\endpspicture
$$
This web exists when there are non-negative integers $x$, $y$, and $z$ such
that $i = x+y$, $j = x+z$, and $k = y+z$.  These two conditions on $i$, $j$,
and $k$ are equivalent.

Now suppose that $|N|>3$ and express the multi-index $N$ as $JK$,
where $J$ and $K$ each have length at least 2.
Since the $V_n$'s constitute all finite-dimensional irreducible
representations, and since they are self-dual, there is a decomposition
$$\Inv(V_J \tensor V_K) \cong \bigoplus_\ell \Inv(V_J \tensor V_\ell) \tensor
\Inv(V_\ell \tensor V_K),$$
where $\ell$ is a single index rather than a multi-index.
It suffices to establish a
bijection
$$f:B(J,K) \longrightarrow \bigcup_\ell \left(B(J,\ell) \times B(\ell,K)\right),$$
where the union is disjoint.  The bijection $f$ is very easy to define:  A
basis web of $w \in W(JK)$ has a minimal cut path, where a {\df cut path} is
a path whose endpoints separate $J$ from $K$.  A cut path is {\df minimal} if
it crosses as few strands as possible:
$$
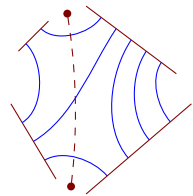
\pspicture(0,0)(2.4,2.5)
\psline[linecolor=darkred](0,1.2)(.6,.2)
\psline[linecolor=darkred](.1,1.9)(.5,2.3)
\psline[linecolor=darkred](1,2.5)(2.2,1.6)
\psline[linecolor=darkred](1,0)(2.4,1.2)
\pnode(1.28,0.24){a1}\pnode(1.56,0.48){a2}\pnode(1.84,0.72){a3}
\pnode(2.12,0.96){a4}
\pnode(2,1.75){b1}\pnode(1.8,1.9){b2}\pnode(1.6,2.05){b3}
\pnode(1.4,2.2){b4}\pnode(1.2,2.35){b5}
\pnode(.4,2.2){c1}\pnode(.2,2.0){c2}
\pnode(.15,.95){d1}\pnode(.3,.7){d2}\pnode(.45,.45){d3}
\nccurve[angleA=130.6,angleB=31]{a1}{d3}
\nccurve[angleA=130.6,angleB=236.3]{a2}{b3}
\nccurve[angleA=130.6,angleB=236.3]{a3}{b2}
\nccurve[angleA=130.6,angleB=236.3]{a4}{b1}
\nccurve[angleA=236.3,angleB=31]{b4}{d2}
\nccurve[angleA=236.3,angleB=-45]{b5}{c1}
\nccurve[angleA=-45,angleB=31]{c2}{d1}
\pcarc[arrows=*-*,linecolor=darkred,linestyle=dashed](.75,2.4)(.8,.1)
\endpspicture
$$
Let $w_1$ and $w_2$ be the two resulting webs. The cut path crosses some $i$
strands, and since it is minimal, there can be no U-turns among the $i$
strands in either $w_1$ or $w_2$.  Thus, the relation $f(w) = (w_1,w_2)$
defines the desired map $f$. It is routine to check that $f$ is both
injective and surjective.
\end{proof}

Note that minimal cut paths are also useful in practice for
generating the basis of a clasped web space.

\begin{exercise} List the elements of $B(2,2,2,2,2)$.
\end{exercise}

\section{Precise definition of a spider}
\label{sprecise}

In this section, we give a precise definition of a spider. The main way in
which the general notion of a spider is more complicated than the example is
that strands may be oriented and there may be more than one type of strand. 
Indeed, even in the $A_1$ case, we may consider $n$ parallel strands as
equivalent to one strand labelled with $n$, as is already suggested by the
notation.

A spider has a {\df strand set} $S$ which is a unital semigroup with unit
$\emptyset$. In most of the examples in the paper, $S$ is a free, non-abelian
semigroup.  The strand set has an anti-involution $*:S \to S$ called {\df
duality} or {\df orientation reversal}.  For each $s \in S$, there is a {\df
web space} $W(s)$, which may be just a set.  For each $a$, there is a
distinguished web $\beta_a \in W(aa^*)$, called a {\df bare strand}, and $1 =
\beta_\emptyset$ is called the {\df empty web}.  Finally, there are three
operations that exist for every $a$ and $b$:
\begin{itemize}
\item[1.] {\df Join} $\j_{a,b}:W(a) \times W(b) \to W(ab)$.
\item[2.] {\df Rotation} $\rho_{a,b}:W(ab) \to W(ba)$.
\item[3.] {\df Stitch} $\sigma_{a,b}:W(aa^*b) \to W(b)$.
\end{itemize}
The subscripts of the operations may be dropped when they are clear from
context. A spider may in addition be defined in a (symmetric tensor) category
other than the category of sets.  For example, an additive spider is one in
which web spaces are abelian groups, rotation and stitch are additive, and join
is additive, while in a linear spider, the web spaces are
vector spaces and the operations are linear or bilinear.

The three spider operations must satisfy the following axioms,
which are divided into groups according to which operations they
involve.  We define 
$\sigma_{b,a,c} = \rho_{c,b}\sigma_{a,cb}\rho_{b,aa^*c}$ for brevity
in the last three axioms.

\begin{itemize}
\item Rotation only:
\begin{itemize}
\item[1.] $\rho_{a,\emptyset} = \rho_{\emptyset,a} = \id$
\item[2.] $\rho_{a,bc}\rho_{c,ab}\rho_{b,ca} = \id$
\item[3.] $\rho_{a,a^*}(\beta_a) = \beta_{a^*}$
\end{itemize}
\item Join only:
\begin{itemize}
\item[4.] $(u \j v) \j w) = u \j (v \j w)$
\item[5.] $u \j 1 = u$
\end{itemize}
\item Join and rotation:
\begin{itemize}
\item[6.] $\rho_{a,b}(u \j_{ab,\emptyset} v)=\rho_{a,b}(u) \j_{ba,\emptyset}v$
\item[7.] $\beta_{ab} = \rho_{a^*,abb^*}(\beta_{a^*} \j \beta_b)$
\item[8.] $\rho_{a,b}(u \j_{a,b} v) = v \j_{b,a} u$
\end{itemize}
\item Rotation and stitch:
\begin{itemize}
\item[9.] $\sigma_{a,bd} \sigma_{aa^*b,c,d} = \sigma_{b,c,d}\sigma_{a,bcc^*d}$
\item[10.] $\sigma_{a,\emptyset} = \sigma_{a^*,\emptyset}\rho_{a,a^*}$.
\end{itemize}
\item Join, rotation, and stitch:
\begin{itemize}
\item[11.] $\sigma_{a,a^*,b}(\beta_a \j_{aa^*,ab} u) = u$
\item[12.] $\sigma_{b,a,c}(u \j_{ba,a^*c} v) =
\rho_{b,c}(\rho_{a^*,c}(v) \j_{ca^*,ab} \rho_{b,a}(u))$.
\end{itemize}
\end{itemize}

The formal spider axioms are embarrassingly complicated, but they can be
phrased in more natural (if less formal) terms.  In particular, the following
compound web principle is equivalent to axioms 1-7 and 10-12. To state the
principle, we first define a {\df pre-spider} to be an algebraic object
satisfying these 10 axioms, but not necessarily axioms 8 and 9.

Let $U$ be the free semigroup generated by a set $X$, and let $*$ be any
involution of $X$ extended to an anti-involution of $U$. Let $L_u$ be an
abstract set of labels for each $u \in U$ and let $L$ be the disjoint union.
Consider a graph $G$ in a disk such that the vertices of $G$ on the boundary
are univalent, such that each edge is labelled by an element of $x \in X$ and
oriented unless $x = x^*$, and such that one of the vertices at the boundary
is distinguished as first, and such that one edge of each internal vertex is
distinguished as first.  The edges incident to an internal vertex $v$ are
then linearly ordered going counterclockwise around $v$; let $x_i$ be the
label of the $i$th edge $e$ of $v$ if the edge is unoriented or oriented
outward, and let $x_i$ be the dual of the label of $e$ if it is oriented
inward.  Each $v$ should be labelled by an element of $L_{x_1\ldots x_n}$.
Then we define the web space $W(x_1\ldots x_n)$ as the set of all graphs $G$
with boundary labelled $x_1,\ldots,x_n$ (following the same convention of
taking the dual when an edge at the boundary is oriented inward), considered
up to isotopy and up to the modification of reversing an edge and
dualizing its label.

We define the free pre-spider ${\cal S}(X,L)$ by defining the spider
operations in the same way as for the combinatorial $A_1$ spider:  Join is
given by band-connected sum, rotation is given by changing which boundary
point is first, and stitch is given by connecting two adjacent boundary
points by an arc.  A tedious but straightforward computation demonstrates
that ${\cal S}(X,L)$ satisfies axioms 1-7 and 10-12.   Conversely, suppose
that an algebraic object $\cal S$ consists of a strand set $S$ and a
collection of web spaces $\{W(s)\}$ with operations of join, rotation, and
stitch.  Then the {\df compound web principle} stipulates that any
$*$-preserving map $f:X \to S$ and any set of maps $\{L_u \to W(f(u))\}$
extend to a morphism ${\cal S}(X,L) \to {\cal S}$ preserving join, rotation,
and stitch.  The compound web principle is equivalent to the statement that
$\cal S$ is a pre-spider.  Informally, in a pre-spider $\cal S$, if a
sequence of joins, rotations, and stitches are denoted by a planar diagram
whose connecting arcs are the stitches, the resulting web depends only on the
diagram and not on the order of the individual operations.  The remaining two
axioms for $\cal S$ can be understood as follows:  Axiom 8 says that if a
disconnected component of a compound web is moved from one face of the
remaining part of the web to another, it does not change the value of the
web.  Axiom 9 says that the value of a boundaryless compound web depends only
on its embedding in the sphere and not on its embedding in the disk.

Although the above definitions are not completely standard, a spider can also
be defined in category-theoretic terms.   A spider is a (small) {\df strict
monoidal category}, which is a category with an associative but not
necessarily commutative tensor product $\tensor$. Moreover, a spider must be
{\df pivotal} (or {\df rigid}), which means that there is a canonical
isomorphism $V^{**} \cong V$, and {\df spherical}, a condition which is
equivalent to axiom 9.  See Barrett and Westbury \cite{BW:spherical} for a
careful exposition of spherical categories.  Given any spherical category,
the strand set of the corresponding spider is precisely the set of objects of
the category, and its web space $W(V)$ is defined as $\Inv(V) = \Hom(T,V)$,
where $T$ is the trivial object, which serves as the identity of the tensor
product operation.  The spider operations are then defined in terms of
monoidal category operations in the same way as for the algebraic $A_1$ spider.

We review the spiders defined so far in terms of these definitions.
The strand set of the combinatorial and algebraic $A_1$ spiders
is a free semigroup with one generator; for $n \in \Z_{\ge 0}$, a
strand with strand type $n$ is synonymous with $n$ parallel strands.
The clasped combinatorial $A_1$ spider is a different object
and its strand set is a free semigroup with countably many generators,
namely the different clasps.

It is useful to treat an additive or linear spider or pre-spider as a ring,
and to consider the usual constructions with rings such as morphisms, ideals,
and quotients. An ideal $\cal I$ in a linear pre-spider $\cal S$ is a
collection of linear subspaces $I(s) \subset W(s)$ which are closed under
rotation and stitch and closed under join with an arbitrary web in $\cal S$. 
Clearly, if $\cal I$ is an ideal, the quotient spaces $W(s)/I(s)$ form a
pre-spider ${\cal S}/{\cal I}$. Given $X$ and $L$ as above, we can form the
free linear pre-spider $\overline{\cal S}(X,L)$ as the linear extension of
${\cal S}(X,L)$.  Many spiders (albeit only those whose strand set is a free
semigroup) can be defined in terms of generators and relations, meaning that
such a spider is a quotient of $\overline{\cal S}(X,L)$ by the ideal
generated by an arbitrary set of relators. Indeed, the combinatorial $A_1$
spider is defined in exactly this fashion.

Note that in an additive pre-spider, the web space $W(\emptyset)$
is a commutative ring, and the other web spaces 
become $W(\emptyset)$-modules under join.  Axiom 8 guarantees
that rotation and stitch are module endomorphisms, so that an
additive spider is automatically a $W(\emptyset)$-module spider,
and in particular it is a linear spider if $W(\emptyset)$ is a field.
See also Barrett and Westbury \cite{BW:spherical}.

The category of finite-dimensional representations of
any quantum Lie group $U_q(\frak g)$ with $\frak g$ a complex simple Lie
algebra is spherical, and therefore yields a spider.  Technically, this is
not a small category, meaning that the collection of all objects is too large
to be a set, but we can obtain an equivalent small category by taking a
single representative of each isomorphism class of finite-dimensional
representations.

In particular, self-dual representations correspond to self-dual or unoriented
strands, but therein lies a technicality and a potential sign error.  An
unoriented strand in a spider can only correspond to a symmetrically
self-dual representation, while many representations (for example the
representation $V_n$ of $\sl(2)$ for $n$ odd) are antisymmetrically
self-dual.  If there is to only one strand for each isomorphism class, each
self-dual representation must be defined as a $\Z/2$-graded vector space in
which the antisymmetric part has an odd grading.  Then any self-dual
representation of any $\frak g$ has a graded-symmetric invariant bilinear
form, and the representation can correspond to an unoriented strand. This is
why rotation in the combinatorial $A_1$ spider differs by a sign from
ordinary cyclic permutation of tensor factors in $\Inv(V_1^{\tensor n})$; it
exactly equals graded cyclic permutation. Among representations of rank two
Lie algebras, the representation $V(a\lambda_1+b\lambda_2)$ of $B_2$ also has
an odd grading when $a$ is odd.

Given a Lie algebra $\frak g$, the subcategory of irreducible representations
and their tensor products also yields a spider, which we will call the {\df
clasped algebraic $\frak g$ spider}.  The {\df unclasped algebraic $\frak g$
spider} is the subcategory whose objects are the fundamental irreducible
representations (those whose heighest weight is a simple weight) and their
tensor products.  In the rest of the paper, we will define combinatorial
$\frak g$ spiders when $\frak g$ has rank 2, and we will show that they are
isomorphic to their algebraic counterparts.

\section{The combinatorial rank 2 spiders}
\label{srank2}

For convenience, let $\C$ with $q \in \C$ or
$\C(q,q^{1/2},q^{1/3},\ldots)$ be the ground field. The unclasped
combinatorial rank 2 spiders describe the invariants of tensor
products of the two fundamental representations $V(\lambda_1)$ and
$V(\lambda_2)$ of each of the Lie algebras $A_2$, $B_2$, and $G_2$.  These
two representations are duals of each other for $A_2$ and are self-dual in
the other two cases.  The strand set for the combinatorial $A_2$
spider is defined as the free semigroup of strings of symbols ``$+$''
and ``$-$'', which correspond to $V(\lambda_1)$ and $V(\lambda_2)$,
respectively.  The dual of a sign string is given by reversing the string and
flipping the signs; for example,
$$(+ - - + +)^* = - - + + -.$$
In the $B_2$ and $G_2$ spiders, the strand set is the free semigroup
of strings of self-dual symbols``$1$'' and ``$2$'', so that
duality is just string reversal.

Given a sign string $s = s_1\ldots s_n$, define the {\df $A_2$ basis
web set} $B(s)$ to be the set of non-elliptic,
bipartite, trivalent graphs properly embedded in a disk with boundary points 
labelled $s_1,\ldots,s_n$ in counter-clockwise order.  By a {\df trivalent graph
properly embedded in a disk}, we mean a 1-dimensional subset of the disk
locally modelled by the following five allowed neighborhoods of a point in
the disk:
$$
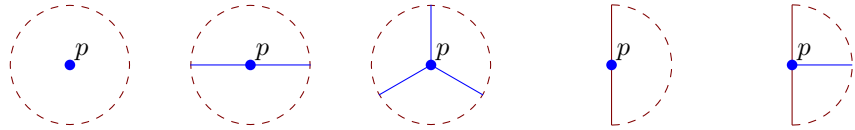
\pspicture(-.8,-.8)(.8,.8)
\pscircle[linecolor=darkred,linestyle=dashed](0,0){.8}
\qdisk(0,0){2pt}\rput[bl](2pt,2pt){$p$}
\endpspicture
\hspace{.8cm}
\pspicture(-.8,-.8)(.8,.8)
\pscircle[linecolor=darkred,linestyle=dashed](0,0){.8}
\qdisk(0,0){2pt}\rput[bl](2pt,2pt){$p$}
\qline(-.8,0)(.8,0)
\endpspicture
\hspace{.8cm}
\pspicture(-.8,-.8)(.8,.8)
\pscircle[linecolor=darkred,linestyle=dashed](0,0){.8}
\qdisk(0,0){2pt}\rput[bl](2pt,2pt){$p$}
\qline(0,0)(.8;210)\qline(0,0)(.8;330)\qline(0,0)(.8;90)
\endpspicture
\hspace{.8cm}
\pspicture(-.8,-.8)(.8,.8)
\psarc[linecolor=darkred,linestyle=dashed](0,0){.8}{270}{90}
\psline[linecolor=darkred](0,-.8)(0,.8)
\qdisk(0,0){2pt}\rput[bl](2pt,2pt){$p$}
\endpspicture
\hspace{.8cm}
\pspicture(-.8,-.8)(.8,.8)
\psarc[linecolor=darkred,linestyle=dashed](0,0){.8}{270}{90}
\psline[linecolor=darkred](0,-.8)(0,.8)
\qdisk(0,0){2pt}\rput[bl](2pt,2pt){$p$}
\qline(0,0)(.8,0)
\endpspicture
$$
The allowed neighborhoods might be called {\df empty disk}, {\df strand}, {\df
trivalent vertex}, {\df empty boundary}, and {\df endpoint}.  Such a graph is
{\df bipartite} if its endpoints are signed and its edges are oriented in
such a way that the in-degree at each vertex is either 0 or 3, and such that
edges point towards positive vertices and away from negative ones. Finally,
such a graph is {\df non-elliptic} if all internal faces have at least six
sides, where an internal face is a component of the complement of the graph
that does not touch the boundary of the disk. These graphs, henceforth called
{\df basis webs}, are considered up to isotopy relative to the boundary. For
example, the 6 elements of $B(+ + + - - -)$ are
$$
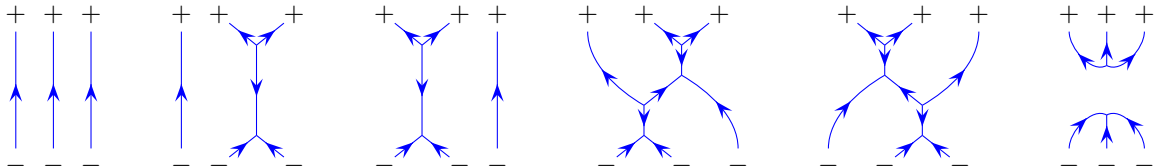
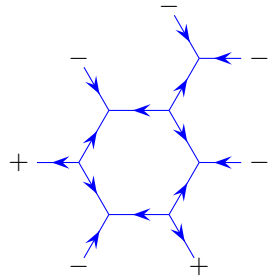
\pspicture(-.6,-1.1)(.6,1.1)
\rput(-.5,1){\rnode{p1}{$+$}}\rput(-.5,-1){\rnode{n1}{$-$}}
\rput(0,1){\rnode{p2}{$+$}}\rput(0,-1){\rnode{n2}{$-$}}
\rput(.5,1){\rnode{p3}{$+$}}\rput(.5,-1){\rnode{n3}{$-$}}
\ncline[nodesep=3pt]{n1}{p1}\middlearrow
\ncline[nodesep=3pt]{n2}{p2}\middlearrow
\ncline[nodesep=3pt]{n3}{p3}\middlearrow
\endpspicture
\hspace{1cm}
\pspicture(-.6,-1.1)(1.1,1.1)
\rput(-.5,1){\rnode{p1}{$+$}}\rput(-.5,-1){\rnode{n1}{$-$}}
\rput(0,1){\rnode{p2}{$+$}}\rput(0,-1){\rnode{n2}{$-$}}
\rput(1,1){\rnode{p3}{$+$}}\rput(1,-1){\rnode{n3}{$-$}}
\pnode(.5,.6){m1}\pnode(.5,-.6){m2}
\ncline[nodesep=3pt]{n1}{p1}\middlearrow\ncline{m1}{p2}\middlearrow
\ncline{m1}{p3}\middlearrow\ncline{m1}{m2}\middlearrow
\ncline{n2}{m2}\middlearrow\ncline{n3}{m2}\middlearrow
\endpspicture
\hspace{1cm}
\pspicture(-1.1,-1.1)(.6,1.1)
\rput(.5,1){\rnode{p1}{$+$}}\rput(.5,-1){\rnode{n1}{$-$}}
\rput(0,1){\rnode{p2}{$+$}}\rput(0,-1){\rnode{n2}{$-$}}
\rput(-1,1){\rnode{p3}{$+$}}\rput(-1,-1){\rnode{n3}{$-$}}
\pnode(-.5,.6){m1}\pnode(-.5,-.6){m2}
\ncline[nodesep=3pt]{n1}{p1}\middlearrow\ncline{m1}{p2}\middlearrow
\ncline{m1}{p3}\middlearrow\ncline{m1}{m2}\middlearrow
\ncline{n2}{m2}\middlearrow\ncline{n3}{m2}\middlearrow
\endpspicture
\hspace{1cm}
\pspicture(-1.1,-1.1)(1.1,1.1)
\rput(-1,1){\rnode{p1}{$+$}}
\rput(-.75,-1){\rnode{n1}{$-$}}
\rput(-.25,1){\rnode{p2}{$+$}}
\rput(.25,-1){\rnode{n2}{$-$}}
\rput(.75,1){\rnode{p3}{$+$}}
\rput(1,-1){\rnode{n3}{$-$}}
\pnode(-.25,-.6){m1}\pnode(-.25,-.2){m2}
\pnode(.25,.2){m3}\pnode(.25,.6){m4}
\ncline{n1}{m1}\middlearrow\ncline{n2}{m1}\middlearrow
\ncline{m2}{m1}\middlearrow\ncline{m2}{m3}\middlearrow
\ncline{m4}{m3}\middlearrow\ncline{m4}{p2}\middlearrow
\ncline{m4}{p3}\middlearrow
\nccurve[angleA=150,angleB=270,nodesepB=3pt]{m2}{p1}\middlearrow
\nccurve[angleB=330,angleA=90,nodesepA=3pt]{n3}{m3}\middlearrow
\endpspicture
\hspace{1cm}
\pspicture(-1.1,-1.1)(1.1,1.1)
\rput(1,1){\rnode{p1}{$+$}}
\rput(.75,-1){\rnode{n1}{$-$}}
\rput(.25,1){\rnode{p2}{$+$}}
\rput(-.25,-1){\rnode{n2}{$-$}}
\rput(-.75,1){\rnode{p3}{$+$}}
\rput(-1,-1){\rnode{n3}{$-$}}
\pnode(.25,-.6){m1}\pnode(.25,-.2){m2}
\pnode(-.25,.2){m3}\pnode(-.25,.6){m4}
\ncline{n1}{m1}\middlearrow\ncline{n2}{m1}\middlearrow
\ncline{m2}{m1}\middlearrow\ncline{m2}{m3}\middlearrow
\ncline{m4}{m3}\middlearrow\ncline{m4}{p2}\middlearrow
\ncline{m4}{p3}\middlearrow
\nccurve[angleA=30,angleB=270,nodesepB=3pt]{m2}{p1}\middlearrow
\nccurve[angleB=210,angleA=90,nodesepA=3pt]{n3}{m3}\middlearrow
\endpspicture
\hspace{1cm}
\pspicture(-.6,-1.1)(.6,1.1)
\rput(-.5,1){\rnode{p1}{$+$}}\rput(-.5,-1){\rnode{n1}{$-$}}
\rput(0,1){\rnode{p2}{$+$}}\rput(0,-1){\rnode{n2}{$-$}}
\rput(.5,1){\rnode{p3}{$+$}}\rput(.5,-1){\rnode{n3}{$-$}}
\pnode(0,-.33){m1}\pnode(0,.33){m2}
\nccurve[nodesepA=3pt,angleA=90,angleB=150]{n1}{m1}\middlearrow
\nccurve[nodesepA=3pt,angleA=90,angleB=270]{n2}{m1}\middlearrow
\nccurve[nodesepA=3pt,angleA=90,angleB=30]{n3}{m1}\middlearrow
\nccurve[nodesepB=3pt,angleB=270,angleA=210]{m2}{p1}\middlearrow
\nccurve[nodesepB=3pt,angleB=270,angleA=90]{m2}{p2}\middlearrow
\nccurve[nodesepB=3pt,angleB=270,angleA=330]{m2}{p3}\middlearrow
\endpspicture
$$
while
$$
\pspicture(-1.7,-1.7)(1.7,2.2)
\pnode(.8;0){a1}\pnode(.8;60){a2}\pnode(.8;120){a3}
\pnode(.8;180){a4}\pnode(.8;240){a5}\pnode(.8;300){a6}
\rput(1.6;0){\rnode{b1}{$-$}}\pnode(1.6;60){b2}
\rput(1.6;120){\rnode{b3}{$-$}}\rput(1.6;180){\rnode{b4}{$+$}}
\rput(1.6;240){\rnode{b5}{$-$}}\rput(1.6;300){\rnode{b6}{$+$}}
\rput(1.6,1.386){\rnode{c1}{$-$}}\rput(.4,2.078){\rnode{c2}{$-$}}
\ncline{a2}{a1}\middlearrow\ncline{a2}{a3}\middlearrow
\ncline{a2}{b2}\middlearrow\ncline{a4}{a3}\middlearrow
\ncline{a4}{a5}\middlearrow\ncline[nodesepB=3pt]{a4}{b4}\middlearrow
\ncline{a6}{a5}\middlearrow\ncline{a6}{a1}\middlearrow
\ncline{a6}{b6}\middlearrow\ncline[nodesepA=3pt]{b1}{a1}\middlearrow
\ncline{b3}{a3}\middlearrow\ncline{b5}{a5}\middlearrow
\ncline[nodesepA=3pt]{c1}{b2}\middlearrow\ncline{c2}{b2}\middlearrow
\endpspicture
$$
is an element of $B(+ - + - - - -)$.

The vector space of formal linear combinations of elements of $B(s_1\ldots
s_n)$ is the {\df $A_2$ web space} $W(s_1\ldots s_n)$. Partly elliptic,
bipartite, trivalent graphs in a disk, will denote webs also, although not
basis webs. Specifically, each type of elliptic face is defined as a linear
combination of basis webs according to the rules
\begin{eqnarray}
\pspicture[.4](-.6,-.5)(.6,.5)
\pscircle(0,0){.4}\psline[arrows=->,arrowscale=1.5](.1,.4)(.11,.4)
\endpspicture
& = & [3] \nonumber \\
\pspicture[.45](-1.5,-.8)(1.5,.8)
\rput(-1.2,0){\rnode{a1}{$-$}}
\pnode(-.4,0){a2}
\pnode(.4,0){a3}
\rput(1.2,0){\rnode{a4}{$+$}}
\ncline[nodesepA=3pt]{a1}{a2}\middlearrow
\nccurve[angleA=120,angleB=60,nodesep=.3pt]{a3}{a2}\middlearrow
\nccurve[angleA=-120,angleB=-60,nodesep=.3pt]{a3}{a2}\middlearrow
\ncline[nodesepB=3pt]{a3}{a4}\middlearrow
\endpspicture
& = & -[2]
\pspicture[.45](-.8,-.8)(.8,.8)
\rput(-.6,0){\rnode{a1}{$-$}}
\rput(.6,0){\rnode{a2}{$+$}}
\ncline[nodesep=3pt]{a1}{a2}\middlearrow
\endpspicture \nonumber \\
\pspicture[.42](-1.5,-1.4)(1.5,1.4)
\pnode(.4; 45){a1}\rput(1.2; 45){\rnode{b1}{$-$}}
\pnode(.4;135){a2}\rput(1.2;135){\rnode{b2}{$+$}}
\pnode(.4;225){a3}\rput(1.2;225){\rnode{b3}{$-$}}
\pnode(.4;315){a4}\rput(1.2;315){\rnode{b4}{$+$}}
\ncline[nodesepA=3pt]{b1}{a1}\middlearrow
\ncline[nodesepB=3pt]{a2}{b2}\middlearrow
\ncline[nodesepA=3pt]{b3}{a3}\middlearrow
\ncline[nodesepB=3pt]{a4}{b4}\middlearrow
\ncarc[arcangle= 15]{a2}{a1}\middlearrow
\ncarc[arcangle=-15]{a2}{a3}\middlearrow
\ncarc[arcangle= 15]{a4}{a3}\middlearrow
\ncarc[arcangle=-15]{a4}{a1}\middlearrow
\endpspicture
& = &
\pspicture[.4](-1,-.9)(1,.9)
\rput(.7; 45){\rnode{a1}{$-$}}\rput(.7;135){\rnode{a2}{$+$}}
\rput(.7;225){\rnode{a3}{$-$}}\rput(.7;315){\rnode{a4}{$+$}}
\ncarc[arcangle=-45]{a3}{a2}\middlearrow
\ncarc[arcangle=-45]{a1}{a4}\middlearrow
\endpspicture
+
\pspicture[.4](-1,-.9)(1,.9)
\rput(.7; 45){\rnode{a1}{$+$}}\rput(.7;135){\rnode{a2}{$-$}}
\rput(.7;225){\rnode{a3}{$+$}}\rput(.7;315){\rnode{a4}{$-$}}
\ncarc[arcangle=-45]{a3}{a2}\middlearrow
\ncarc[arcangle=-45]{a1}{a4}\middlearrow
\endpspicture
\label{fa2elliptic}
\end{eqnarray}
The value of a larger graph which contains an elliptic face is
inductively defined by the same rules.  For example,
\begin{eqnarray*}
\pspicture[.4](-1,-1)(2,1)
\pnode(.5; 30){a1}\pnode(.5; 90){a2}
\pnode(.5;150){a3}\pnode(.5;210){a4}
\pnode(.5;270){a5}\pnode(.5;330){a6}
\rput(1; 90){\rnode{b2}{$+$}}
\rput(1;150){\rnode{b3}{$-$}}
\rput(1;210){\rnode{b4}{$+$}}
\rput(1;270){\rnode{b5}{$-$}}
\pnode([nodesep=.5,angle=0]a1){c1}
\pnode([nodesep=.5,angle=0]a6){c2}
\rput([nodesep=.5,angle=45]c1){\rnode{d1}{$-$}}
\rput([nodesep=.5,angle=-45]c2){\rnode{d2}{$+$}}
\ncline[nodesepB=3pt]{a2}{b2}
\ncline{a3}{b3}\ncline{a4}{b4}
\ncline[nodesepB=3pt]{a5}{b5}
\ncline{a1}{a2}\ncline{a2}{a3}\ncline{a3}{a4}
\ncline{a4}{a5}\ncline{a5}{a6}\ncline{a6}{a1}
\ncline{a1}{c1}\ncline{a6}{c2}\ncline{c1}{c2}
\ncline{c1}{d1}\ncline{c2}{d2}
\endpspicture
& = &
\pspicture[.4](-1.3,-1.1)(1.3,1.1)
\pnode(.5; 90){a2}\pnode(.5;150){a3}
\pnode(.5;210){a4}\pnode(.5;270){a5}
\rput(1; 30){\rnode{b1}{$-$}}\rput(1; 90){\rnode{b2}{$+$}}
\rput(1;150){\rnode{b3}{$-$}}\rput(1;210){\rnode{b4}{$+$}}
\rput(1;270){\rnode{b5}{$-$}}\rput(1;330){\rnode{b6}{$+$}}
\ncline[nodesepB=3pt]{a2}{b2}
\ncline{a3}{b3}\ncline{a4}{b4}
\ncline[nodesepB=3pt]{a5}{b5}
\ncline{a2}{a5}\ncline{a2}{a3}\ncline{a3}{a4}\ncline{a4}{a5}
\nccurve[angleA=210,angleB=150,ncurv=1]{b1}{b6}
\endpspicture
+
\pspicture[.4](-1.3,-1.1)(1.3,1.1)
\pnode(.5; 90){a2}\pnode(.5;150){a3}
\pnode(.5;210){a4}\pnode(.5;270){a5}
\rput(1; 30){\rnode{b1}{$-$}}\rput(1; 90){\rnode{b2}{$+$}}
\rput(1;150){\rnode{b3}{$-$}}\rput(1;210){\rnode{b4}{$+$}}
\rput(1;270){\rnode{b5}{$-$}}\rput(1;330){\rnode{b6}{$+$}}
\ncline[nodesepB=3pt]{a2}{b2}
\ncline{a3}{b3}\ncline{a4}{b4}
\ncline[nodesepB=3pt]{a5}{b5}
\ncline{a2}{a3}\ncline{a3}{a4}\ncline{a4}{a5}
\nccurve[angleA=-30,angleB=210]{a2}{b1}
\nccurve[angleA=30,angleB=150]{a5}{b6}
\endpspicture
\\ & = &
\pspicture[.4](-1.3,-1.1)(1.3,1.1)
\rput(1; 30){\rnode{b1}{$-$}}\rput(1; 90){\rnode{b2}{$+$}}
\rput(1;150){\rnode{b3}{$-$}}\rput(1;210){\rnode{b4}{$+$}}
\rput(1;270){\rnode{b5}{$-$}}\rput(1;330){\rnode{b6}{$+$}}
\nccurve[angleA=270,nodesepA=3pt,angleB=330]{b2}{b3}
\nccurve[angleA= 30,angleB=90,nodesepB=3pt]{b4}{b5}
\nccurve[angleA=150,angleB=210]{b6}{b1}
\endpspicture
+
\pspicture[.4](-1.3,-1.1)(1.3,1.1)
\rput(1; 30){\rnode{b1}{$-$}}\rput(1; 90){\rnode{b2}{$+$}}
\rput(1;150){\rnode{b3}{$-$}}\rput(1;210){\rnode{b4}{$+$}}
\rput(1;270){\rnode{b5}{$-$}}\rput(1;330){\rnode{b6}{$+$}}
\nccurve[angleA=150,angleB=210]{b6}{b1}
\nccurve[angleA=330,angleB=30]{b3}{b4}
\ncline[nodesep=3pt]{b2}{b5}
\endpspicture
+
\pspicture[.4](-1.3,-1.1)(1.3,1.1)
\pnode(.5; 90){a2}\pnode(.5;150){a3}
\pnode(.5;210){a4}\pnode(.5;270){a5}
\rput(1; 30){\rnode{b1}{$-$}}\rput(1; 90){\rnode{b2}{$+$}}
\rput(1;150){\rnode{b3}{$-$}}\rput(1;210){\rnode{b4}{$+$}}
\rput(1;270){\rnode{b5}{$-$}}\rput(1;330){\rnode{b6}{$+$}}
\ncline[nodesepB=3pt]{a2}{b2}
\ncline{a3}{b3}\ncline{a4}{b4}
\ncline[nodesepB=3pt]{a5}{b5}
\ncline{a2}{a3}\ncline{a3}{a4}\ncline{a4}{a5}
\nccurve[angleA=-30,angleB=210]{a2}{b1}
\nccurve[angleA=30,angleB=150]{a5}{b6}
\endpspicture
\end{eqnarray*}
(Here and below, we may omit orientations of edges when they are clear from
context.) Rotation, join, and stitch operations are defined in the usual
graphical way:  For basis webs, rotation is rotation of the disk, join is
band-connected sum, and stitch is the operation of connecting two adjacent
boundary points by an arc.  Stitch is only defined when the two adjacent
boundary points have opposite sign, to preserve the orientation structure of
$A_2$ webs, and it may produce elliptic faces, which must be reduced to
obtain a linear combination of basis webs.  Compound webs, and in particular
concatenation, are also defined either by extension of the three basic
operations, or directly by the principle of reduction of elliptic faces.

Another way to phrase the definition of the $A_2$ spider is that it
is generated by the two webs
$$
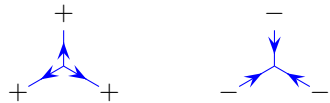
\pspicture(-.9,-.9)(.9,.9)
\rput(.7; 90){\rnode{b1}{$+$}}\rput(.7;210){\rnode{b2}{$+$}}
\rput(.7;330){\rnode{b3}{$+$}}
\pnode(0,0){a1}
\ncline[nodesepB=3pt]{a1}{b1}\middlearrow
\ncline{a1}{b2}\middlearrow
\ncline{a1}{b3}\middlearrow
\endpspicture
\hspace{1cm}
\pspicture(-.9,-.9)(.9,.9)
\rput(.7; 90){\rnode{b1}{$-$}}\rput(.7;210){\rnode{b2}{$-$}}
\rput(.7;330){\rnode{b3}{$-$}}
\pnode(0,0){a1}
\ncline[nodesepA=3pt]{b1}{a1}\middlearrow
\ncline{b2}{a1}\middlearrow
\ncline{b3}{a1}\middlearrow
\endpspicture
$$
with equations~(\ref{fa2elliptic}) as relations.  Similarly, the $A_1$ spider
is trivially generated with the sole relation that a closed loop yields a
factor of $a$.

The $B_2$ and $G_2$ spiders are also most conveniently defined by generators
and relations.  The $B_2$ spider is generated by the single web
$$
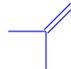
\pspicture(-.7,-.7)(.7,.7)
\qline(-.5,0)(0,0)\qline(0,-.5)(0,0)\psline[doubleline=true](0,0)(.35,.35)
\endpspicture
$$
with the relations
\begin{eqnarray}
\singleloop & = & -(q^2+q+q^{-1}+q^{-2}) \nonumber \\
\doubleloop & = & q^3+q+1+q^{-1}+q^{-3} \nonumber \\
\pspicture[.4](-.6,-.5)(.6,.5)
\psline[doubleline=true](-.5,0)(0,0)
\psbezier(0,0)(.7,.7)(.7,-.7)(0,0)
\endpspicture
& = & 0 \nonumber \\
\pspicture[.4](-.8,-.5)(.8,.5)
\psline[doubleline=true](-.7,0)(-.3,0)
\pcarc[arcangle=45](-.3,0)(.3,0)
\pcarc[arcangle=-45](-.3,0)(.3,0)
\psline[doubleline=true](.3,0)(.7,0)
\endpspicture
& = & -(q+2+q^{-1})\pspicture[.4](-.6,-.5)(.6,.5)
\psline[doubleline=true](-.5,0)(.5,0)
\endpspicture \nonumber \\
\pspicture[.4](-.9,-.9)(.9,.9)
\psline[doubleline=true](.4;90)(.8;90)
\psline[doubleline=true](.4;210)(.8;210)
\psline[doubleline=true](.4;330)(.8;330)
\pcarc[arcangle=-15](.4;90)(.4;210)
\pcarc[arcangle=-15](.4;210)(.4;330)
\pcarc[arcangle=-15](.4;330)(.4;90)
\endpspicture
& = & 0 \nonumber \\
\btwohvert & - & \btwohhoriz = \hh - \vv \label{fb2elliptic}
\end{eqnarray}
where a strand is denoted by a double edge if its strand type is ``2'' and by a
plain single edge if its strand type is ``1''.  The $G_2$ spider is generated
by the webs
$$
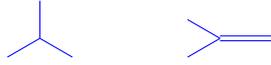
\pspicture(-.7,-.7)(.7,.7)
\qline(.5;90)(0,0)\qline(.5;210)(0,0)\qline(.5;330)(0,0)
\endpspicture
\hspace{1cm}
\pspicture(-.7,-.7)(.7,.7)
\qline(.5;150)(0,0)\qline(.5;210)(0,0)
\psline[doubleline=true](0,0)(.7,0)
\endpspicture
$$
with the relations
\begin{eqnarray}
\singleloop & = & q^5 + q^4 + q + 1 + q^{-1} + q^{-4} + q^{-5} \nonumber \\
\doubleloop & = & q^9 + q^6 + q^5 + q^4 + q^3 + q + 2 + q^{-1}
+ q^{-3} + q^{-4} + q^{-5} + q^{-6} + q^{-9} \nonumber \\
\pspicture[.4](-.6,-.5)(.6,.5)
\qline(-.5,0)(0,0)
\psbezier(0,0)(1;60)(1;300)(0,0)
\endpspicture
& = & 0 \nonumber \\
\pspicture[.4](-.9,-.5)(.9,.5)
\qline(-.8,0)(-.4,0)\qline(.4,0)(.8,0)
\pcarc[arcangle=60](-.4,0)(.4,0)
\pcarc[arcangle=-60](-.4,0)(.4,0)
\endpspicture
& = & -(q^3+q^2+q+q^{-1}+q^{-2}+q^{-3})\pspicture[.4](-.6,-.5)(.6,.5)
\qline(-.5,0)(.5,0) \endpspicture \nonumber \\
\pspicture[.4](-.9,-.9)(.9,.9)
\psline(.4;90)(.8;90)\psline(.4;210)(.8;210)
\psline(.4;330)(.8;330)
\pcarc[arcangle=-30](.4;90)(.4;210)
\pcarc[arcangle=-30](.4;210)(.4;330)
\pcarc[arcangle=-30](.4;330)(.4;90)
\endpspicture
& = & (q^2 + 1 + q^{-2})
\pspicture[.4](-.9,-.9)(.9,.9)
\qline(.7;90)(0,0)\qline(.7;210)(0,0)\qline(.7;330)(0,0)
\endpspicture \nonumber \\
\pspicture[.4](-.9,-.9)(.9,.9)
\psline(.4;45)(.8;45)\psline(.4;135)(.8;135)
\psline(.4;225)(.8;225)\psline(.4;315)(.8;315)
\pcarc[arcangle=-15](.4;45)(.4;135)\pcarc[arcangle=-15](.4;135)(.4;225)
\pcarc[arcangle=-15](.4;225)(.4;315)\pcarc[arcangle=-15](.4;315)(.4;45)
\endpspicture
& = & - (q + q^{-1})\left(\gtwohvert + \gtwohhoriz\right)
+ (q + 1 + q^{-1})\left(\hh + \vv\right) \nonumber \\
\pspicture[.4](-1,-1)(1,1)\pentanode
\ncline{a1}{b1}\ncline{a2}{b2}\ncline{a3}{b3}\ncline{a4}{b4}
\ncline{a5}{b5}
\ncarc[arcangle=6]{a1}{a5}\ncarc[arcangle=6]{a2}{a1}
\ncarc[arcangle=6]{a3}{a2}\ncarc[arcangle=6]{a4}{a3}
\ncarc[arcangle=6]{a5}{a4}
\endpspicture
& = & -\left(
\pspicture[.4](-1,-1)(1,1)\pentanode
\ncline{a5}{b5}\ncline{a1}{b1}\ncline{a2}{b2}
\ncarc[arcangle=6]{a2}{a1}\ncarc[arcangle=6]{a1}{a5}
\nccurve[angleA=312,angleB= 54]{a2}{b3}
\nccurve[angleA=258,angleB=126]{a5}{b4}
\endpspicture +
\pspicture[.4](-1,-1)(1,1)\pentanode
\ncline{a1}{b1}\ncline{a2}{b2}\ncline{a3}{b3}
\ncarc[arcangle=6]{a3}{a2}\ncarc[arcangle=6]{a2}{a1}
\nccurve[angleA= 24,angleB=126]{a3}{b4}
\nccurve[angleA=330,angleB=198]{a1}{b5}
\endpspicture +
\pspicture[.4](-1,-1)(1,1)\pentanode
\ncline{a2}{b2}\ncline{a3}{b3}\ncline{a4}{b4}
\ncarc[arcangle=6]{a4}{a3}\ncarc[arcangle=6]{a3}{a2}
\nccurve[angleA= 96,angleB=198]{a4}{b5}
\nccurve[angleA= 42,angleB=270]{a2}{b1}
\endpspicture +
\pspicture[.4](-1,-1)(1,1)\pentanode
\ncline{a3}{b3}\ncline{a4}{b4}\ncline{a5}{b5}
\ncarc[arcangle=6]{a5}{a4}\ncarc[arcangle=6]{a4}{a3}
\nccurve[angleA=168,angleB=270]{a5}{b1}
\nccurve[angleA=114,angleB=342]{a3}{b2}
\endpspicture +
\pspicture[.4](-1,-1)(1,1)\pentanode
\ncline{a4}{b4}\ncline{a5}{b5}\ncline{a1}{b1}
\ncarc[arcangle=6]{a1}{a5}\ncarc[arcangle=6]{a5}{a4}
\nccurve[angleA=240,angleB=342]{a1}{b2}
\nccurve[angleA=186,angleB= 54]{a4}{b3}
\endpspicture\right) \nonumber \\
& & + \left(
\pspicture[.4](-1,-1)(1,1)\pentanode\ncline{b1}{a1}
\nccurve[angleA=210,angleB=342]{a1}{b2}
\nccurve[angleA= 54,angleB=126]{b3}{b4}
\nccurve[angleA=330,angleB=198]{a1}{b5}
\endpspicture +
\pspicture[.4](-1,-1)(1,1)\pentanode\ncline{b2}{a2}
\nccurve[angleA=282,angleB= 54]{a2}{b3}
\nccurve[angleA=126,angleB=198]{b4}{b5}
\nccurve[angleA= 42,angleB=270]{a2}{b1}
\endpspicture +
\pspicture[.4](-1,-1)(1,1)\pentanode\ncline{b3}{a3}
\nccurve[angleA=354,angleB=126]{a3}{b4}
\nccurve[angleA=198,angleB=270]{b5}{b1}
\nccurve[angleA=114,angleB=342]{a3}{b2}
\endpspicture +
\pspicture[.4](-1,-1)(1,1)\pentanode\ncline{b4}{a4}
\nccurve[angleA= 66,angleB=198]{a4}{b5}
\nccurve[angleA=270,angleB=342]{b1}{b2}
\nccurve[angleA=186,angleB= 54]{a4}{b3}
\endpspicture +
\pspicture[.4](-1,-1)(1,1)\pentanode\ncline{b5}{a5}
\nccurve[angleA=138,angleB=270]{a5}{b1}
\nccurve[angleA=342,angleB= 54]{b2}{b3}
\nccurve[angleA=258,angleB=126]{a5}{b4}
\endpspicture\right) \nonumber \\
\pspicture[.4](-.9,-.5)(.9,.5)
\pnode(-.25,0){a1}\pnode(.25,0){a2}
\qline([nodesep=.5,angle=150]a1)(a1)
\qline([nodesep=.5,angle=210]a1)(a1)
\qline([nodesep=.5,angle=30]a2)(a2)
\qline([nodesep=.5,angle=330]a2)(a2)
\psline[doubleline=true](a1)(a2)
\endpspicture
& = & \hh - \gtwohvert - {1 \over q^2 - 1 + q^{-2}}\vv 
+ {1 \over q + 1 + q^{-1}}\gtwohhoriz \label{fg2elliptic}
\end{eqnarray}
To form a basis, we first note in the $G_2$ case that the last relation allows
the elimination of all internal double edges, leaving only those with at
least one vertex at the boundary.  In the $B_2$ spider, we define a tetravalent
vertex to achieve the same end:
$$
\btwox = \btwohhoriz - \vv
$$
The tetravalent vertex then satisfies the relations
\begin{eqnarray*}
\pspicture[.4](-1.1,-.5)(1.1,.5)
\pccurve[angleA=-45,angleB=-135,ncurv=1](-.85,.264)(.85,.264)
\pccurve[angleA=45,angleB=135,ncurv=1](-.85,-.264)(.85,-.264)
\endpspicture
& = & - (q+2+q^{-1})\btwox - (q^2 + 2q + 2 + 2q^{-1} + q^{-2})\vv \\
\pspicture[.42](-1.1,-1)(1.1,1)
\pcarc[arcangle=30](1;0)(1;180)
\pcarc[arcangle=30](1;120)(1;300)
\pcarc[arcangle=30](1;240)(1;60)
\endpspicture
& = &
(q + 2 + q^{-1})\left(
\pspicture[.42](-1.1,-1)(1.1,1)
\pccurve[angleA=240,angleB=300](1; 60)(1;120)
\pccurve[angleA=  0,angleB=120](1;180)(1;300)
\pccurve[angleA= 60,angleB=180](1;240)(1;  0)
\endpspicture +
\pspicture[.42](-1.1,-1)(1.1,1)
\pccurve[angleA=  0,angleB= 60](1;180)(1;240)
\pccurve[angleA=120,angleB=240](1;300)(1; 60)
\pccurve[angleA=180,angleB=300](1;  0)(1;120)
\endpspicture +
\pspicture[.42](-1.1,-1)(1.1,1)
\pccurve[angleA=120,angleB=180](1;300)(1;  0)
\pccurve[angleA=240,angleB=  0](1; 60)(1;180)
\pccurve[angleA=300,angleB= 60](1;120)(1;240)
\endpspicture
\right) \\ 
& & + (q^2 + 4q + 6 + 4q^{-1} + q^{-2})
\pspicture[.42](-1.1,-1)(1.1,1)
\psbezier(1;60)(.5;60)(.5;120)(1;120)
\psbezier(1;180)(.5;180)(.5;240)(1;240)
\psbezier(1;300)(.5;300)(.5;0)(1;0)
\endpspicture
\end{eqnarray*}
In both spiders, we can say that non-elliptic webs with no internal double
edges are a basis, provided that we define formal angles of
$$
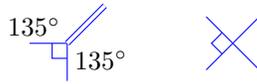
\pspicture(-.6,-.5)(.6,.5)
\psline(-.2,0)(-.2,-.2)(0,-.2)
\qline(-.5,0)(0,0)\qline(0,-.5)(0,0)
\psline[doubleline=true](0,0)(.5,.5)
\rput[rb](-.1,.1){$\scriptsize 135^\circ$}
\rput[lt](.1,-.1){$\scriptsize 135^\circ$}
\endpspicture
\hspace{1cm}
\pspicture(-.6,-.5)(.6,.5)
\psline(.2;135)(.2828;180)(.2;225)
\qline(.5;45)(.5;225)\qline(.5;135)(.5;315)
\endpspicture
$$
for the vertices in the $B_2$ spider and angles of
$$
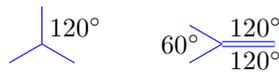
\pspicture(-.7,-.7)(.7,.7)
\qline(.5;90)(0,0)\qline(.5;210)(0,0)\qline(.5;330)(0,0)
\rput[lb](.1,.1){$\scriptsize 120^\circ$}
\endpspicture
\hspace{1cm}
\pspicture(-.7,-.7)(.7,.7)
\qline(.5;150)(0,0)\qline(.5;210)(0,0)
\rput[r](-.3,0){$\scriptsize 60^\circ$}
\rput[lb](.1,.1){$\scriptsize 120^\circ$}
\rput[lt](.1,-.1){$\scriptsize 120^\circ$}
\psline[doubleline=true](0,0)(.7,0)
\endpspicture
$$
in the $G_2$ spider, and we declare that a face is elliptic if its total
exterior angle is less than 360 degrees.  For example, in the $G_2$
spider, the pentagon
$$
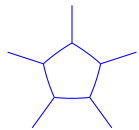
\pspicture(-1,-1)(1,1)
\pnode(.4;90){a1}\pnode(.4;162){a2}\pnode(.4;234){a3}
\pnode(.4;306){a4}\pnode(.4;18){a5}
\pnode(.9;90){b1}\pnode(.9;162){b2}\pnode(.9;234){b3}
\pnode(.9;306){b4}\pnode(.9;18){b5}
\ncline{a1}{b1}\ncline{a2}{b2}\ncline{a3}{b3}\ncline{a4}{b4}
\ncline{a5}{b5}
\ncarc[arcangle=6]{a1}{a5}\ncarc[arcangle=6]{a2}{a1}
\ncarc[arcangle=6]{a3}{a2}\ncarc[arcangle=6]{a4}{a3}
\ncarc[arcangle=6]{a5}{a4}
\endpspicture
$$
is elliptic, but the pentagon
$$
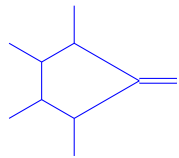
\pspicture(-1,-1)(1.5,1)
\pspolygon(.866;0)(.5;90)(.5;150)(.5;210)(.5;270)
\psline(.5;90)(1;90)\psline(.5;150)(1;150)
\psline(.5;210)(1;210)\psline(.5;270)(1;270)
\psline[doubleline=true](.866;0)(1.366;0)
\endpspicture
$$
is flat rather than elliptic.

If we understand the rank two spiders in terms of generators and relations,
then it is clear that non-elliptic webs {\df linearly span all webs}, but it
is not obvious that there are linearly independent. Alternatively, if we
understand the rank two spiders in terms of the non-elliptic bases, then it
is not clear that the elliptic reduction equations are consistent.  Put a
third way, do two different reductions of a partly elliptic web always give
the same linear combination of non-elliptic webs?  The
author~\cite{Kuperberg:g2}, and independently Jaeger \cite{Jaeger:confluent},
established that the given coefficients are, up to trivial normalization, the
only choices for which the equations are consistent.

Rank 2 spiders also admit crossings that lead to link invariants,
but the link invariants will not be discussed in this paper.  There are
two types of crossings in the $A_2$ spider:
\begin{eqnarray*}
\pspicture[.4](-1,-.9)(1,.9)
\rput(.7;45){\rnode{a1}{$+$}}
\rput(.7;135){\rnode{a2}{$+$}}
\rput(.7;225){\rnode{a3}{$-$}}
\rput(.7;315){\rnode{a4}{$-$}}
\ncline{a2}{a4}\ncline[border=.1]{a3}{a1}
\endpspicture
& = & q^{1/6}
\pspicture[.4](-1,-.9)(1,.9)
\pnode(0,.35){a1}\pnode(0,-.35){a2}
\rput([nodesep=.7,angle= 30]a1){\rnode{b1}{$+$}}
\rput([nodesep=.7,angle=150]a1){\rnode{b2}{$+$}}
\rput([nodesep=.7,angle=330]a2){\rnode{b3}{$-$}}
\rput([nodesep=.7,angle=210]a2){\rnode{b4}{$-$}}
\ncline{a1}{b1}\ncline{a1}{b2}\ncline{a1}{a2}
\ncline{a2}{b3}\ncline{a2}{b4}
\endpspicture
+ q^{-1/3}
\pspicture[.4](-1,-.9)(1,.9)
\rput(.7; 45){\rnode{a1}{$+$}}\rput(.7;135){\rnode{a2}{$+$}}
\rput(.7;225){\rnode{a3}{$-$}}\rput(.7;315){\rnode{a4}{$-$}}
\ncarc[arcangle=-45]{a3}{a2}
\ncarc[arcangle=-45]{a1}{a4}
\endpspicture
\\
\pspicture[.4](-1,-.9)(1,.9)
\rput(.7;45){\rnode{a1}{$+$}}
\rput(.7;135){\rnode{a2}{$+$}}
\rput(.7;225){\rnode{a3}{$-$}}
\rput(.7;315){\rnode{a4}{$-$}}
\ncline{a3}{a1}
\ncline[border=.1]{a2}{a4}
\endpspicture
& = & q^{-1/6}
\pspicture[.4](-1,-.9)(1,.9)
\pnode(0,.35){a1}\pnode(0,-.35){a2}
\rput([nodesep=.7,angle= 30]a1){\rnode{b1}{$+$}}
\rput([nodesep=.7,angle=150]a1){\rnode{b2}{$+$}}
\rput([nodesep=.7,angle=330]a2){\rnode{b3}{$-$}}
\rput([nodesep=.7,angle=210]a2){\rnode{b4}{$-$}}
\ncline{a1}{b1}\ncline{a1}{b2}\ncline{a1}{a2}
\ncline{a2}{b3}\ncline{a2}{b4}
\endpspicture
+ q^{1/3}
\pspicture[.4](-1,-.9)(1,.9)
\rput(.7; 45){\rnode{a1}{$+$}}\rput(.7;135){\rnode{a2}{$+$}}
\rput(.7;225){\rnode{a3}{$-$}}\rput(.7;315){\rnode{a4}{$-$}}
\ncarc[arcangle=-45]{a3}{a2}
\ncarc[arcangle=-45]{a1}{a4}
\endpspicture,
\end{eqnarray*}
four types in the $B_2$ spider:
\begin{eqnarray*}
\rcrossing & = & -q^{1/2}\vv - q^{-1/2}\hh
+ {1 \over q^{1/2} + q^{-1/2}}\btwox
\\
\pspicture[.4](-.6,-.5)(.6,.5)
\psline(.5;135)(.5;315)
\psline[doubleline=true,border=.1](.5;45)(.5;225)
\endpspicture
& = &
{q^{-1/2} \over q^{1/2} + q^{-1/2}}\doublehvert
+ {q^{1/2} \over q^{-1/2} + q^{1/2}}\doublehhoriz
\\
\pspicture[.4](-.6,-.5)(.6,.5)
\psline[doubleline=true](.5;45)(.5;225)
\psline[border=.1](.5;135)(.5;315)
\endpspicture
& = &
{q^{1/2} \over q^{1/2} + q^{-1/2}}\doublehvert
+ {q^{-1/2} \over q^{-1/2} + q^{1/2}}\doublehhoriz
\\
\pspicture[.4](-.6,-.5)(.6,.5)
\psline[doubleline=true](.5;135)(.5;315)
\psline[doubleline=true,border=.1](.5;45)(.5;225)
\endpspicture
& = & q
\pspicture[.4](-.6,-.5)(.6,.5)
\psbezier[doubleline=true](.5;45)(.25;45)(.25;315)(.5;315)
\psbezier[doubleline=true](.5;225)(.25;225)(.25;135)(.5;135)
\endpspicture
+ q^{-1}
\pspicture[.4](-.6,-.5)(.6,.5)
\psbezier[doubleline=true](.5;45)(.25;45)(.25;135)(.5;135)
\psbezier[doubleline=true](.5;225)(.25;225)(.25;315)(.5;315)
\endpspicture
+ {1 \over q + 2 +q^{-1}} \doublesquare
\end{eqnarray*}
and four types in the $G_2$ spider
\begin{eqnarray*}
\rcrossing & = &
{q^{-1/2} \over q^{1/2} + q^{-1/2}}\gtwohvert
+ {q^{1/2} \over q^{1/2} + q^{-1/2}}\gtwohhoriz
+ {q^{-3/2} \over q^{1/2} + q^{-1/2}}\hh
+ {q^{3/2} \over q^{1/2} + q^{-1/2}}\vv
\\
\pspicture[.4](-.6,-.5)(.6,.5)
\psline(.5;135)(.5;315)
\psline[doubleline=true,border=.1](.5;45)(.5;225)
\endpspicture
& = &
-{q^{-3/2} \over q^{1/2} + q^{-1/2}}\doublehvert
-{q^{3/2} \over q^{1/2} + q^{-1/2}}\doublehhoriz
-{1\over q+2+q^{-1}}\semidoublesquare
\\
\pspicture[.4](-.6,-.5)(.6,.5)
\psline[doubleline=true](.5;45)(.5;225)
\psline[border=.1](.5;135)(.5;315)
\endpspicture
& = &
-{q^{3/2} \over q^{1/2} + q^{-1/2}}\doublehvert
-{q^{-3/2} \over q^{1/2} + q^{-1/2}}\doublehhoriz
-{1\over q+2+q^{-1}}\semidoublesquare
\\
\pspicture[.4](-.6,-.5)(.6,.5)
\psline[doubleline=true](.5;135)(.5;315)
\psline[doubleline=true,border=.1](.5;45)(.5;225)
\endpspicture
& = & {q^5 - q^3 + q - q^{-3} \over q^2 - 1 + q^{-2}}
\pspicture[.4](-.6,-.5)(.6,.5)
\psbezier[doubleline=true](.5;45)(.25;45)(.25;315)(.5;315)
\psbezier[doubleline=true](.5;225)(.25;225)(.25;135)(.5;135)
\endpspicture
+ {- q^3 + q - q^{-3} - q^{-5} \over q^2 - 1 + q^{-2}}
\pspicture[.4](-.6,-.5)(.6,.5)
\psbezier[doubleline=true](.5;45)(.25;45)(.25;135)(.5;135)
\psbezier[doubleline=true](.5;225)(.25;225)(.25;315)(.5;315)
\endpspicture
+ {1 \over q + 2 +q^{-1}} \doublesquare
\end{eqnarray*}
Note that crossings are the first webs
whose coefficients are not symmetric in $q$ and $q^{-1}$.

\subsection{Clasps and clasped web spaces}

The rank 2 spiders also admit clasps and clasped web spaces. As before, the
type of a clasp is the same as an unclasped strand type, but in the rank 2
cases there are many more combinations. Our notation will be that if $s$ an
unclasped strand type, $c = [s]$ is the corresponding clasp. We will also
consider sequences of clasps $C = c_1c_2...c_n = [s_1][s_2]\ldots[s_n]$.

Clasps and clasped web spaces in the $A_2$ spider are the easiest
to describe:  Define the weight $\wt(s)$ of a sign string $s$
$n\lambda_1 + k\lambda_2$ if there are $n$ plusses and $k$ minuses.  Recall
the usual partial ordering of the weight lattice of lattice of $A_2$:  it
is generated by
\begin{eqnarray*}
a\lambda_1 + b\lambda_2 & \succ & (a+1)\lambda_1 + (b-2)\lambda_2 \\
a\lambda_1 + b\lambda_2 & \succ & (a-2)\lambda_1 + (b+1)\lambda_2.
\end{eqnarray*}
There is a clasped web space $W(C)$ for each possible clasp sequence $C$.
For example, $W([+ - +][- + -][- +])$ denotes a clasped web space.  The basis
$B([s_1][s_2]\ldots[s_k])$ of $W([s_1][s_2]\ldots[s_k])$ is a subset of the
set of unclasped basis webs $B(s_1 s_2\ldots s_k)$ consisting of those
non-elliptic webs with non-convex clasps. Here an (external) clasp is
{\df non-convex} if it has the property that the weight of any path between
endpoints of the clasp that is transverse to the web (a {\df cut path}) is greater
than or equal to the weight of the clasp; the weight of such a path is
defined as the number of strands that cross it in the direction of the clasp
and the number that cross away from the clasp. For example, the following web
has a partly convex clasp, because its weight is $2\lambda_1 + \lambda_2$,
but there is an arc with weight $2\lambda_2$ that cuts off the clasp:
$$
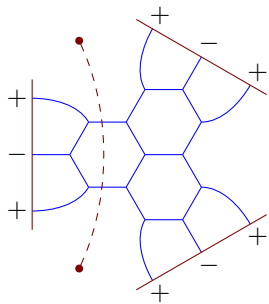
\pspicture(-2,-2)(2,2)
\pnode(0,.866){a1}\pnode(.5,.866){a2}
\pnode(-.75,.433){a3}\pnode(-.25,.433){a4}\pnode(.75,.433){a5}
\pnode(-1,0){a6}\pnode(0,0){a7}\pnode(.5,0){a8}
\pnode(-.75,-.433){a9}\pnode(-.25,-.433){a10}\pnode( .75,-.433){a11}
\pnode(0,-.866){a12}\pnode(.5,-.866){a13}
\pnode(1.5; 60){b2}\pnode([nodesep=.75,angle=150]b2){b3}
\pnode([nodesep=.75,angle=330]b2){b1}
\pnode(1.5;180){b5}\pnode([nodesep=.75,angle=270]b5){b6}
\pnode([nodesep=.75,angle= 90]b5){b4}
\pnode(1.5;300){b8}\pnode([nodesep=.75,angle= 30]b8){b9}
\pnode([nodesep=.75,angle=210]b8){b7}
\ncline{a1}{a2}\ncline{a1}{a4}\ncline{a2}{a5}
\ncline{a3}{a4}\ncline{a3}{a6}\ncline{a4}{a7}\ncline{a5}{a8}
\ncline{a6}{a9}\ncline{a7}{a8}\ncline{a7}{a10}\ncline{a8}{a11}
\ncline{a9}{a10}\ncline{a10}{a12}\ncline{a11}{a13}\ncline{a12}{a13}
\ncline{a2}{b2}\nccurve[angleA=120,angleB=240]{a1}{b3}
\nccurve[angleA=0,angleB=240]{a5}{b1}
\ncline{a6}{b5}\nccurve[angleA=120,angleB=0]{a3}{b4}
\nccurve[angleA=240,angleB=0]{a9}{b6}
\ncline{a13}{b8}\nccurve[angleA=240,angleB=120]{a12}{b7}
\nccurve[angleA=0,angleB=120]{a11}{b9}
\psline[linecolor=darkred]([nodesep=.25,angle=150]b3)([nodesep=.25,angle=330]b1)
\psline[linecolor=darkred]([nodesep=.25,angle=270]b6)([nodesep=.25,angle= 90]b4)
\psline[linecolor=darkred]([nodesep=.25,angle= 30]b9)([nodesep=.25,angle=210]b7)
\rput([nodesep=.2,angle= 60]b1){$+$}
\rput([nodesep=.2,angle= 60]b2){$-$}
\rput([nodesep=.2,angle= 60]b3){$+$}
\rput([nodesep=.2,angle=180]b4){$+$}
\rput([nodesep=.2,angle=180]b5){$-$}
\rput([nodesep=.2,angle=180]b6){$+$}
\rput([nodesep=.2,angle=300]b7){$+$}
\rput([nodesep=.2,angle=300]b8){$-$}
\rput([nodesep=.2,angle=300]b9){$+$}
\pcarc[arrows=*-*,linecolor=darkred,linestyle=dashed,arcangle=25]
    (1.75;120)(1.75;240)
\endpspicture
$$
It is therefore not a basis web of $W([+ - +][+ - +][+ - +])$. Clasped web
spaces can again be interpreted as both quotients and subspaces of unclasped
web spaces.  Any web with at least one partly convex clasp is understood as
the zero vector.  By this convention and equations~(\ref{fa2elliptic}), any
trivalent graph in a disk with suitable boundary can be interpreted as some
vector in the clasped web space.

The non-convexity condition for clasps may seems unnecessarily strong. One
might alternatively stipulate that every cut path cross at least as many
strands as the number of strands in the clasp, or that the weight of no cut
path be strictly less than the weight of the clasp.  However, by
Lemma~\ref{lcomparable}, the two condictions are both equivalent to the one
given. Say that a cut path of a web is {\df minimal} if its weight is minimal
with respect to the partial ordering.  Then in particular,
Lemma~\ref{lcomparable} implies that all minimal cut paths with a fixed pair
of endpoints have the same weight, the minimal cut weight.

As usual, join and rotation in the clasped $A_2$ spider are straightforward,
but stitch involves internal clasps.  An {\df internal clasp of type $s$} is
defined as an idempotent in the unclasped web space $W(s^*s)$ that
annihilates any web in
$W(s^*t)$ with $\wt(t) \prec \wt(s)$.  For example,
$$
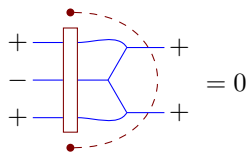
\pspicture[.42](-.5,-.9)(2.2,.9)
\rput(-.2,-.5){$+$}\rput(-.2,0){$-$}\rput(-.2,.5){$+$}
\qline(0,-.5)(.4,-.5)\qline(0,0)(.4,0)\qline(0,.5)(.4,.5)
\psframe[linecolor=darkred](.4,-.7)(.6,.7)
\qline(.6,0)(1,0)
\pccurve[angleA=0,angleB=240](.6,-.5)(1.25,-.433)
\pccurve[angleA=0,angleB=120](.6,.5)(1.25,.433)
\qline(1,0)(1.25,.433)\qline(1,0)(1.25,-.433)
\qline(1.25,.433)(1.75,.433)\qline(1.25,-.433)(1.75,-.433)
\rput(1.95,.433){$+$}\rput(1.95,-.433){$+$}
\pccurve[arrows=*-*,linestyle=dashed,linecolor=darkred,
    angleA=0,angleB=0,ncurv=1.7](.5,-.9)(.5,.9)
\endpspicture = 0
$$
It is not clear that internal clasps exist; without them, we must understand
the clasped spider not as a spider but as a collection of web spaces with the
operations of rotation and join only.  We will use this incomplete structure
in an indirect argument that internal clasps must exist for all three rank 2
cases (as before, they are highest-weight projections), but in the $A_2$ case
we will also give a more explicit construction.

The construction of the clasped $B_2$ spider is entirely analagous to that of
the clasped $A_2$ spider, except that the definition of a cut path and its
weight are slightly different. A cut path may cut diagonally through a
tetravalent vertex, and its weight is defined
as $n\lambda_1 + (k + k')\lambda_2$, where
$n$ is the number of type ``1'' strands that it cuts, $k$ is the number of
type ``2'' strands that it cuts, and $k'$ is the number of tetravalent
vertices that it bisects.  For example, the following cut path has weight
$\lambda_1 + 2\lambda_2$:
$$
\pspicture(0,-.75)(2.5,2.25)
\pnode(0,1.5){a1}\pnode(0,1){a2}\pnode(0,.5){a3}\pnode(0,0){a4}
\pnode(.5,1){b2}\pnode([nodesep=.5,angle=45]b2){b1}
\pnode([nodesep=.5,angle=315]b1){b3}
\pnode([nodesep=.5,angle=315]b2){b4}
\pnode([nodesep=.5,angle=315]b3){b5}
\pnode([nodesep=.5,angle=315]b4){b6}
\nccurve[angleA=0,angleB=135]{a1}{b1}
\ncline[doubleline=true]{a2}{b2}
\nccurve[angleA=0,angleB=225]{a3}{b4}
\nccurve[angleA=0,angleB=225]{a4}{b6}
\qline(b2)([nodesep=.5,angle=45]b1)
\qline(b4)([nodesep=.5,angle=45]b3)
\qline(b6)([nodesep=.5,angle=45]b5)
\qline(b1)([nodesep=.5,angle=315]b5)
\qline(b2)([nodesep=.5,angle=315]b6)
\pscurve[linecolor=darkred,linestyle=dashed]{*-*}(0,2)(b1)(b4)(0,-.5)
\psline[linecolor=darkred](0,-.25)(0,1.75)
\endpspicture
$$
Recall that there is a natural partial ordering of the $B_2$ weight lattice
given by
\begin{eqnarray}
a\lambda_1 + b\lambda_2 &\succ &(a-2)\lambda_1 + (b+1)\lambda_2 \\
a\lambda_1 + b\lambda_2 &\succ &(a+2)\lambda_1 + (b-2)\lambda_2.
\end{eqnarray}

The clasped $G_2$ spider has a more significant difference.  A cut
path may contain a type ``1'' strand in its interior, and its weight is
$n\lambda_1 + (k + k')\lambda_2$, where $n$ is the number of type ``1''
strands that it cuts, $k$ is the number of type ``2'' strands that it cuts,
and $k'$ is the number of type ``1'' strands that it contains.  For example,
the following cut path has weight $\lambda_1 + 2\lambda_2$:
$$
\pspicture(0,-1.4)(1.6,1.4)
\pccurve[angleA=0,angleB=150](0,1)(.933,.75)
\pccurve[angleA=0,angleB=210](0,-.5)(.933,-.25)
\psline[doubleline=true](0,0)(.5,0)
\psline[doubleline=true](.933,-.25)(1.433,-.25)
\psline(.933,-.25)(.5,0)(.933,.25)(.933,.75)(1.366,1)
\qline(.933,.25)(1.366,0)\qline(0,-1)(1.4,-1)
\pccurve[arrows=*-,angleA=0,angleB=90,linestyle=dashed,linecolor=darkred]
    (0,1.4)(.933,.75)
\pccurve[angleA=270,angleB=90,linestyle=dashed,linecolor=darkred]
    (.933,.25)(1.183,-.35)
\pccurve[arrows=-*,angleA=270,angleB=0,linestyle=dashed,linecolor=darkred]
    (1.183,-.35)(0,-1.4)
\psline[linecolor=darkred](0,-1.2)(0,1.2)
\endpspicture
$$
As before, the webs with non-convex clasps, where non-convexity is 
defined using the partial ordering
\begin{eqnarray*}
a\lambda_1 + b\lambda_2 \succ (a-2)\lambda_1 + (b+1)\lambda_2 \\
a\lambda_1 + b\lambda_2 \succ (a+3)\lambda_1 + (b-2)\lambda_2
\end{eqnarray*}
in the $G_2$ weight lattice, form a basis of each clasped web space. The more
significant difference is that a basis element of the unclasped web space
with a convex clasp is not necessarily zero.  Rather, if a web, whether
non-elliptic or not, has a cut path which cuts of a clasp, which does not
contain any type ``1'' strands and whose weight is less than that of the
clasp, then the web is zero.  For example,
$$
\pspicture(0,-.5)(3,2)
\pnode(0,1){a1}\pnode(0,.5){a2}\pnode(0,0){a3}
\pnode(.5,.5){b3}\pnode([nodesep=.5,angle=30]b3){b1}
\pnode([nodesep=.5,angle=0]b1){b2}
\pnode([nodesep=.5,angle=330]b2){b4}
\pnode([nodesep=.5,angle=330]b3){b5}
\pnode([nodesep=.5,angle=0]b5){b6}
\ncline[doubleline=true]{b1}{b2}
\ncline[doubleline=true]{a2}{b3}
\ncline[doubleline=true]{b5}{b6}
\psline[doubleline=true](b4)([nodesep=.5,angle=0]b4)
\ncline{b1}{b3}\ncline{b2}{b4}
\ncline{b3}{b5}\ncline{b4}{b6}
\qline(b2)([nodesep=.5,angle=30]b2)
\qline(b6)([nodesep=.5,angle=330]b6)
\nccurve[angleA=0,angleB=150]{a1}{b1}
\nccurve[angleA=0,angleB=210]{a3}{b5}
\psline[linecolor=darkred](0,-.25)(0,1.25)
\pscurve[linecolor=darkred,linestyle=dashed]{*-*}(0,-.5)
   ([nodesep=.25,angle=0]b5)([nodesep=.25,angle=0]b1)(0,1.5)
\endpspicture
$$
It may not be immediate that the kernel of this quotienting operation
complements the subspace defined as the clasped web space.  This
will be shown in Section~\ref{sb2g2}.

This concludes the definition of the combinatorial rank 2 spiders
and the definition of clasped web spaces.  Only the operation
of stitch, which depends on the existence of internal clasps,
remains to be fully defined.

\section{The morphism from combinatorial to algebraic}
\label{smorphism}

Let $V_+$ and $V_-$ be the two fundamental representations of $A_2 \cong
\sl(3)$, with $V_- \cong (V_+)^*$.  Let $V_1$ and $V_2$ be the two
fundamental representations of $B_2 \cong \sp(4) \cong \so(5)$, and give the two
fundamental representations of $G_2$ the same names with $\dim V_1 < \dim
V_2$ in both cases.  Then the vector spaces
$$\begin{array}{l}
\Inv_{A_2}(V_+^{\tensor 3}) \\ \Inv_{A_2}(V_-^{\tensor 3}) \\
\Inv_{B_2}(V_1 \tensor V_1 \tensor V_2) \\ \Inv_{G_2}(V_1^{\tensor 3}) \\
\Inv_{G_2}(V_1 \tensor V_1 \tensor V_2)
\end{array}$$
are all 1-dimensional.  When $q=1$, and for each of $A_2$, $B_2$, and
$G_2$, there exists a morphism $\Phi$ from the unclasped combinatorial spider
to the unclasped algebraic spider with the following property.  If $s$ a
generator of the strand set, then $\Phi(s) = V_s$, and if $T$ is a trivalent
vertex of some type, then $\Phi(T)$ is a non-zero element in one of the above
invariant spaces.  We sketch the argument for the existence of $\Phi$ in the
$A_2$ case (see Reference \citen{Kuperberg:g2} for details):  Pick any two
non-zero elements $x \in \Inv_{A_2}(V_+^{\tensor 3})$ and $x^* \in
\Inv_{A_2}(V_-^{\tensor 3})$.  By counting dimensions of invariant spaces,
each of the left sides of equations~(\ref{fa2elliptic}),
(\ref{fb2elliptic}), and (\ref{fg2elliptic}) must be some
linear combinations of the right sides in the algebraic $A_2$ spider, if $x$
and $x^*$ are denoted by the usual trivalent vertices.  But at the same time,
a computation shows that the right sides are linearly independent in the
algebraic $A_2$ spider, and that, up to normalization of $x$ and $x^*$, the
given coefficients are the only ones that respect this normalization. Thus,
after rescaling, the invariant tensors $x$ and $x^*$ of $U_q(A_2)$ must
satisfy the relations of the combinatorial $A_2$ spider, although perhaps
with a different choice of $q$.  Another simple computation shows that the
choice of $q$ is the same.  Thus we can set $\Phi(T) = x$ and $\Phi(T^*) =
x^*$, if $T$ and $T^*$ are the trivalent vertices that generate $A_2$.

\begin{theorem} The morphism $\Phi$ from the combinatorial to the algebraic
$A_2$, $B_2$, or $G_2$ spider is surjective when $q=1$, and therefore
for generic $q$.
\end{theorem}
This theorem is proved in Reference~\citen{Kuperberg:g2}, but we 
give a more conceptual argument here:

\begin{proof}
Let $q=1$, let $\frak g$ be $A_2$, $B_2$, or $G_2$, and let $G$ be the
compact, simply-connected Lie group whose complexified Lie algebra is $\frak
g$.  If $s = s_1\ldots s_n$ is a string, define
$$V_s = V_{s_1} \tensor \ldots V_{s_n}.$$
In each case, the image $\cal X$ of $\Phi$ is some subspider of the algebraic
spider of $\frak g$.  We interpret $\cal X$ as a category of linear
transformations between tensor products of fundamental representations. The
category $\cal X$ contains switching maps $x \tensor y \mapsto y \tensor x$,
since they are the images of crossings under the map $\Phi$.  Moreover, each
$\End_{\cal X}(V_s)$ is a semisimple algebra, by the following construction:
We define the Hermitian adjoint $w^*$ of $w \in W(ss^*)$ by reflecting $w$
about a vertical axis, taking the complex conjugate of all coefficients, and,
in the $A_2$ spider, reversing all orientations and signs:
$$
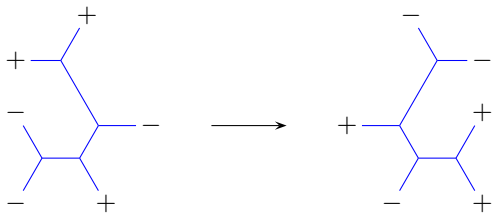
\pspicture[.5](-1.5,-1.5)(1.5,1.5)
\pnode(.5;0){a1}\pnode(.5;120){a3}
\pnode(.5;180){a4}\pnode(.5;240){a5}\pnode(.5;300){a6}
\pnode(1;0){b1}\pnode(1;120){b3}
\pnode(1;180){b4}\pnode(1;240){b5}\pnode(1;300){b6}
\pnode(-.4,.866){c}\pnode(0,.866){a2}\pnode([nodesep=.5,angle=60]a2){b2}
\qline(c)(a2)\qline(a1)(a2)\qline(a1)(b1)\qline(a2)(b2)
\qline(a4)(a5)\qline(a5)(a6)\qline(a6)(a1)
\qline(a5)(b5)\qline(a6)(b6)
\rput([nodesep=.2,angle=0]b1){$-$}
\rput([nodesep=.2,angle=60]b2){$+$}
\rput([nodesep=.2,angle=240]b5){$-$}
\rput([nodesep=.2,angle=300]b6){$+$}
\rput([nodesep=.2,angle=120]a4){$-$}
\rput([nodesep=.2,angle=180]c){$+$}
\endpspicture
\psgoesto
\pspicture[.5](-1.5,-1.5)(1.5,1.5)
\pnode(.5;180){a1}\pnode(.5;60){a3}
\pnode(.5;0){a4}\pnode(.5;300){a5}\pnode(.5;240){a6}
\pnode(1;180){b1}\pnode(1;60){b3}
\pnode(1;0){b4}\pnode(1;300){b5}\pnode(1;240){b6}
\pnode(.4,.866){c}\pnode(0,.866){a2}\pnode([nodesep=.5,angle=120]a2){b2}
\qline(c)(a2)\qline(a1)(a2)\qline(a1)(b1)\qline(a2)(b2)
\qline(a4)(a5)\qline(a5)(a6)\qline(a6)(a1)
\qline(a5)(b5)\qline(a6)(b6)
\rput([nodesep=.2,angle=180]b1){$+$}
\rput([nodesep=.2,angle=120]b2){$-$}
\rput([nodesep=.2,angle=300]b5){$+$}
\rput([nodesep=.2,angle=240]b6){$-$}
\rput([nodesep=.2,angle=60]a4){$+$}
\rput([nodesep=.2,angle=0]c){$-$}
\endpspicture
$$
The morphism $\Phi$ intertwines the combinatorial Hermitian adjoint
with the usual Hermitian adjoint in the algebraic spider.  Since
$\End_{\cal X}(V_s)$ is closed under Hermitian adjoint, it must
be semisimple.

The category $\cal X$ does not contain the kernels and co-kernels of its
morphisms, but it may be completed to a bigger category $\cal X'$, in which
$\cal X$ is a full subcategory, by adding these vector spaces as new objects.
The category $\cal X'$ then satisfies the hypotheses of the Tannaka-Krein
duality theorem \cite[p. 177]{Kirillov}, and must be the category of
finite-dimensional representations of some compact group $H$. On the one
hand, $H \subseteq G$, since everything in $\cal X'$ is invariant under $G$. 
On the other hand, $H$ cannot be any bigger than $G$, because for each choice
of $\frak g$, $G$ is a maximal compact subgroup of the symmetry group of
$\Phi(t)$ for a vertex $t$. (For example, the symmetry group of a non-zero
element of $\Inv_{B_2}(V_1 \tensor V_1 \tensor V_2)$ is $\sp(4,\C)$, with
maximal compact subgroup $\Spin(5)$.)  Therefore $\cal X'$ is equivalent to
the representation category of $U(\frak g)$, and $\cal X$ coincides with the
algebraic $\frak g$ spider.
\end{proof}

It is more difficult to show that $\Phi$ is injective.  We will prove this in
Section~\ref{sequinumeration} by demonstrating that the clasped web spaces
and the web spaces of the clasped algebraic spider have equinumerous bases. 
If this is so for clasped web spaces, then it is also true for unclasped web
spaces, which demonstrates that $\Phi$ is an isomorphism between unclasped
spiders. We can then define an internal clasp as $\Phi^{-1}(\pi)$, where
$\pi$ is the highest-weight projection from any strand in the algebraic
$\frak g$ spider to itself. This completes the definition of the clasped
spiders, provided we verify the following lemma to show that $\Phi$ maps
stitch to contraction.

\begin{lemma} Suppose that $s$ and $t$ are strand types and
$w \in W(st)$ is zero in the clasped web space $W([s]t)$.  Then
$w$ is annihilated by an internal clasp of type $[s]$.
\end{lemma}
\begin{proof}
We assume the injectivity of $\Phi$ for unclasped web spaces, and we assume
that $w \in B(st)$ is a basis web.  Let $\lambda$ be the weight of $s$, and
let $w'$ be $w$ with an internal clasp attached The tensor $\Phi(w)$ may be
interpreted as a homomorphism, in particular as a composotion $V_s \to
V(\lambda) \to V_u \to V_t$, where $u$ is the transverse strand type of a
minimal cut path that separates $s$ from $t$ in $w$. By hypothesis, the
weight of $u$ is lower than that of $s$. Therefore any map $V(\lambda) \to
V_u$ must vanish, so $w'$ must vanish.
\end{proof}

The following result is then a corollary of the discussion of this
section:

\begin{theorem} The morphism $\Phi$ from the combinatorial to algebraic
rank two spiders extends to clasped spiders.
\end{theorem}

\section{Equinumeration}
\label{sequinumeration}

\subsection{The $A_2$ case}

\begin{theorem} Let $C$ be an $A_2$ clasp sequence and
let $\wt C$ be the corresponding sequence of weights.
Then the vector spaces $W(C)$ and $\Inv(V(\wt C))$ have
the same dimension. \label{tha2equinum}
\end{theorem}

One of the main steps of the proof of Theorem~\ref{tha2equinum} is
the same as that of the proof of Theorem~\ref{tha1isom}:  Given
two clasp sequences $C$ and $D$, there is a decomposition of
vector spaces
$$\Inv(V(\wt C) \tensor V(\wt D)) \cong \bigoplus_\mu \Inv(V(\wt S) \tensor V(\mu))
\tensor \Inv(V(\mu)^* \tensor V(\wt D)).$$
We wish to prove the corresponding decomposition of sets of basis webs:

\begin{theorem} Given two clasp sequences $C$ and $D$,
$$B(CD) \cong \bigcup_\lambda \left(B(Cc_\lambda) \times
B(c_\lambda^*D)\right),$$ where for each weight $\lambda$, $c_\lambda$ is
some clasp with weight $\lambda$. \label{tha2splitting}
\end{theorem}

Theorem~\ref{tha2splitting} is more complicated than its analogue for the
$A_1$ spider, because there is no longer always a unique minimal cut path
separating the clasp sequences $C$ and $D$:
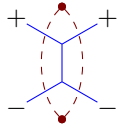
$$
\pspicture[.4](-1,-.9)(1,.9)
\pnode(0,.25){a1}\pnode(0,-.25){a2}
\rput([nodesep=.7,angle= 30]a1){\rnode{b1}{$+$}}
\rput([nodesep=.7,angle=150]a1){\rnode{b2}{$+$}}
\rput([nodesep=.7,angle=330]a2){\rnode{b3}{$-$}}
\rput([nodesep=.7,angle=210]a2){\rnode{b4}{$-$}}
\ncline{a1}{b1}\ncline{a1}{b2}\ncline{a1}{a2}
\ncline{a2}{b3}\ncline{a2}{b4}
\pcarc[arrows=*-*,linecolor=darkred,linestyle=dashed,arcangle=45]
    (0,.75)(0,-.75)
\pcarc[arrows=*-*,linecolor=darkred,linestyle=dashed,arcangle=-45]
    (0,.75)(0,-.75)
\endpspicture
\label{fnonunique}
$$
Although the two cut paths in this example have the same weight, they have
different transverse strand types and they yield webs in different clasped
web spaces.  Yet not every ordering of the strands is always possible.
Nevertheless, there is a way to reconcile these different decompositions and
ameliorate the ordering problem.

Indeed, order independence is a large part of the combinatorial content of
Theorem~\ref{tha2equinum}.  For example, it implies the following result,
which has an independent proof that is another warm-up to the proof
of Theorem~\ref{tha2equinum}.

\begin{theorem} If two sign strings $s_1$ and $s_2$ of the
unclasped $A_2$ spider differ only in order, then $B(s_1)$
and $B(s_2)$ are equinumerous.  \label{lhmap}
\end{theorem}
\begin{proof} The web
$$
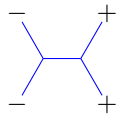
\pspicture(-.9,-.7)(.9,.7)
\pnode(.25,0){a1}\pnode(-.25,0){a2}
\rput([nodesep=.7,angle=60]a1){\rnode{b1}{$+$}}
\rput([nodesep=.7,angle=300]a1){\rnode{b2}{$+$}}
\rput([nodesep=.7,angle=120]a2){\rnode{b3}{$-$}}
\rput([nodesep=.7,angle=240]a2){\rnode{b4}{$-$}}
\ncline{a1}{a2}\ncline{a1}{b1}
\ncline{a1}{b2}\ncline{b3}{a2}
\ncline{b4}{a2}
\endpspicture
$$
is the H-web.  It suffices to consider the case in which $s_1$ and
$s_2$ are the same except for one pair of transposed signs at adjacent
positions $p$ and $q$. Then there is a bijection $h:B(s_1) \to B(s_2)$ that
has the following effect on non-elliptic webs.  If a web $w$ connects $p$ to
$q$ by a bare strand, $h(w)$ is the same web with the orientation of the
strand reversed:
$$
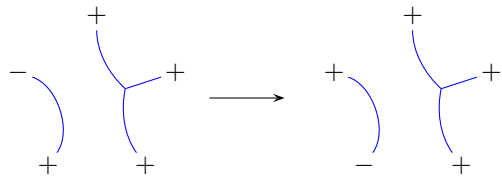
\pspicture[.5](-1.2,-1.2)(1,1.2)\pentanode\ncline{b5}{a5}
\nccurve[angleA=138,angleB=270]{a5}{b1}
\nccurve[angleA=342,angleB= 54]{b2}{b3}
\nccurve[angleA=258,angleB=126]{a5}{b4}
\rput([nodesep=.2,angle= 90]b1){$+$}
\rput([nodesep=.2,angle=162]b2){$-$}
\rput([nodesep=.2,angle=234]b3){$+$}
\rput([nodesep=.2,angle=306]b4){$+$}
\rput([nodesep=.2,angle= 18]b5){$+$}
\endpspicture
\psgoesto
\pspicture[.5](-1.2,-1.2)(1,1.2)\pentanode\ncline{b5}{a5}
\nccurve[angleA=138,angleB=270]{a5}{b1}
\nccurve[angleA=342,angleB= 54]{b2}{b3}
\nccurve[angleA=258,angleB=126]{a5}{b4}
\rput([nodesep=.2,angle= 90]b1){$+$}
\rput([nodesep=.2,angle=162]b2){$+$}
\rput([nodesep=.2,angle=234]b3){$-$}
\rput([nodesep=.2,angle=306]b4){$+$}
\rput([nodesep=.2,angle= 18]b5){$+$}
\endpspicture
$$
If there is an H-web attached at $p$ and $q$, $h(w)$ is $w$
with the H-web removed:
$$
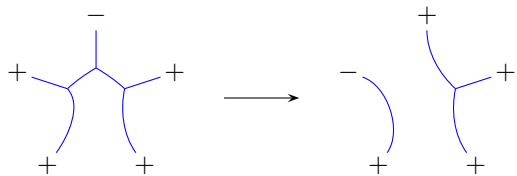
\pspicture[.5](-1.2,-1.2)(1.2,1.2)\pentanode
\ncline{a5}{b5}\ncline{a1}{b1}\ncline{a2}{b2}
\ncarc[arcangle=6]{a2}{a1}\ncarc[arcangle=6]{a1}{a5}
\nccurve[angleA=312,angleB= 54]{a2}{b3}
\nccurve[angleA=258,angleB=126]{a5}{b4}
\rput([nodesep=.2,angle= 90]b1){$-$}
\rput([nodesep=.2,angle=162]b2){$+$}
\rput([nodesep=.2,angle=234]b3){$+$}
\rput([nodesep=.2,angle=306]b4){$+$}
\rput([nodesep=.2,angle= 18]b5){$+$}
\endpspicture
\psgoesto
\pspicture[.5](-1.2,-1.2)(1,1.2)\pentanode\ncline{b5}{a5}
\nccurve[angleA=138,angleB=270]{a5}{b1}
\nccurve[angleA=342,angleB= 54]{b2}{b3}
\nccurve[angleA=258,angleB=126]{a5}{b4}
\rput([nodesep=.2,angle= 90]b1){$+$}
\rput([nodesep=.2,angle=162]b2){$-$}
\rput([nodesep=.2,angle=234]b3){$+$}
\rput([nodesep=.2,angle=306]b4){$+$}
\rput([nodesep=.2,angle= 18]b5){$+$}
\endpspicture
$$
If there is no H-web, then $h(w)$ is $w$ with an H-web attached:
$$
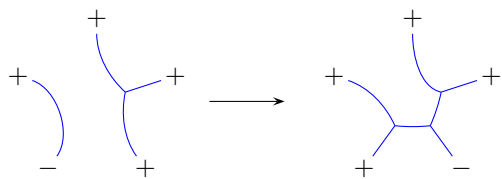
\pspicture[.5](-1.2,-1.2)(1,1.2)\pentanode\ncline{b5}{a5}
\nccurve[angleA=138,angleB=270]{a5}{b1}
\nccurve[angleA=342,angleB= 54]{b2}{b3}
\nccurve[angleA=258,angleB=126]{a5}{b4}
\rput([nodesep=.2,angle= 90]b1){$+$}
\rput([nodesep=.2,angle=162]b2){$+$}
\rput([nodesep=.2,angle=234]b3){$-$}
\rput([nodesep=.2,angle=306]b4){$+$}
\rput([nodesep=.2,angle= 18]b5){$+$}
\endpspicture
\psgoesto
\pspicture[.5](-1.2,-1.2)(1,1.2)\pentanode
\ncline{a3}{b3}\ncline{a4}{b4}\ncline{a5}{b5}
\ncarc[arcangle=6]{a5}{a4}\ncarc[arcangle=6]{a4}{a3}
\nccurve[angleA=168,angleB=270]{a5}{b1}
\nccurve[angleA=114,angleB=342]{a3}{b2}
\rput([nodesep=.2,angle= 90]b1){$+$}
\rput([nodesep=.2,angle=162]b2){$+$}
\rput([nodesep=.2,angle=234]b3){$+$}
\rput([nodesep=.2,angle=306]b4){$-$}
\rput([nodesep=.2,angle= 18]b5){$+$}
\endpspicture
$$
It is easy to check that $h(w)$ is non-elliptic if $w$ is,
and that $h$ has an inverse.  In fact, the inverse is also an $h$-map.
\end{proof}

In the notation of the above proof, we also say that $w$
and $h(w)$ differ by an H-move.

\begin{lemma} If $\alpha$ is a minimal cut path of $w \in B(ST)$ separating
$S$ from $T$,
it divides $w$ into two parts $w_1 \in B(Sc)$ and $w_2 \in B(c^*T)$
with non-convex clasps $c$ and $c^*$.
\end{lemma}

\begin{proof} If $c$ were convex in $w_1$, $w_1$ would have a cut
path $\alpha'$ whose weight is lower than that of $c$, the same
as the weight of $\alpha$.  But $\alpha'$ is also a cut path in $w$ and
has the same endpoints as $\alpha$, contradicting the hypothesis
that $\alpha$ is minimal:
$$
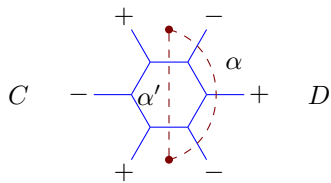
\pspicture[.42](-2.5,-1.5)(2.5,1.5)
\pnode(.5;  0){a1}\pnode(.5; 60){a2}\pnode(.5;120){a3}
\pnode(.5;180){a4}\pnode(.5;240){a5}\pnode(.5;300){a6}
\pnode(1;  0){b1}\pnode(1; 60){b2}\pnode(1;120){b3}
\pnode(1;180){b4}\pnode(1;240){b5}\pnode(1;300){b6}
\pnode([nodesep=.5,angle=60]a3){c1}
\pnode([nodesep=.5,angle=240]a6){c2}
\rput([nodesep=.2,angle=  0]b1){$+$}
\rput([nodesep=.2,angle= 60]b2){$-$}
\rput([nodesep=.2,angle=120]b3){$+$}
\rput([nodesep=.2,angle=180]b4){$-$}
\rput([nodesep=.2,angle=240]b5){$+$}
\rput([nodesep=.2,angle=300]b6){$-$}
\rput(-2,0){$C$}\rput(2,0){$D$}
\qline(a1)(b1)\qline(a2)(b2)\qline(a3)(b3)
\qline(a4)(b4)\qline(a5)(b5)\qline(a6)(b6)
\qline(a1)(a2)\qline(a2)(a3)\qline(a3)(a4)
\qline(a4)(a5)\qline(a5)(a6)\qline(a6)(a1)
\ncline[arrows=*-*,linecolor=darkred,linestyle=dashed]{c2}{c1}
\Aput[.1]{$\alpha'$}
\ncarc[linecolor=darkred,linestyle=dashed,arcangle=75,ncurv=1]{c1}{c2}
\aput(.35){$\alpha$}
\endpspicture
$$
\end{proof}

\begin{lemma} If $\alpha$ and $\beta$ are cut paths from $p$ to $q$ of a basis
web $w$ and $\alpha$ is minimal, then the weight of $\alpha$ is less than or
equal to (and not incomparable to) the weight of $\beta$.  If $\beta$ is also
minimal, the two parts of $w$ cut by $\beta$ are the same as those of $w$ cut
by $\alpha$ up to H-moves.  \label{lcomparable}
\end{lemma}
\begin{proof}
The proof is by induction on the complexity of $w$; assume that $w$ is a
minimal counterexample.  Assume also that, having chosen $w$, $\alpha$ and
$\beta$ are transverse and intersect minimally.  First, we can discard
any structure of $w$ not bounded by $\alpha \cup \beta$; the new web has no
clasps and is non-elliptic since $w$ is non-elliptic. Second, we claim that
$\alpha$ and $\beta$ do not intersect except at their endpoints.  For
otherwise, let $\alpha'$ be a segment of $\alpha$ between two consecutive
intersections $x_1$ and $x_2$, and let $\beta'$ be the arc of $\beta$ from
$x_1$ to $x_2$ (which are not necessarily consecutive along $\beta$). Then
the region between $\alpha'$ and $\beta'$ is either empty, in which case
$\alpha'$ and $\beta'$ do not intersect minimally; or it constitutes a
smaller counterexample than $w$ for some choice of $\alpha'$; or $\beta'$; or
$\alpha$ and $\beta$ have comparable weight:
$$
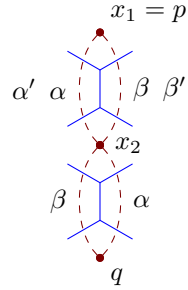
\pspicture(-1,-2)(1,2)
\psline(-.433,.25)(0,.5)(.433,.25)\qline(0,.5)(0,1)
\psline(-.433,1.25)(0,1)(.433,1.25)
\psline(-.433,-.25)(0,-.5)(.433,-.25)\qline(0,-.5)(0,-1)
\psline(-.433,-1.25)(0,-1)(.433,-1.25)
\pcarc[arrows=*-*,linestyle=dashed,linecolor=darkred,arcangle= 45]
    (0,1.5)(0,0)\Aput{$\beta\,\,\,\beta'$}
\pcarc[arrows=*-*,linestyle=dashed,linecolor=darkred,arcangle=-45]
    (0,1.5)(0,0)\Bput{$\alpha'\,\,\,\alpha$}
\pcarc[arrows=*-*,linestyle=dashed,linecolor=darkred,arcangle= 45]
    (0,-1.5)(0,0)\Aput{$\beta$}
\pcarc[arrows=*-*,linestyle=dashed,linecolor=darkred,arcangle=-45]
    (0,-1.5)(0,0)\Bput{$\alpha$}
\uput{.2}[0](0,0){$x_2$}
\uput{.2}[45](0,1.5){$x_1 = p$}
\uput{.2}[315](0,-1.5){$q$}
\endpspicture
$$
Third, $w$ must be connected, for if one of its connected components meets
$\alpha$ but not $\beta$, $\alpha$ is not a minimal cut path; if one of its
components meets $\beta$ but not $\alpha$, it may be discarded to produce a
smaller counterexample; and if all components meet both $\alpha$ and $\beta$,
one of them is a smaller counterexample:
$$
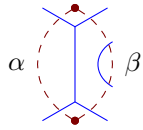
\pspicture(-1,-1.5)(1,1.5)
\qline(0,.5)(0,-.5)
\psline(-.433,.75)(0,.5)(.433,.75)
\psline(-.433,-.75)(0,-.5)(.433,-.75)
\pccurve[angleA=210,angleB=150,ncurv=1](.5,.3)(.5,-.3)
\pcarc[arrows=*-*,linestyle=dashed,linecolor=darkred,arcangle= 60,ncurv=1]
   (0,.75)(0,-.75)\Aput{$\beta$}
\pcarc[arrows=*-*,linestyle=dashed,linecolor=darkred,arcangle=-60,ncurv=1]
   (0,.75)(0,-.75)\Bput{$\alpha$}
\endpspicture
$$

Finally, if $e_1$ and $e_2$ are adjacent endpoints of $w$, define the
exterior curvature of the arc $\overline{e_1e_2}$ to be
$180^\circ - n60^\circ$, where $n$ is the number of vertices
of $w$ connecting $e_1$ to $e_2$:
$$
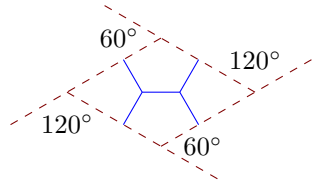
\pspicture(-2,-1.5)(2,1.5)
\pnode(.25,0){a1}\pnode(-.25,0){a2}
\psline([nodesep=.5,angle=60]a1)(a1)([nodesep=.5,angle=300]a1)
\psline([nodesep=.5,angle=120]a2)(a2)([nodesep=.5,angle=240]a2)
\qline(a1)(a2)
\psline[linestyle=dashed,linecolor=darkred](1.25,0)(-.75,1.152)
\psline[linestyle=dashed,linecolor=darkred](-1.25,0)(.75,-1.152)
\psline[linestyle=dashed,linecolor=darkred](0,.72)(-2,-.432)
\psline[linestyle=dashed,linecolor=darkred](0,-.72)(2,.432)
\uput{.3}[90](1.25,0){$\scriptsize 120^\circ$}
\uput{.3}[270](-1.25,0){$\scriptsize 120^\circ$}
\uput{.3}[180](0,.72){$\scriptsize 60^\circ$}
\uput{.3}[0](0,-.72){$\scriptsize 60^\circ$}
\endpspicture
$$
The total exterior curvature of $w$ is at least $360^\circ$ since $w$ is
non-elliptic.  Moreover, the curvature at the arcs containing $p$ and $q$ is
at most $120^\circ$, for if it were $180^\circ$, either $w$ would be
disconnected or it would be a single strand.  Therefore there must be a
segment $\gamma$ of either $\alpha$ or $\beta$ with positive curvature. 
There are three possibilities for the edges of $w$ that bound a face together
with $\gamma$:
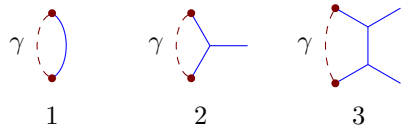
$$
\begin{array}{ccc}
\pspicture[.5](-.75,-.5)(.75,.5)
\pcarc[arcangle=60](0,.433)(0,-.433)
\pcarc[arrows=*-*,linestyle=dashed,linecolor=darkred,arcangle=-60]
    (0,.433)(0,-.433)\Bput{$\gamma$}
\endpspicture &
\pspicture[.5](-.75,-.5)(1,.5)
\qline(0,.433)(.25,0)
\qline(0,-.433)(.25,0)\qline(.25,0)(.75,0)
\pcarc[arrows=*-*,linestyle=dashed,linecolor=darkred,arcangle=-60]
    (0,.433)(0,-.433)\Bput{$\gamma$}
\endpspicture &
\pspicture[.5](-1,-.75)(.75,.75)
\psline(-.433,-.5)(0,-.25)(.433,-.5)\qline(0,-.25)(0,.25)
\psline(-.433,.5)(0,.25)(.433,.5)
\pcarc[arrows=*-*,linestyle=dashed,linecolor=darkred,arcangle=-30]
    (-.433,.5)(-.433,-.5)\Bput{$\gamma$}
\endpspicture \\
1 & 2 & 3
\end{array}
$$
\begin{itemize}
\item[1.] A ``U''.  In this case, $w$ is either disconnected or
a single strand.
\item[2.] A ``Y''.  If $\gamma$ lies on $\alpha$, then $\alpha$ 
is not minimal.  If $\gamma$ lies on $\beta$, then an isotopy of
$\beta$ across the ``Y'' produces a smaller counterexample.
\item[3.] An ``H''.  In this case, an isotopy of either $\alpha$
or $\beta$ produces a smaller counterexample by an H-move.
\end{itemize}
This eliminates all possibilities for the least counterexample.
\end{proof}

\begin{lemma}
Let $C$ be a sequence of clasps, let $c$ be an arbitrary clasp, and let $w \in
B(Cc)$ be a basis web.  There is a cut path $\gamma$ that separates $c$
such that any other such cut path lies between $\gamma$ and $c$.
\label{lecore}
\end{lemma}
\begin{proof}
If $\alpha$ and $\beta$ are two transverse, minimal
cut paths that separate $c$ and $p$ and $q$ are two consecutive transverse
intersection points along either path, then the weight of the arc of $\alpha$
from $p$ to $q$ must equal the weight of the arc of $\beta$ from $p$ to $q$,
for otherwise one path would provide a short-cut for the other and either
$\alpha$ or $\beta$ would not be minimal.  Thus the path $\gamma$ that
follows the perimeter of $\alpha \cup \beta$ must also be minimal:
$$
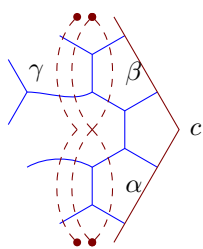
\pspicture(-1,-2)(1.5,2)
\psline(-.433,1.25)(0,1)(.433,1.25)
\psline(0,1)(0,.5)(.433,.25)(.433,-.25)(0,-.5)(0,-1)
\psline(-.433,-1.25)(0,-1)(.433,-1.25)
\qline(.433,.25)(.866,.5)\qline(.433,-.25)(.866,-.5)
\pccurve[angleA=150,angleB=30](0,-.5)(1;210)
\pccurve[angleA=210,angleB=0](0,.5)(1;150)
\pnode(1;150){c}
\psline([nodesep=.5,angle=120]c)(c)([nodesep=.5,angle=240]c)
\psline[linecolor=darkred](.289,1.5)(1.155,0)(.289,-1.5)
\rput[l](1.3,0){$c$}
\pcarc[arrows=*-,linestyle=dashed,linecolor=darkred,arcangle= 45]
    (0,1.5)(0,0)\Aput{$\beta$}
\pcarc[arrows=*-,linestyle=dashed,linecolor=darkred,arcangle=-45]
    (0,1.5)(0,0)
\pcarc[arrows=*-,linestyle=dashed,linecolor=darkred,arcangle= 45]
    (0,-1.5)(0,0)
\pcarc[arrows=*-,linestyle=dashed,linecolor=darkred,arcangle=-45]
    (0,-1.5)(0,0)\Bput{$\alpha$}
\pcarc[arrows=*-,linestyle=dashed,linecolor=darkred,arcangle=-45]
    (-.2,1.5)(-.2,0)\Bput{$\gamma$}
\pcarc[arrows=*-,linestyle=dashed,linecolor=darkred,arcangle= 45]
    (-.2,-1.5)(-.2,0)
\endpspicture
$$
If we partially order cut paths by their distance from $c$, there is
a unique maximal element.
\end{proof}

Given $w$ as in Lemma~\ref{lecore}, define the {\df core} of $w$ relative to
$c$ to be the web $w' \in B(Cc')$ obtained by cutting away $c$ along the cut
path $\gamma$ guaranteed by the lemma. Here $c'$ is another clasp with the
same weight as $c$.  Let $B(C;\lambda)$ be the set of all cores of webs $w
\in B(Cc)$ for all clasps $c$ of weight $\lambda$. Note that, since $c$ is
non-convex in $w$, one of the minimal cut paths that separates $c$ is
parallel to $c$.  Therefore, by Lemma~\ref{lcomparable}, $w$ differs from its
core by H-moves.  Note also that the core $w'$ has the property that no
H-webs are attached at $c'$, and that any web with this property is its own
core.

We describe how an arbitrary basis web extends from one of its cores. Consider
the following {\df stair-step construction} of a basis web from a core.
Suppose that $w \in B(C[s])$ is a core relative to a clasp $[s]$ of weight
$\lambda = a\lambda_1 + b\lambda_2$, and suppose that $[s']$ is an arbitrary
clasp of the same weight. Then the strings $s$ and $s'$ each represent paths
$p$ and $p'$ from the the upper left corner to the lower right corner of an
$a \times b$ rectangle, where each ``$+$'' in $s$ or $s'$ is a step to the
right and each ``$-$'' is a step down:
$$
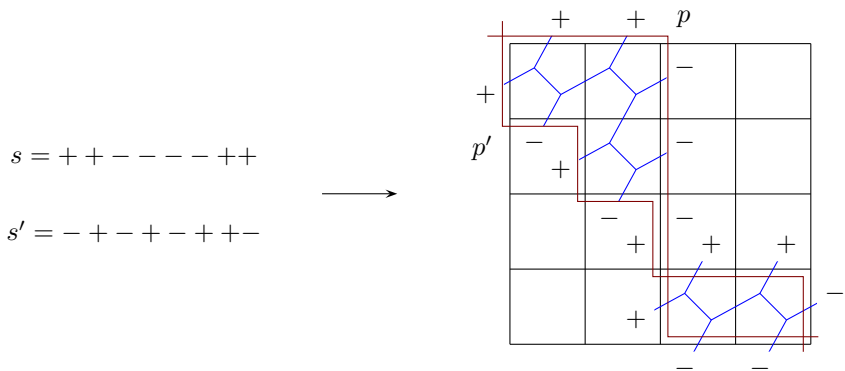
\pspicture[.5](-2,-1)(2,1)
\rput(0,.5){$ s = ++----++ $}
\rput(0,-.5){$ s' = -+-+-++- $}
\endpspicture
\psgoesto
\pspicture[.5](-1,-.5)(4.5,4.5)
\psline[linecolor=black](0,0)(4,0)\psline[linecolor=black](0,0)(0,4)
\psline[linecolor=black](0,1)(4,1)\psline[linecolor=black](1,0)(1,4)
\psline[linecolor=black](0,2)(4,2)\psline[linecolor=black](2,0)(2,4)
\psline[linecolor=black](0,3)(4,3)\psline[linecolor=black](3,0)(3,4)
\psline[linecolor=black](0,4)(4,4)\psline[linecolor=black](4,0)(4,4)
\pnode(0.5,3.5){a11}\pnode(1.5,3.5){a12}\pnode(2.5,3.5){a13}\pnode(3.5,3.5){a14}
\pnode(0.5,2.5){a21}\pnode(1.5,2.5){a22}\pnode(2.5,2.5){a23}\pnode(3.5,2.5){a24}
\pnode(0.5,1.5){a31}\pnode(1.5,1.5){a32}\pnode(2.5,1.5){a33}\pnode(3.5,1.5){a34}
\pnode(0.5,0.5){a41}\pnode(1.5,0.5){a42}\pnode(2.5,0.5){a43}\pnode(3.5,0.5){a44}
\pnode(-.5,3.5){a10}\pnode(0.5,4.5){a01}\pnode(1.5,4.5){a02}
\pnode(4.5,0.5){a45}\pnode(2.5,-.5){a53}\pnode(3.5,-.5){a54}
\pnode([nodesep=.25,angle=135]a11){b11}\pnode([nodesep=.25,angle=315]a11){c11}
\pnode([nodesep=.25,angle=135]a12){b12}\pnode([nodesep=.25,angle=315]a12){c12}
\pnode([nodesep=.25,angle=135]a13){b13}\pnode([nodesep=.25,angle=315]a13){c13}
\pnode([nodesep=.25,angle=135]a14){b14}\pnode([nodesep=.25,angle=315]a14){c14}
\pnode([nodesep=.25,angle=135]a21){b21}\pnode([nodesep=.25,angle=315]a21){c21}
\pnode([nodesep=.25,angle=135]a22){b22}\pnode([nodesep=.25,angle=315]a22){c22}
\pnode([nodesep=.25,angle=135]a23){b23}\pnode([nodesep=.25,angle=315]a23){c23}
\pnode([nodesep=.25,angle=135]a24){b24}\pnode([nodesep=.25,angle=315]a24){c24}
\pnode([nodesep=.25,angle=135]a31){b31}\pnode([nodesep=.25,angle=315]a31){c31}
\pnode([nodesep=.25,angle=135]a32){b32}\pnode([nodesep=.25,angle=315]a32){c32}
\pnode([nodesep=.25,angle=135]a33){b33}\pnode([nodesep=.25,angle=315]a33){c33}
\pnode([nodesep=.25,angle=135]a34){b34}\pnode([nodesep=.25,angle=315]a34){c34}
\pnode([nodesep=.25,angle=135]a41){b41}\pnode([nodesep=.25,angle=315]a41){c41}
\pnode([nodesep=.25,angle=135]a42){b42}\pnode([nodesep=.25,angle=315]a42){c42}
\pnode([nodesep=.25,angle=135]a43){b43}\pnode([nodesep=.25,angle=315]a43){c43}
\pnode([nodesep=.25,angle=135]a44){b44}\pnode([nodesep=.25,angle=315]a44){c44}
\pnode([nodesep=.25,angle=135]a10){b10}\pnode([nodesep=.25,angle=315]a10){c10}
\pnode([nodesep=.25,angle=135]a01){b01}\pnode([nodesep=.25,angle=315]a01){c01}
\pnode([nodesep=.25,angle=135]a02){b02}\pnode([nodesep=.25,angle=315]a02){c02}
\pnode([nodesep=.25,angle=135]a45){b45}\pnode([nodesep=.25,angle=315]a45){c45}
\pnode([nodesep=.25,angle=135]a53){b53}\pnode([nodesep=.25,angle=315]a53){c53}
\pnode([nodesep=.25,angle=135]a54){b54}\pnode([nodesep=.25,angle=315]a54){c54}

\psline(b11)(c11)(b12)(c12)(b22)(c22)
\psline(b43)(c43)(b44)(c44)
\qline(b21)(c11)\qline(c21)(b22)\qline(c12)(b13)
\psline(b32)(c22)(b23)
\psline(c42)(b43)(c33)\qline(b44)(c34)
\psline(c10)(b11)(c01)\qline(b12)(c02)
\qline(b53)(c43)\psline(b54)(c44)(b45)
\rput*(c01){$+$}\rput*(c02){$+$}\rput*(b13){$-$}\rput*(b23){$-$}
\rput*(b33){$-$}\rput*(c33){$+$}\rput*(c34){$+$}\rput*(b45){$-$}
\rput*(c10){$+$}\rput*(b21){$-$}\rput*(c21){$+$}\rput*(b32){$-$}
\rput*(c32){$+$}\rput*(c42){$+$}\rput*(b53){$-$}\rput*(b54){$-$}
\psline[linecolor=darkred]
    (-.1,4.3)(-.1,2.9)(.9,2.9)(.9,1.9)(1.9,1.9)(1.9,.9)(3.9,.9)(3.9,-.1)
\psline[linecolor=darkred]
    (-.3,4.1)(2.1,4.1)(2.1,.1)(4.1,.1)
\uput[225](-.1,2.9){$p'$}
\uput[45](2.1,4.1){$p$}
\endpspicture
$$
If an H-web is placed in each square as indicated, the two paths $p$
and $p'$
delineate a sequence of connected webs, separated by points or
segments where $p$ and $p'$ meet and may or may not cross.  We attach each
web bounded below by $p$ to $w$, and we invert and reverse the arrows
of each web bounded above by $p$ and then attach it to $w$.  The
result is a web $w' \in B(C[s'])$ with core $w$:
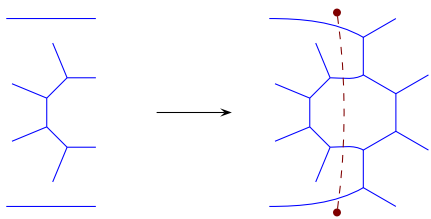
$$
\pspicture[.5](-1,-1.25)(.5,1.25)
\qline(-1,1.25)(.191,1.25)\qline(-1,-1.25)(.191,-1.25)
\qline(.5;112.5)(1;112.5)\qline(.5;157.5)(1;157.5)
\qline(.5;202.5)(1;202.5)\qline(.5;247.5)(1;247.5)
\psline(.5;67.5)(.5;112.5)(.5;157.5)(.5;202.5)(.5;247.5)(.5;292.5)
\endpspicture
\psgoesto
\pspicture[.5](-1,-1.4)(1,1.4)
\qline(.5;112.5)(1;112.5)
\qline(.5;157.5)(1;157.5)
\qline(.5;202.5)(1;202.5)
\qline(.5;247.5)(1;247.5)
\psline(.5;112.5)(.5;157.5)(.5;202.5)(.5;247.5)
\pccurve[angleA=0,angleB=210](.5;112.5)(.25,.5)
\pccurve[angleA=0,angleB=150](.5;247.5)(.25,-.5)
\pccurve[angleA=0,angleB=150](-1,1.25)(.25,1)
\pccurve[angleA=0,angleB=210](-1,-1.25)(.25,-1)
\psline(.683,1.25)(.25,1)(.25,.5)(.683,.25)(.683,-.25)(.25,-.5)
    (.25,-1)(.683,-1.25)
\qline(.683,.25)(1.116,.5)\qline(.683,-.25)(1.116,-.5)
\pcarc[arrows=*-*,linestyle=dashed,linecolor=darkred](-.1,1.33)(-.1,-1.33)
\endpspicture
$$

If $s_1$ and $s_2$ are arbitrary strings of the
same weight, we can set $s'$ to
either in the stair-step construction, which implies that any core $w$
of a web in $B(C[s_1])$ is also a core of a web in $B(C[s_2])$.
Moreover, it is easy to check that, if $w$ is fixed, the
collection of webs so produced is closed under H-moves.
Since every basis web is related to its core by H-moves,
the stair-step construction produces all webs in $B(Cc)$ for
all $c$ with weight $\lambda$.  In particular,
we have demonstrated the following lemma:

\begin{lemma} Let $C$ be a sequence of clasps and let
$c$ be a clasp of weight $\lambda$.  If $w \in B(C[c])$,
let $f(w)$ be its core.  Then the map $f:B(C[c]) \to B(C;\lambda)$
is a bijection. \label{lcorebijec}
\end{lemma}

Using the above lemmas, we can prove Theorem~\ref{tha2splitting}:

\begin{proof} Let $w_1 \in B(C;\lambda)$ and $w_2 \in B(D;\lambda^*)$,
and let $c_1$ and $c_2$ be the clasps of $w_1$ and $w_2$ of
weight $\lambda$ and $\lambda^*$.  Then we can either extend $w_1$
to a basis web in $B(Cc_2^*)$ or extend $w_2$ to a basis web
in $B(Dc_1^*)$ and then sew the two webs together;
by the symmetry of the stair-step construction,
the resulting web $w \in B(CD)$ is the same in both cases.
Define
$$f:\bigcup_\lambda
\left(B(C;\lambda) \times B(D;\lambda^*) \right)\to B(CD)$$
by the above operation on pairs of cores.  By
Lemma~\ref{lcomparable} and the fact that H-moves preserve
cores, there is also a map
$$g:B(CD) \to \bigcup_\lambda \left(B(C;\lambda)
\times B(D;\lambda^*)\right)$$
given by splitting $w$ along a minimal cut path and taking
the cores of the two halves.  By Lemma~\ref{lcorebijec},
$f$ and $g$ are inverses.  Using Lemma~\ref{lcorebijec} again,
we can see $g$ as the bijection claimed by the theorem.
\end{proof}

To prove Theorem~\ref{tha2equinum}, we need one additional lemma.
Say that a sign string $s$ is segregated if it is of the form
$$+ + \ldots + - - \ldots -;$$
a sign string with only $+$'s or only $-$'s is automatically segregated.
Likewise, say that a clasp is segregated if it is of the form $[s]$ for a
segregated sign string $s$. 

\begin{lemma}  Let $c$ and $d$ be segregated clasps of weight
$\lambda$ and $\mu$.  Then the web basis set $B([+]cd)$ has
has one element if and only if $\lambda^* = \mu + \lambda_1$
$\mu^* = \lambda + \lambda_1$, or $\lambda^* = \mu + \lambda_2 - \lambda_1$,
and is empty otherwise. \label{lthreeclasps}
\end{lemma}
\begin{proof}
Consider first three segregated clasps $c$, $d$, and $e$ (one of which
might be empty) of
arbitrary weight, and let $w \in B(cde)$.  We claim that $w$
must consist of a number of bare strands plus a flat component,
as in the following example:
$$
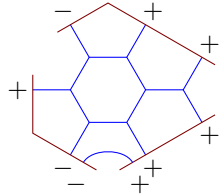
\pspicture(-2,-1.5)(2,1.5)
\pnode(.5;  0){a1}\pnode(1;  0){b1}
\pnode(.5; 60){a2}\pnode(1; 60){b2}
\pnode(.5;120){a3}\pnode(1;120){b3}
\pnode(.5;180){a4}\pnode(1;180){b4}
\pnode(.5;240){a5}\pnode(1;240){b5}
\pnode(.5;300){a6}\pnode(1;300){b6}
\pnode([nodesep=.5,angle=60]b1){c1}
\pnode([nodesep=.5,angle=300]b1){c2}
\pnode([nodesep=.2,angle=330]b5){d1}
\pnode([nodesep=.2,angle=210]b6){d2}
\psline(c2)(b1)(c1)
\pspolygon(a1)(a2)(a3)(a4)(a5)(a6)
\qline(a1)(b1)\qline(a2)(b2)\qline(a3)(b3)
\qline(a4)(b4)\qline(a5)(b5)\qline(a6)(b6)
\psline[linecolor=darkred]([nodesep=.2,angle=90]b4)
	([nodesep=.577,angle=270]b4)([nodesep=.4,angle=330]b5)
\psline[linecolor=darkred]([nodesep=.4,angle=210]b6)
	([nodesep=.2,angle=30]c2)
\psline[linecolor=darkred]([nodesep=.2,angle=330]c1)
	([nodesep=.577,angle=150]b2)([nodesep=.2,angle=210]b3)
\nccurve[angleA=60,angleB=120]{d1}{d2}
\rput([nodesep=.2,angle=60]c1){$+$}
\rput([nodesep=.2,angle=60]b2){$+$}
\rput([nodesep=.2,angle=120]b3){$-$}
\rput([nodesep=.2,angle=180]b4){$+$}
\rput([nodesep=.2,angle=240]b5){$-$}
\rput[t]([nodesep=.2,angle=240]d1){$-$}
\rput[t]([nodesep=.2,angle=300]d2){$+$}
\rput([nodesep=.2,angle=300]b6){$+$}
\rput([nodesep=.2,angle=300]c2){$+$}
\endpspicture
$$
The argument is similar to the proof of Lemma~\ref{lcomparable}.
Assume first that $w$ is not connected.
If $p$ and $q$ are adjacent endpoints of $w$, we define the curvature
of the segment of the boundary of $w$ from $p$ to $q$ as in
Lemma~\ref{lcomparable}.
Since $w$ is non-elliptic, it must
one the one hand have total exterior curvature at most $360^\circ$.
On the other hand, there can be no ``U'' or ``Y'' attached along $c$,
$d$, or $e$, and in there is only one place along each clasp where
an $H$ is possible.  The curvature at this $H$, if it is present,
is at most $60^\circ$, and the curvature elsewhere along the clasps
is at most $0^\circ$.  Moreover, unless $w$ is a bare strand,
the curvature at the segments between the clasps is at most
$60^\circ$ also.  The largest possible total is $360^\circ$ exactly,
so that $w$ is flat and its boundary is qualitatively a hexagon:
$$
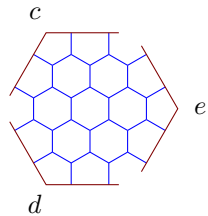
\pspicture(-1.5,-1.5)(1.5,1.5)
\rput(-.5,.5774){\littley}\rput(0,.5774){\littley}\rput(.5,.5774){\littley}
\rput(-.75,.1443){\littley}\rput(-.25,.1443){\littley}
\rput(.25,.1443){\littley}\rput(.75,.1443){\littley}
\rput(-.5,-.2887){\littley}\rput(0,-.2887){\littley}\rput(.5,-.2887){\littley}
\rput(-.25,-.7217){\littley}\rput(.25,-.7217){\littley}
\qline(-.25,1.0104)(-.25,.7217)\qline(.25,1.0104)(.25,.7217)
\rput{120}{\qline(-.25,1.0104)(-.25,.7217)\qline(.25,1.0104)(.25,.7217)}
\rput{240}{\qline(-.25,1.0104)(-.25,.7217)\qline(.25,1.0104)(.25,.7217)}
\pnode(1.167;0){c1}\pnode(1.167;60){c2}\pnode(1.167;120){c3}
\pnode(1.167;180){c4}\pnode(1.167;240){c5}\pnode(1.167;300){c6}
\psline[linecolor=darkred]([nodesep=.2,angle=180]c2)(c3)([nodesep=.2,angle=60]c4)
\psline[linecolor=darkred]([nodesep=.2,angle=300]c4)(c5)([nodesep=.2,angle=180]c6)
\psline[linecolor=darkred]([nodesep=.2,angle=60]c6)(c1)([nodesep=.2,angle=300]c2)
\rput([nodesep=.3,angle=120]c3){$c$}
\rput([nodesep=.3,angle=240]c5){$d$}
\rput([nodesep=.3,angle=0]c1){$e$}
\endpspicture
$$
If $w$ is not connected, then the above argument applies to each
connected component of $w$, with the additional observation
that at most one component of $w$ can meet all three clasps $c$,
$d$, and $e$.  A component that only meets two clasps is necessarily
a bare strand.

In the case of interest, $e = [+]$.  In this case, $w$ might
consist entirely of bare strands, or it might have a component
in the shape of a parallelogram or a trapezoid:
$$
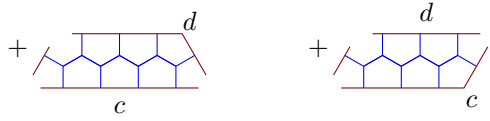
\pspicture(-1.5,-.75)(1.5,1)
\rput(-.5,0){\littley}\rput(0,0){\littley}\rput(.5,0){\littley}
\rput(1,0){\littley}
\rput(.75,.1443){\littlelam}
\rput(.25,.1443){\littlelam}
\rput(-.25,.1443){\littlelam}
\pnode(1.083,.433){c1}\pnode(-.583,.433){c2}
\pnode(-1,-.2887){c3}\pnode(1.5,-.2887){c4}
\ncline[linecolor=darkred,nodesep=.2]{c3}{c4}\Bput{$c$}
\ncline[linecolor=darkred,nodesep=.2]{c3}{c2}\Aput{$+$}
\psline[linecolor=darkred]
	([nodesep=.2,angle=0]c2)(c1)([nodesep=.2,angle=120]c4)
\rput([nodesep=.2,angle=60]c1){$d$}
\endpspicture
\hspace{1cm}
\pspicture(-1.5,-.75)(1.5,1)
\rput(-.5,0){\littley}\rput(0,0){\littley}\rput(.5,0){\littley}
\rput(.75,.1443){\littlelam}
\rput(.25,.1443){\littlelam}
\rput(-.25,.1443){\littlelam}
\pnode(1.25,.433){c1}\pnode(-.583,.433){c2}
\pnode(-1,-.2887){c3}\pnode(.833,-.2887){c4}
\ncline[linecolor=darkred,nodesep=.2]{c2}{c1}\Aput{$d$}
\ncline[linecolor=darkred,nodesep=.2]{c3}{c2}\Aput{$+$}
\psline[linecolor=darkred]
	([nodesep=.2,angle=0]c3)(c4)([nodesep=.2,angle=240]c1)
\rput([nodesep=.2,angle=300]c4){$c$}
\endpspicture
$$
These possibilities exactly match the restrictions on the
weights of $c$ and $d$.
\end{proof}

Finally, we prove Theorem~\ref{tha2equinum}:

\begin{proof} Let $\Lambda_+$ be the set of all weights $\lambda = a\lambda_1
+ b\lambda_2$ with $a,b \in \Z_{\ge 0}$. Consider the abelian group
$\Z[\Lambda_+]$ of formal sums of elements $v(\lambda)$ for each $\lambda \in
\Lambda_+$.  Then we can define a product in $\Z[\Lambda_+]$ using the
dimensions of the clasped web spaces:
$$v(\lambda_1)v(\lambda_2) = \sum_{\lambda_3} |B(c_1c_2c_3^*)|v(\lambda_3),$$
where in the sum each $c_i$ is some clasp of weight $\lambda_i$.
By Theorem~\ref{tha2splitting}, there exist bijections
$$\bigcup_d \left(B(c_1c_2d) \times B(c_3c_4d^*) \right)\cong B(c_1c_2c_3c_4)
\cong \left(\bigcup_d B(c_2c_3d) \times B(c_4c_1d^*)\right).$$
These bijections imply that
multiplication is associative, and therefore $\Z[\Lambda_+]$
is a ring.  To establish Theorem~\ref{tha2equinum}, it suffices
to check that the map $v(\lambda) \mapsto V(\lambda)$ induces
an isomorphism from $\Z[\Lambda_+]$ to the Grothendieck ring of
$A_2$.  Using induction, it suffices to check that for all
$\lambda \in \Lambda_+$,
\begin{eqnarray*}
v(0)v(\lambda) & = & v(\lambda) \\
v(\lambda_1)v(\lambda) & = & v(\lambda+\lambda_1) +
v(\lambda-\lambda_1+\lambda_2) + v(\lambda-\lambda_2) \\
v(\lambda_2)v(\lambda) & = & v(\lambda+\lambda_2) +
v(\lambda-\lambda_2+\lambda_1) + v(\lambda-\lambda_1),
\end{eqnarray*}
where $v(\lambda)$ is defined as 0 for $\lambda \notin \Lambda_+$,
because similar relations hold in the Grothendieck ring.
(Note that we cannot argue from an {\em a priori} hypothesis that
$\Z[\Lambda_+]$ is the Grothendieck ring of a category,
because such a construction assumes the existence of internal clasps,
which in turn depends on Theorem~\ref{tha2equinum}.)
In terms of webs, the first relation states that if $c$ and $d$
are two segregated clasps, then $B(cd)$ is empty unless $c = d^*$,
in which case it has one element.  This follows from Lemma~\ref{lcomparable}
or from arguments similar to those of Lemma~\ref{lthreeclasps}.  Similarly,
the other relations are equivalent to Lemma~\ref{lthreeclasps}.
\end{proof}

\subsection{The $B_2$ and $G_2$ cases}
\label{sb2g2}

\begin{theorem} Let $C$ be a $B_2$ clasp sequence and
let $\wt C$ be the corresponding sequence of weights.
Then the vector spaces $W(C)$ and $\Inv(V(\wt C))$ have
the same dimension. \label{thb2equinum}
\end{theorem}

\begin{theorem} Let $C$ be an $G_2$ clasp sequence and
let $\wt C$ be the corresponding sequence of weights.
Then the vector spaces $W(C)$ and $\Inv(V(\wt C))$ have
the same dimension. \label{thg2equinum}
\end{theorem}

To establish these two results we mainly need to alter various technical
definitions given for the $A_2$ case; most of the lemmas leading up to
Theorem~\ref{tha2splitting} and the proof of the theorem then carry over
word-for-word.  An H-web in the $B_2$ or $G_2$ spider is the web:
$$\doublehhoriz$$
In both the $B_2$ and $G_2$ spiders, one can define an H-move
on a basis web $w$ with two adjacent endpoints labelled 1 and 2.
An H-move at two such endpoints consists of attaching
an H-web to $w$, provided that this operation does not create
an elliptic face, and then replacing an internal double edge
by a tetravalent vertex in the $B_2$ spider
or by a perpendicular single edge in the $G_2$ spider:
$$ \btwohhoriz \psgoesto \btwox \hspace{2cm} \mbox{($B_2$ case)} 
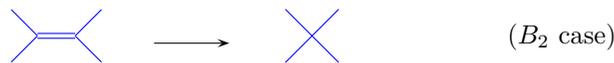
$$
\begin{equation}
\btwohhoriz \psgoesto \gtwohvert \hspace{2cm} \mbox{($G_2$ case)} \label{fcontract}
\end{equation}
If attaching an H-web would result in an elliptic face, then
there are two alternatives:  Either an H-web can be removed
(possibly after introducing an internal type 2 strand by the above
operations), or the two endpoints make a ``Y'':
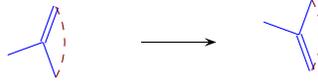
$$
\pspicture[.5](-.6,-.6)(.8,.6)
\psline(.5;290)(0,0)(.5;200)
\psline[doubleline=true](.5;70)(0,0)
\pcarc[linecolor=darkred,linestyle=dashed,arcangle=30](.5;70)(.5;290)
\endpspicture
\psgoesto
\pspicture[.5](-.6,-.8)(.6,.6)
\psline(.5;70)(0,0)(.5;160)
\psline[doubleline=true](.5;290)(0,0)
\pcarc[linecolor=darkred,linestyle=dashed,arcangle=30](.5;70)(.5;290)
\endpspicture
$$
or, in the $G_2$ case, they make an ``H'' with only one type
2 strand:
$$
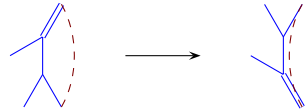
\pspicture[.5](-.6,-.6)(.6,1.1)
\psline(.25,-.433)(0,0)(0,.5)(-.433,.25)
\qline(0,0)(-.25,-.433)
\psline[doubleline=true](0,.5)(.25,.933)
\pcarc[linecolor=darkred,linestyle=dashed,arcangle=30](.25,.933)(.25,-.433)
\endpspicture
\psgoesto
\pspicture[.5](-.6,-1.1)(.6,.6)
\psline(.25,.433)(0,0)(0,-.5)(-.433,-.25)
\qline(0,0)(-.25,.433)
\psline[doubleline=true](0,-.5)(.25,-.933)
\pcarc[linecolor=darkred,linestyle=dashed,arcangle=30](.25,-.933)(.25,.433)
\endpspicture
$$
In each of these exceptional cases, the operation indicated is
the H-move. 

The operation of cutting along a minimal path is slightly more complicated
than before, because the path might cut diagonally across a vertex in the
$B_2$ spider or contain a type 1 edge in the $G_2$ spider.  The operation is
defined by first introducing a type 2 edge, as is sometimes also necessary
for an H-move, by reversing the contraction operation of
Figure~(\ref{fcontract}). A core is defined in the same way for all three rank
2 spiders. The stair-step construction is essentially the same as before.
Each square is an H-web, with the result that sewing together two cores
with stair steps in between results in many internal type 2 strands, which
are removed by the operation of Figure~(\ref{fcontract}).

The biggest difference between the three rank 2 cases is in the
statement and proof of the analogues of Lemma~\ref{lthreeclasps}.
We wish to check that, if $s \in \{1,2\}$ and $c$ and $d$
are arbitrary clasps of weight $\lambda$ and
$\mu$, then
\begin{equation}
|B([s]cd)| = \dim \Inv(V_s \tensor V(\lambda) \tensor V(\mu)).
\label{eb2g2rels}
\end{equation}
For this purpose, one would like to argue that an element of 
$B([s]cd)$ has no hyperbolic faces.  As before, it is convenient
to consider segregated clasps, where here a clasp is segregated
if it is of the form $11\ldots122\ldots2$, as well as reverse
segregated clasps of the form $22\ldots211\ldots1$.  In the $B_2$
case, it is easy to enumerate the elements of $B([1]cd)$ if
$c$ is segregated and $d$ is reverse segregated, as well
as the elements of $B([2]cd)$ if $c$ is reverse segregated and 
$d$ is segregated.  A web $w$ in such a basis web set cannot
have negative curvature, and after removing bare strands, it
is either empty or its the shape of its boundary is one of the
following (possibly with $c$ and $d$ switched):
$$
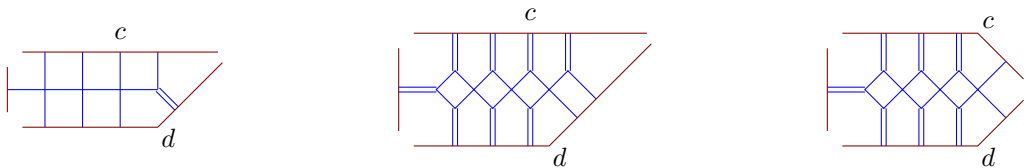
\pspicture[.5](-.1,-.6)(3.1,.6)
\psline(0,0)(2,0)(2,.5)
\psline[doubleline=true](2,0)(2.25,-.25)
\qline(.5,-.5)(.5,.5)
\qline(1,-.5)(1,.5)
\qline(1.5,-.5)(1.5,.5)
\pnode(2,-.5){c1}
\pnode(3,.5){c2}
\psline[linecolor=darkred](0,-.3)(0,.3)
\psline[linecolor=darkred](.2,-.5)(c1)([nodesep=.2,angle=225]c2)
\rput([nodesep=.2,angle=315]c1){$d$}
\pcline[linecolor=darkred](.2,.5)([nodesep=.2,angle=180]c2)\Aput{$c$}
\endpspicture
\hspace{2cm}
\pspicture[.5](-.1,-.5)(3.6,2)
\psline[doubleline=true](0,.75)(.5,.75)
\psline[doubleline=true](.75,0)(.75,.5)
\psline[doubleline=true](.75,1)(.75,1.5)
\psline[doubleline=true](1.25,0)(1.25,.5)
\psline[doubleline=true](1.25,1)(1.25,1.5)
\psline[doubleline=true](1.75,0)(1.75,.5)
\psline[doubleline=true](1.75,1)(1.75,1.5)
\psline[doubleline=true](2.25,1)(2.25,1.5)
\psline(2.375,.375)(1.75,1)(1.25,.5)(.75,1)(.5,.75)
	(.75,.5)(1.25,1)(1.75,.5)(2.25,1)(2.625,.625)
\pnode(2,0){c1}\pnode(3.5,1.5){c2}
\psline[linecolor=darkred](0,.2)(0,1.3)
\psline[linecolor=darkred](.2,0)(c1)([nodesep=.2,angle=225]c2)
\pcline[linecolor=darkred](.2,1.5)([nodesep=.2,angle=180]c2)\Aput{$c$}
\rput([nodesep=.2,angle=315]c1){$d$}
\endpspicture
\hspace{2cm}
\pspicture[.5](-.1,-.5)(3.6,2)
\psline[doubleline=true](0,.75)(.5,.75)
\psline[doubleline=true](.75,0)(.75,.5)
\psline[doubleline=true](.75,1)(.75,1.5)
\psline[doubleline=true](1.25,0)(1.25,.5)
\psline[doubleline=true](1.25,1)(1.25,1.5)
\psline[doubleline=true](1.75,0)(1.75,.5)
\psline[doubleline=true](1.75,1)(1.75,1.5)
\psline(2.375,.375)(1.75,1)(1.25,.5)(.75,1)(.5,.75)
     (.75,.5)(1.25,1)(1.75,.5)(2.375,1.125)
\pnode(2,0){c1}\pnode(2.75,.75){c2}\pnode(2,1.5){c3}
\psline[linecolor=darkred](0,.2)(0,1.3)
\psline[linecolor=darkred](.2,0)(c1)([nodesep=.2,angle=225]c2)
\psline[linecolor=darkred](.2,1.5)(c3)([nodesep=.2,angle=135]c2)
\rput([nodesep=.2,angle=315]c1){$d$}
\rput([nodesep=.2,angle=45]c3){$c$}
\endpspicture
$$
In the $G_2$ case, the computation is again simpler if one of $c$
and $d$ is segregated and the other is reverse segregated.
No webs in $B([1]cd)$ have any negative
curvature, but a web in $B([2]cd)$ might have one negatively
curved face, namely one with three type 2 strands and one type
1 strand incident to it:
$$
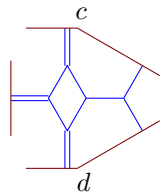
\pspicture(-.1,-1.5)(2.5,1.5)
\psline[doubleline=true](0,0)(.5,0)
\psline[doubleline=true](.75,.433)(.75,.933)
\psline[doubleline=true](.75,-.433)(.75,-.933)
\pspolygon(1,0)(.75,-.433)(.5,0)(.75,.433)
\qline(1,0)(1.5,0)\qline(1.5,0)(1.75,.433)\qline(1.5,0)(1.75,-.433)
\pnode(.884,-.933){c1}\pnode(2.5,0){c2}\pnode(.884,.933){c3}
\psline[linecolor=darkred](0,-.5)(0,.5)
\psline[linecolor=darkred](.2,-.933)(c1)([nodesep=.6,angle=210]c2)
\psline[linecolor=darkred](.2,.933)(c3)([nodesep=.6,angle=150]c2)
\rput([nodesep=.2,angle=295]c1){$d$}
\rput([nodesep=.2,angle=75]c3){$c$}
\endpspicture
$$
All other $G_2$ webs in $B([2]cd)$ are flat.  An exhaustive enumeration
of such webs in both spiders verifies equation~(\ref{eb2g2rels}).
The details are left to the reader.

Finally, in the $G_2$ spider, unlike in the $A_2$ and $B_2$
spiders, web with a partly convex clasp does not
necessarily vanish in the clasped web space $W([s_1][s_2]\ldots[s_n])$.
Recall that the clasped web space is defined as the subspace
of the unclasped web space $W(s_1s_2\ldots s_n)$ spanned by
non-elliptic webs with non-convex clasps and no internal type 2
edges, and that there is also a kernel $I([s_1][s_2]\ldots [s_n])$
spanned by webs $w$ such that there is a cut path transverse to
all edges of $w$ which cuts off a clasp of lower weight.
In order to show that the $G_2$ spider is well-defined, we need
the following lemma:

\begin{lemma} The spaces $I([s_1][s_2]\ldots [s_n])$ and
$W([s_1][s_2]\ldots [s_n])$ are transverse in $W(s_1s_2\ldots s_n)$.
\end{lemma}
\begin{proof}
Rank the webs in $W(s_1s_2\ldots s_n)$ first by the number 
of vertices, and second by the sum of the weights of minimal cut
paths cutting off each string $s_i$.  We consider a change of basis
which is lower-triangular with respect to this ranking:  For each
web $w \in W(s_1s_2\ldots s_n)$, choose a minimal cut path $p_i$ cutting
off each $s_i$ as close as possible to the boundary.  The cut
paths are unique by the $G_2$ analogue of Lemma~\ref{lecore},
and they cannot cross, although they could in principle share
segments.  Recall the equation for a double edge from the
equations~(\ref{fg2elliptic}), rearranged slightly:
$$
\gtwohvert
=
- 
\pspicture[.4](-.9,-.5)(.9,.5)
\pnode(-.25,0){a1}\pnode(.25,0){a2}
\qline([nodesep=.5,angle=150]a1)(a1)
\qline([nodesep=.5,angle=210]a1)(a1)
\qline([nodesep=.5,angle=30]a2)(a2)
\qline([nodesep=.5,angle=330]a2)(a2)
\psline[doubleline=true](a1)(a2)
\endpspicture
+
\hh
- {1 \over q^2 - 1 + q^{-2}}\vv 
+ {1 \over q + 1 + q^{-1}}\gtwohhoriz
$$
Applying this equation, we can convert every single edge contained
in a path $p_i$ to a transverse double edge; the other terms
are all lower with respect to the ranking of webs.  This produces
a new basis for $W(s_1s_2\ldots s_n)$ in which $I([s_1][s_2]\ldots [s_n])$
and $W([s_1][s_2]\ldots [s_n])$ are manifestly complements.
\end{proof}

\section{Explicit formulas for $A_2$ clasps}
\label{sexplicit}

Although internal clasps must exist by Section~\ref{smorphism}, the argument
given there is too indirect for practical computations.  In this section, we
give an explicit formulas for $A_2$ clasps.  Note first that clasps are sent
to each other by composition with $H$
webs:
$$
\pspicture[.42](-.5,-1.4)(2.5,1.4)
\psframe[linecolor=darkred](-.2,-.1)(1.7,.1)
\qline(0,.1)(0,.966)\qline(.5,.1)(.5,.966)
\qline(1,.1)(1,.966)\qline(1.5,.1)(1.5,.966)
\qline(0,-.1)(0,-.966)\qline(.5,-.1)(.5,-.966)
\qline(1,-.1)(1,-.966)\qline(1.5,-.1)(1.5,-.966)
\uput{3pt}[90]( 0,.966){$+$}\uput{3pt}[90](.5,.966){$-$}
\uput{3pt}[90](1,.966){$+$}\uput{3pt}[90](1.5,.966){$-$}
\uput{3pt}[270](0,-.966){$-$}\uput{3pt}[270](.5,-.966){$+$}
\uput{3pt}[270](1,-.966){$-$}\uput{3pt}[270](1.5,-.966){$+$}
\endpspicture
=
\pspicture[.42](-.5,-1.4)(3,1.4)
\psframe[linecolor=darkred](-.2,-.1)(2.2,.1)
\qline(0,.1)(0,.966)\qline(2,.1)(2,.966)
\psline(.5,.1)(.75,.533)(.5,.966)
\psline(1.5,.1)(1.25,.533)(1.5,.966)
\qline(.75,.533)(1.25,.533)
\qline(0,-.1)(0,-.966)\qline(2,-.1)(2,-.966)
\psline(.5,-.1)(.75,-.533)(.5,-.966)
\psline(1.5,-.1)(1.25,-.533)(1.5,-.966)
\qline(.75,-.533)(1.25,-.533)
\uput{3pt}[90](0,.966){$+$}\uput{3pt}[120](.5,.966){$-$}
\uput{3pt}[60](1.5,.966){$+$}\uput{3pt}[90](2,.966){$-$}
\uput{3pt}[270](0,-.966){$-$}\uput{3pt}[240](.5,-.966){$+$}
\uput{3pt}[300](1.5,-.966){$-$}\uput{3pt}[270](2,-.966){$+$}
\endpspicture
$$

Thus, it suffices to derive a formula for segregated clasps.

\begin{lemma} If $s$ is a segregated sign string and $t^*$ is another sign
string of lower or incomparable weight, then a basis web $w \in B(st)$ has
either a ``U'' or a ``Y'' attached to $s$. \label{lsegregated}
\end{lemma}
\begin{proof}
The proof is also similar to that of Lemma~\ref{lcomparable}. We assume the
minimal counterexample.  First, the web $w$ must be connected, for if a
connected component meets $t$ but not $s$, it may be discarded to produce a
smaller counterexample; if it meets $s$ but not $t$, then it is a smaller
counterexample; and if all components meet both $s$ and $t$, then one of them
must be a smaller counterexample.  Second, there can be no ``Y'' or ``H''
attached at $t$, for otherwise they may be discarded to produce a smaller
counterexample.  (Note that a ``U'' attached to $t$ has already been eliminated.)

Since $w$ is non-elliptic, the total
exterior curvature is at least $360^\circ$.  Since there 
is no ``U'', ``Y'', or ``H'' attached at $t^*$, the total
curvature along $t$ is at most $0^\circ$.  Moreover,
unless $w$ is a bare strand (which is not a counterexample),
the two arcs connecting $s$ to $t$ have curvature at most
$120^\circ$.  By the hypothesis that $s$ has no ``Y'',
the curvature in each segregated segment of $S$ is at most $0^\circ$
also.  Finally, since $s$ has no ``U'', the curvature at the single
arc connecting opposite signs of $s$ is at most $60^\circ$, for
a total of at most $300^\circ$, a contradiction:
$$
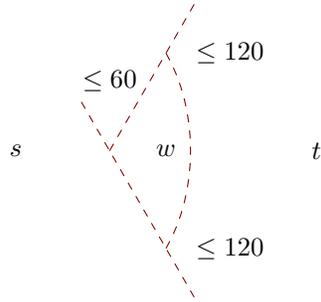
\pspicture(-2.5,-2)(2.5,2)
\pnode(0,1.299){a1}\pnode(-.75,0){a2}
\pnode(0,-1.299){a3}
\psline[linestyle=dashed,linecolor=darkred]([nodesep=.75,angle=300]a3)
     ([nodesep=.75,angle=120]a2)
\psline[linestyle=dashed,linecolor=darkred](a2)([nodesep=.75,angle=60]a1)
\pcarc[linestyle=dashed,linecolor=darkred,arcangle=30](a1)(a3)
\rput(-2,0){$s$}
\rput(2,0){$t$}
\rput(0,0){$w$}
\uput{.4}[0](a1){$\scriptsize \le 120$}
\uput{.4}[0](a3){$\scriptsize\le 120$}
\uput{.8}[90](a2){$\scriptsize \le 60$}
\endpspicture
$$
\end{proof}

By Lemma~\ref{lsegregated}, if a web $w \in W(ss^*)$
annihilates any ``U'' or ``Y'', it satisfies the annihilation axiom
of a clasp.  Moreover, among terms
in $w$, only $\beta_s$ (the web consisting solely of
parallel strands) has the weight of $s$ as its cut weight.
Therefore all other terms are annihilated by such $w$, and the
annihilation property implies that $w$ is an idempotent under
concatenation if the coefficient of the leading term $\beta_s$
is 1.  Given these facts, a routine argument by induction
(see Wenzl \cite{Wenzl:canada1987} and Ohtsuki and Yamada \cite{OY:quantum})
establishes that
$$
\pspicture[.4](-.1,-.3)(1.1,1.3)
\rput[tr](.7,-.1){$n\,-$}\qline(.5,0)(.5,.4)
\psframe[linecolor=darkred](0,.4)(1,.6)
\qline(.5,.6)(.5,1)\rput[br](.7,1.1){$n\,+$}
\endpspicture
=
\pspicture[.4](-.1,-.3)(1.3,1.3)
\rput[tr](.7,-0.1){$n-1\,-$}\qline(.5,0)(.5,.4)
\psframe[linecolor=darkred](0,.4)(1,.6)
\qline(.5,.6)(.5,1)\rput[br](.7,1.1){$n-1\,+$}
\qline(1.3,0)(1.3,1)
\endpspicture
+
{[n-1] \over [n]}
\pspicture[.4](-1,-.3)(1.8,2.3)
\rput[tr](.7,-.1){$n-1\,-$}
\qline(.5,0)(.5,.4)
\psframe[linecolor=darkred](0,.4)(1,.6)
\pcline(.25,.6)(.25,1.4)\middlearrow\rput[r](.15,1){$n-2$}
\psframe[linecolor=darkred](0,1.4)(1,1.6)
\qline(.5,1.6)(.5,2)
\rput[br](.7,2.1){$n-1\,+$}
\psline(.75,.6)(1.096,.8)(1.096,1.2)(.75,1.4)
\pccurve[angleA=330,angleB=90](1.096,.8)(1.443,0)
\pccurve[angleA=30,angleB=270](1.096,1.2)(1.443,2)
\rput[t](1.443,-.1){$-$}
\rput[b](1.443,2.1){$+$}
\endpspicture
$$
recursively defines a clasp of weight $n\lambda_1$, and a
computation shows that
$$
\pspicture[.4](-.5,-.3)(1.1,1.3)
\rput[tr](.2,-.1){$a\,-$}\qline(0,0)(0,.4)
\rput[tr](.9,-.1){$b\,+$}\qline(.7,0)(.7,.4)
\psframe[linecolor=darkred](-.2,.4)(.9,.6)
\qline(0,.6)(0,1)\qline(.7,.6)(.7,1)
\rput[br](.2,1.1){$a\,+$}
\rput[br](.9,1.1){$b\,-$}
\endpspicture
= 
\sum_{k=0}^{\min(a,b)}
(-1)^k{[a]![b]![a+b+k+1]!\over [a-k]![b-k]![a+b+1]![k]!}
\pspicture[.45](-2.7,-2)(2.7,2)
\psframe[linecolor=darkred](-1.5,-1.1)(-.5,-.9)
\psframe[linecolor=darkred]( 1.5,-1.1)( .5,-.9)
\psframe[linecolor=darkred](-1.5, 1.1)(-.5, .9)
\psframe[linecolor=darkred]( 1.5, 1.1)( .5, .9)
\pcline(-1.3,-.9)(-1.3,.9)\middlearrow\Aput{$a-k$}
\pcline(1.3,.9)(1.3,-.9)\middlearrow\Aput{$b-k$}
\pccurve[angleA=270,angleB=270,ncurv=1](.7,.9)(-.7,.9)
\middlearrow
\pccurve[angleA=90,angleB=90,ncurv=1](-.7,-.9)(.7,-.9)
\middlearrow\rput(0,0){$k$}
\qline(-1,-1.5)(-1,-1.1)\qline(-1,1.1)(-1,1.5)
\qline( 1,-1.5)( 1,-1.1)\qline( 1,1.1)( 1,1.5)
\rput[br](-.8,1.7){$a\,+$}\rput[br](1.2,1.7){$b\,-$}
\rput[tr](-.8,-1.7){$a\,-$}\rput[tr](1.2,-1.7){$b\,+$}
\endpspicture
$$
defines a segregated clasp of weight $a\lambda_1+b\lambda_2$.

\section{Applications and problems}
\label{slast}

\subsection{Higher rank}

The main open problem related to the combinatorial rank 2 spiders is how to
generalize them to higher rank.  A proper generalization would consist of a
complete set of generators and relations for the higher-rank spiders; the
strand sets would correspond to the fundamental representations and their
tensor products.  It is easy to make a HOMFLY spider which corresponds to the
HOMFLY polynomial, but this spider describes the invariant theory of $A_n$
only in the stable limit of large $n$; the $A_n$ web spaces for any fixed $n$
are quotients of the HOMFLY web spaces.   (Recently, Murakami, Ohtsuki, and
Yamada \cite{MOY:homfly} have defined the HOMFLY spider in terms of trivalent
graphs; this is a step toward an explicit description of the unstable
truncation.)  The generalization to a higher-rank Lie algebra $\frak g$, if
it exists, would also likely involve formal angles related to the Coxeter
geometry of $\frak g$; note that in both the $A_2$ and $B_2$ spiders, a
trivalent vertex is dual to a Weyl alcove.  Note also that in all three rank
2 spiders, a large, flat basis web (one with neither elliptic nor hyperbolic
faces) coincides with the Voronoi tiling of the plane given by the weight
lattice of the corresponding Lie algebra.

The bases given by the combinatorial rank 1 spider are dual to Lusztig's
canonical bases \cite{FK:canonical}.  Those of the rank 2 spiders are almost
certainly also dual to canonical bases or are closely related, because
canonical bases have the same symmetry of cyclic permutation of tensor
factors, and because of the integrality and positivity properties of the
coefficients in equations~(\ref{fa2elliptic}), (\ref{fb2elliptic}), and
(\ref{fg2elliptic}).

\subsection{Generalized $6j$ symbols}

Given four webs $w_1$, $w_2$, $w_3$, and $w_4$ and six clasps
$c_{12}$, $c_{13}$, $c_{14}$, $c_{23}$, $c_{24}$, and $c_{34}$ in some
spider, their
tetrahedron symbol is defined as the value of the compound web
$$
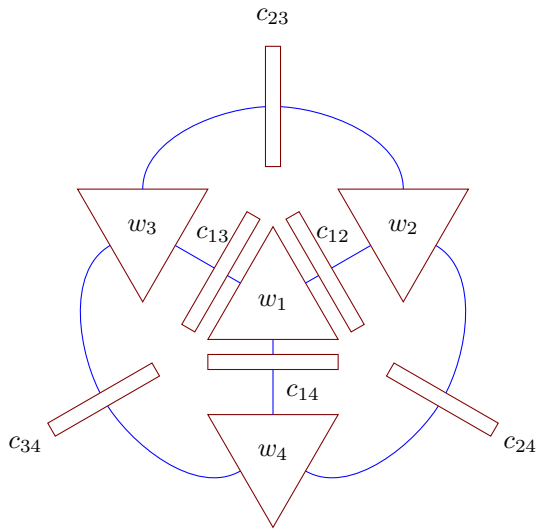
\pspicture(-4.5,-3.5)(4.5,4.5)
\pnode(1; 90){a1}\pnode(1.8; 90){b1}\pnode(2.6; 90){c1}\pnode(3.4; 90){d1}
\pnode(1;210){a2}\pnode(1.8;210){b2}\pnode(2.6;210){c2}\pnode(3.4;210){d2}
\pnode(1;330){a3}\pnode(1.8;330){b3}\pnode(2.6;330){c3}\pnode(3.4;330){d3}
\pnode(2.0;30){w2}\pnode(2.0;150){w3}\pnode(2.0;270){w4}
\pspolygon[linecolor=darkred](a1)(a2)(a3)
\pspolygon[linecolor=darkred]
    ([nodesep=1,angle=30]a1)([nodesep=1,angle=30]a3)(3.0;30)
\pspolygon[linecolor=darkred]
    ([nodesep=1,angle=150]a1)([nodesep=1,angle=150]a2)(3.0;150)
\pspolygon[linecolor=darkred]
    ([nodesep=1,angle=270]a2)([nodesep=1,angle=270]a3)(3.0;270)
\pspolygon[linecolor=darkred]
    ([nodesep=.2,angle=30]a1)([nodesep=.4,angle=30]a1)
    ([nodesep=.4,angle=30]a3)([nodesep=.2,angle=30]a3)
\pspolygon[linecolor=darkred]
    ([nodesep=.2,angle=150]a2)([nodesep=.4,angle=150]a2)
    ([nodesep=.4,angle=150]a1)([nodesep=.2,angle=150]a1)
\pspolygon[linecolor=darkred]
    ([nodesep=.2,angle=270]a3)([nodesep=.4,angle=270]a3)
    ([nodesep=.4,angle=270]a2)([nodesep=.2,angle=270]a2)
\pccurve[angleA=  0,angleB= 90](c1)([nodesep=.5,angle= 90]w2)
\pccurve[angleA=180,angleB= 90](c1)([nodesep=.5,angle= 90]w3)
\pccurve[angleA=120,angleB=210](c2)([nodesep=.5,angle=210]w3)
\pccurve[angleA=300,angleB=210](c2)([nodesep=.5,angle=210]w4)
\pccurve[angleA=240,angleB=330](c3)([nodesep=.5,angle=330]w4)
\pccurve[angleA= 60,angleB=330](c3)([nodesep=.5,angle=330]w2)
\pspolygon[linecolor=darkred,fillstyle=solid,fillcolor=white]
    ([nodesep=.1,angle=  0]b1)([nodesep=.1,angle=180]b1)
    ([nodesep=.1,angle=180]d1)([nodesep=.1,angle=  0]d1)
\pspolygon[linecolor=darkred,fillstyle=solid,fillcolor=white]
    ([nodesep=.1,angle=120]b2)([nodesep=.1,angle=300]b2)
    ([nodesep=.1,angle=300]d2)([nodesep=.1,angle=120]d2)
\pspolygon[linecolor=darkred,fillstyle=solid,fillcolor=white]
    ([nodesep=.1,angle=240]b3)([nodesep=.1,angle= 60]b3)
    ([nodesep=.1,angle= 60]d3)([nodesep=.1,angle=240]d3)
\rput(0,0){$w_1$}\rput(w2){$w_2$}\rput(w3){$w_3$}\rput(w4){$w_4$}
\qline(.5; 30)(.7; 30)\pcline(.9; 30)(1.5; 30)\Aput{$c_{12}$}
\qline(.5;150)(.7;150)\pcline(.9;150)(1.5;150)\Bput{$c_{13}$}
\qline(.5;270)(.7;270)\pcline(.9;270)(1.5;270)\Aput{$c_{14}$}
\rput(3.8; 90){$c_{23}$}\rput(3.8;210){$c_{34}$}\rput(3.8;330){$c_{24}$}
\endpspicture
$$
provided that the four webs $w_i$ are members of the appropriate clasped web
spaces.  The tetrahedron symbol is closely related to the $6j$ symbol,
which expresses the change of basis from
$$\bigoplus_{c_{14}} W(c_{12}c_{13}c_{14}) \tensor W(c_{43}c_{42}c_{41})$$
to 
$$\bigoplus_{c_{23}} W(c_{42}c_{12}c_{32}) \tensor W(c_{13}c_{43}c_{23})$$
via the identification of both with $W(c_{12}c_{13}c_{43}c_{42})$, where in
general $c_{ij} = c_{ji}^*$. Up to normalization, the $A_1$ tetrahedron
symbol at $q=1$ equals the Racah-Wigner $6j$ symbol used in mathematical
physics \cite{MS:cmp1989}. Using the $A_1$ spider, Masbaum and Vogel have
found a new proof of the Racah formula for the $6j$ symbol and its quantum
generalization \cite{MV:pjm1994}. The combinatorial rank 2 spiders could be
equally useful for understanding the rank 2 generalization of the $6j$ symbol.

\subsection{Practical computation of rank 2 link invariants}

If one is interested in computing link invariants,
equations~(\ref{fa2elliptic}), (\ref{fb2elliptic}), and
(\ref{fg2elliptic}) can be interpreted as inductive rules for
evaluating links and knotted graphs without boundary in the web space
$W(\emptyset)$. Mollard \cite{Mollard} and Sinha \cite{Sinha} have
independently written computer implementations of this algorithm for $G_2$,
and there are well-known computer programs to compute the Jones, HOMFLY, and
Kauffman polynomials using the same basic strategy.

Spiders suggest an alternative method for computing the same invariants which
is more efficient than a direct application of the above rules to closed
links and graphs.  The method consists of assembling a link projection, as a
web, from indiviual crossings using spider operations,  and reducing
intermediate webs to linear combinations of basis webs.  For example, to
evaluate the quantum $A_1$ link invariant (the Jones polynomial) of a figure
eight knot, we can decompose the knot as a nested sequence of three
tangles $a$, $b$, and $c$:
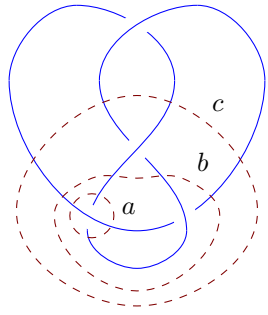
$$
\pspicture(-2.25,-.25)(2.25,4.25)
\pnode(-.8,4){a4}\pnode(.8,4){a8}\pnode(-1.7,3){a5}
\pnode(-.5,3){a9}\pnode(.5,3){a3}\pnode(1.7,3){a7}
\pnode(-.66,1){a2}\pnode(0,1){a6}\pnode(.66,1){a10}\pnode(0,.5){a1}
\nccurve[border=.15,angleA=0,angleB=270]{a6}{a7}
\nccurve[border=.15,angleA=270,angleB=90]{a9}{a10}
\nccurve[border=.15,angleA=270,angleB=0]{a10}{a1}
\nccurve[border=.15,angleA=180,angleB=270]{a1}{a2}
\nccurve[border=.15,angleA=90,angleB=270]{a2}{a3}
\nccurve[border=.15,angleA=270,angleB=180]{a5}{a6}
\nccurve[border=.15,angleA=90,angleB=0]{a3}{a4}
\nccurve[border=.15,angleA=180,angleB=90]{a4}{a5}
\nccurve[border=.15,angleA=90,angleB=0]{a7}{a8}
\nccurve[border=.15,angleA=180,angleB=90]{a8}{a9}
\pscircle[linestyle=dashed,linecolor=darkred](-.6,1.2){.3}
\psccurve[linestyle=dashed,linecolor=darkred](-1.1,1.2)(-.6,1.7)
    (0,1.7)(.6,1.7)(1.1,1.2)(0,.2)
\psccurve[linestyle=dashed,linecolor=darkred](-1.6,1.2)(0,2.8)(1.6,1.2)(0,0)
\rput[l](-.2,1.3){$a$}\rput[lb](.8,1.8){$b$}\rput[lb](1,2.6){$c$}
\endpspicture
$$
We recall some basic identities in the $A_1$ spider:
\begin{eqnarray*}
\singleloop & = & -q^{1/2} - q^{-1/2} \\
\rcrossing & = & - q^{1/4}\vv - q^{-1/4}\hh \\
\pspicture[.4](-.6,-.5)(.6,.5)
\pccurve[angleA=315,angleB=270,ncurv=1](0,0)(.5,0)
\qline(.5;135)(0,0)
\pccurve[border=.1,angleA=45,angleB=90,ncurv=1](0,0)(.5,0)
\psline[border=.1](.5;225)(0,0)
\endpspicture 
& = & q^{3/4}\pspicture[.4](-.6,-.5)(.6,.5)
\psbezier(.5;225)(.25;225)(.25;135)(.5;135)\endpspicture
\end{eqnarray*}
The tangle $a$ is a right-handed crossing.  Its expansion leads to an
expansion of the tangle $b$:
$$
\pspicture[.4](-1.1,-.5)(1.1,.5)
\pnode(-.5,0){a1}\pnode(.5,0){a2}
\pccurve[angleA=315,angleB=180]([nodesep=.5,angle=135]a1)(0,-.25)
\pccurve[angleA=135,angleB=0]([nodesep=.5,angle=315]a2)(0,.25)
\pccurve[border=.1,angleA=45,angleB=180]([nodesep=.5,angle=225]a1)(0,.25)
\pccurve[border=.1,angleA=225,angleB=0]([nodesep=.5,angle=45]a2)(0,-.25)
\endpspicture
= -q^{1/4}\pspicture[.4](-.9,-.5)(.9,.5)
\pnode(-.3,0){a1}\pnode(.3,0){a2}
\psbezier([nodesep=.5,angle=225]a1)([nodesep=.25,angle=225]a1)
    ([nodesep=.25,angle=135]a1)([nodesep=.5,angle=135]a1)
\pccurve[angleA=135,angleB=90,ncurv=1](a2)(-.2,0)
\qline([nodesep=.5,angle=315]a2)(a2)
\pccurve[border=.1,angleA=225,angleB=270,ncurv=1](a2)(-.2,0)
\psline[border=.1]([nodesep=.5,angle=45]a2)(a2)
\endpspicture
- q^{-1/4}\rcrossing = (1-q)\vv + q^{-1/2}\hh
$$
which leads to an expansion of tangle $c$:
$$
\pspicture[.43](-1,-.5)(1,1.1)
\pnode(-.3,1){a1}\pnode(.3,1){a2}
\pnode(0,.2){a3}\pnode(-.6,-.2){a4}
\pnode(0,0){a5}\pnode(.6,-.2){a6}
\nccurve[border=.1,angleA=0,angleB=135]{a3}{a6}
\nccurve[border=.1,angleA=315,angleB=0,ncurv=1]{a1}{a5}
\nccurve[border=.1,angleA=225,angleB=180,ncurv=1]{a2}{a5}
\nccurve[border=.1,angleA=45,angleB=180]{a4}{a3}
\endpspicture
=
(1-q)\rcrossing
+ q^{-1/2}\pspicture[.4](-.6,-.8)(.8,.6)
\pnode(0,-.3){a1}\pnode(0,.3){a2}
\psbezier([nodesep=.5,angle=225]a1)([nodesep=.25,angle=225]a1)
    ([nodesep=.25,angle=315]a1)([nodesep=.5,angle=315]a1)
\pccurve[angleA=315,angleB=0,ncurv=1](a2)(0,-.2)
\qline([nodesep=.5,angle=135]a2)(a2)
\pccurve[border=.1,angleA=225,angleB=180,ncurv=1](a2)(0,-.2)
\psline[border=.1]([nodesep=.5,angle=45]a2)(a2)
\endpspicture
 = (q^{5/4} - q^{1/4})\vv +
(q^{3/4} - q^{-1/4} + q^{-5/4})\hh
$$
which leads to the evaluation of the entire knot projection:
\begin{eqnarray*}
\pspicture[.45](-1.25,-.25)(1.25,2.25)
\pnode(-.4,2){a4}\pnode(.4,2){a8}\pnode(-.85,1.5){a5}
\pnode(-.25,1.5){a9}\pnode(.25,1.5){a3}\pnode(.85,1.5){a7}
\pnode(-.33,.5){a2}\pnode(0,.5){a6}\pnode(.33,.5){a10}\pnode(0,.25){a1}
\nccurve[border=.1,angleA=0,angleB=270]{a6}{a7}
\nccurve[border=.1,angleA=270,angleB=90]{a9}{a10}
\nccurve[border=.1,angleA=270,angleB=0]{a10}{a1}
\nccurve[border=.1,angleA=180,angleB=270]{a1}{a2}
\nccurve[border=.1,angleA=90,angleB=270]{a2}{a3}
\nccurve[border=.1,angleA=270,angleB=180]{a5}{a6}
\nccurve[border=.1,angleA=90,angleB=0]{a3}{a4}
\nccurve[border=.1,angleA=180,angleB=90]{a4}{a5}
\nccurve[border=.1,angleA=90,angleB=0]{a7}{a8}
\nccurve[border=.1,angleA=180,angleB=90]{a8}{a9}
\endpspicture
& = & (q^{5/4}-q^{1/4})
\pspicture[.4](-.6,-.5)(.6,.5)
\pccurve[angleA=315,angleB=270,ncurv=1](0,0)(.5,0)
\pccurve[angleA=135,angleB= 90,ncurv=1](0,0)(-.5,0)
\pccurve[border=.1,angleA= 45,angleB= 90,ncurv=1](0,0)(.5,0)
\pccurve[border=.1,angleA=225,angleB=270,ncurv=1](0,0)(-.5,0)
\endpspicture
 + (q^{3/4} - q^{-1/4} + q^{-5/4})
\pspicture[.42](-.9,-.9)(.9,.6)
\pccurve[angleA=0,angleB=0,ncurv=1](-.3,.5)(0,-.2)
\pccurve[border=.1,angleA=180,angleB=180,ncurv=1](.3,.5)(0,-.2)
\pccurve[angleA=180,angleB=180,ncurv=1](-.3,.5)(0,-.75)
\pccurve[angleA=0,angleB=0,ncurv=1](.3,.5)(0,-.75)
\endpspicture \\
& = & -(q^{1/2}-q^{-1/2})(q^2-q+1-q^{-1}+q^{-2})
\end{eqnarray*}

Given a knot $K$ and a clasp $c$ in one of the combinatorial spiders presented
here, consider
the clasped cabling of $K$.  \Ie, replace the strand of $K$ by several
strands tied together with $c$, for example:
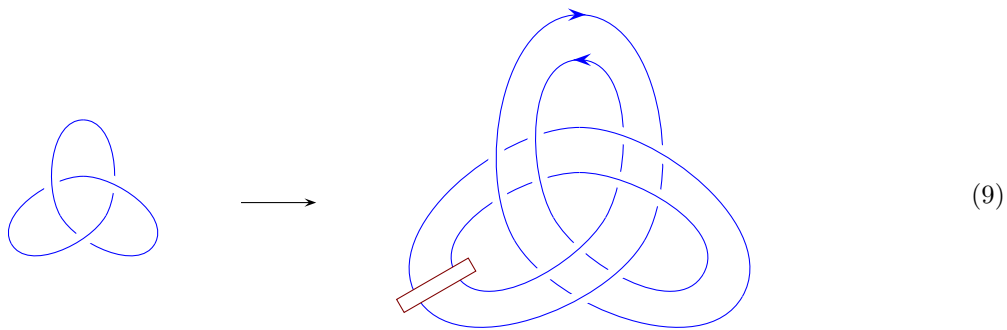
\begin{equation}
\pspicture[.5](-1.6,-1.6)(1.6,1.6)
\pccurve[border=.1,angleA=0,angleB=60](1.1;90)(.35;330)
\pccurve[border=.1,angleA=120,angleB=180](1.1;210)(.35;90)
\pccurve[border=.1,angleA=240,angleB=300](1.1;330)(.35;210)
\pccurve[border=.1,angleA=180,angleB=120](1.1;90)(.35;210)
\pccurve[border=.1,angleA=300,angleB=240](1.1;210)(.35;330)
\pccurve[border=.1,angleA=60,angleB=0](1.1;330)(.35;90)
\endpspicture
\psgoesto
\pspicture[.5](-3,-3)(3,3)
\pccurve[border=.1,angleA=0,angleB=60](2.5;90)(1;330)
\pccurve[border=.1,angleA=0,angleB=60](1.9;90)(.4;330)
\pccurve[border=.1,angleA=120,angleB=180](2.5;210)(1;90)
\pccurve[border=.1,angleA=120,angleB=180](1.9;210)(.4;90)
\pccurve[border=.1,angleA=240,angleB=300](2.5;330)(1;210)
\pccurve[border=.1,angleA=240,angleB=300](1.9;330)(.4;210)
\pccurve[border=.1,angleA=180,angleB=120](2.5;90)(1;210)
\pccurve[border=.1,angleA=180,angleB=120](1.9;90)(.4;210)
\pccurve[border=.1,angleA=300,angleB=240](2.5;210)(1;330)
\pccurve[border=.1,angleA=300,angleB=240](1.9;210)(.4;330)
\pccurve[border=.1,angleA=60,angleB=0](2.5;330)(1;90)
\pccurve[border=.1,angleA=60,angleB=0](1.9;330)(.4;90)
\pnode(2.75;210){a1}\pnode(1.65;210){a2}
\psline[arrows=->,arrowscale=1.5](2.2pt,2.5)(2.3pt,2.5)
\psline[arrows=->,arrowscale=1.5](-2.2pt,1.9)(-2.3pt,1.9)
\pspolygon[fillstyle=solid,fillcolor=white,linecolor=darkred]
    ([nodesep=.1,angle=120]a1)([nodesep=.1,angle=300]a1)
    ([nodesep=.1,angle=300]a2)([nodesep=.1,angle=120]a2)
\endpspicture\label{fccable}
\end{equation}
The value of a web such as the one in Figure~(\ref{fccable}) is also a regular
isotopy invariant of $K$; it is the Reshetikhin-Turaev ribbon graph invariant
of $K$ colored by an irreducible representation whose weight is the highest
weight of $c$ \cite{RT:cmp1990}.  By the isomorphism between combinatorial
and algebraic spiders, clasps are a complete set in the sense that these
clasped invariants yield all of the information about $K$ that can be
obtained from applying the unclasped link invariants to all cablings of of
$K$ and all other satellites of $K$.  Moreover, the value of a clasped
cabling of $K$ can be computed more easily than an unclasped cabling,
because clasped web spaces are much smaller vector spaces than their
unclasped counterparts.

\subsection{Combinatorial consequences}

Although Theorems~\ref{tha2equinum},\ref{thb2equinum}, and \ref{thg2equinum}
are properly results in representation theory, they are also interesting as
results in enumerative combinatorics.  John Stembridge and Richard Stanley
\cite{Stanley:personal} noted that a $B_2$ basis web $w \in B(11\ldots1)$ is
equivalent to a matching of $2n$ cyclically ordered points with no 6-point
star, meaning that among the $2n$ points, there are no six in cyclic order no
six points $p_1,p_2,p_3,q_1,q_2,q_3$ with each $p_i$ matched to $q_i$.  On
other hand, the number of such matchings with no $2k$-point star is known to
be $\dim \Inv(V^{\tensor 2n})$, where $V$ is the defining representation of
$\sp(2k-2)$ \cite{Sundaram:ima1990}. Thus, this special case of
Theorem~\ref{thb2equinum} was known previously.

It is interesting that ordinary trivalent graphs are related to the
exceptional Lie algebra $G_2$. One corollary of this surprising connection is
the following enumerative result:

\begin{theorem} For $n \ge 3$, let $a_n$ be the number
of triangulations of a fixed convex $n$-gon such that
at least six triangles meet at each internal vertex, and let
$$A(x) = 1 + x^2 + \sum_n a_n x^n$$
be a generating function.  For $n \ge 0$, let
$$b_n = \dim \Inv_{G_2}(V(\lambda_1)^{\tensor n})$$
for $n \ge 0$, and let
$$B(x) = \sum_n b_n x^n$$
be a generating function.  Then
$$B(x) = A(xB(x)). \label{eotrees}$$
\label{thtriangs}
\end{theorem}
\begin{proof}(Sketch) By Theorem~\ref{thg2equinum}, $b_n$ may equally well be
defined as the number of $G_2$ basis webs with $n$ endpoints.  Each
triangulation is dual to a connected basis web.  Equation~(\ref{eotrees})
becomes a standard relation between the number of connected graphs of a
certain type and the number of disjoint unions of such graphs.  More
precisely, consider $n$ boundary points on a disk with one marked, and
consider a basis web $w \in B(11\ldots1)$ with this boundary.  Then $w$ is
given by an {\em ordered tree} whose vertices are its connected components:
The root is the component with the marked vertex, and the children of each
component are the adjacent components other than the parent, if there is one.
 The children of each component are ordered going counterclockwise, with the
distiguished boundary point separating the first and last children of the
root. Equation~(\ref{eotrees}) is the usual generating function for ordered
trees \cite[p. 11-12]{Goulden-Jackson}.
\end{proof}

The sequence $\{b_n\}$ is easy to compute using character theory.  Thus,
Theorem~\ref{thtriangs} produces a fast algorithm for computing $\{a_n\}$.

It is easy to show that the radius of convergence of $B(x)$ is 1/7.  Then by
Theorem~\ref{thtriangs}, the radius of convergence of $A(x)$ is at least
$B(1/7)/7$.  Furthermore, numerical evidence supports the conjecture that it
is exactly $B(1/7)/7$, \ie, that

\begin{conjecture} If $a_n$ and $b_n$ are defined as in
Theorem~\ref{thtriangs}, then
$$\lim_{n \to \infty} \sqrt[n]{a_n} = {7\over \sum_{n=0}^\infty b_n 7^{-n}}
= 6.811\ldots$$
\end{conjecture}

Our equinumeration theorems count all basis webs with a fixed boundary. 
However, for each $n$ and $k$, one can also consider the set of $G_2$ basis
webs, for example, with $n$ endpoints of type 1 and $k$ internal vertices. 
The role of these sets and their cardinality in representation theory is not
known.

Given an $A_2$ basis-web set $B(s)$ for some sign string $s$,
there are various H-maps corresponding to different adjacent
pairs of signs.  But if $s'$ is another sign string of the same
weight as $s$, then any two sequences of H-maps permuting $s$
into $s'$ yield the same bijection $B(s) \to B(s')$, provided
that the first and the last sign of $s$ are not considered adjacent.
On the other hand, if a sign is moved all the way around
the boundary by H-maps, holonomy arises:
$$
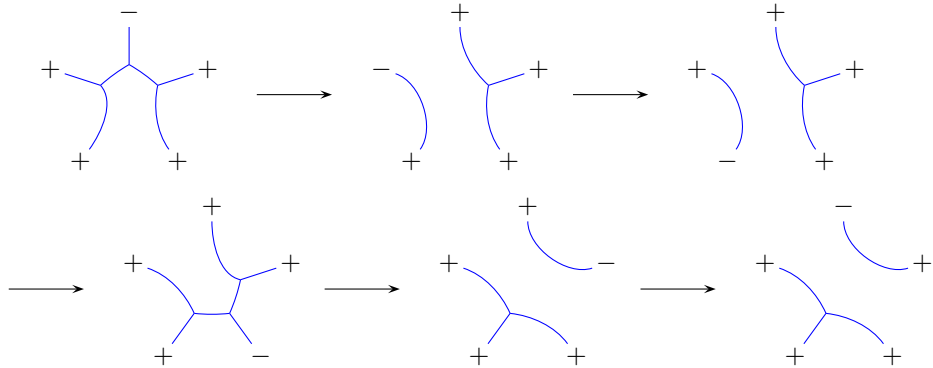
\pspicture[.5](-1.2,-1.2)(1.2,1.2)\pentanode
\ncline{a5}{b5}\ncline{a1}{b1}\ncline{a2}{b2}
\ncarc[arcangle=6]{a2}{a1}\ncarc[arcangle=6]{a1}{a5}
\nccurve[angleA=312,angleB= 54]{a2}{b3}
\nccurve[angleA=258,angleB=126]{a5}{b4}
\rput([nodesep=.2,angle= 90]b1){$-$}
\rput([nodesep=.2,angle=162]b2){$+$}
\rput([nodesep=.2,angle=234]b3){$+$}
\rput([nodesep=.2,angle=306]b4){$+$}
\rput([nodesep=.2,angle= 18]b5){$+$}
\endpspicture
\psgoesto
\pspicture[.5](-1.2,-1.2)(1,1.2)\pentanode\ncline{b5}{a5}
\nccurve[angleA=138,angleB=270]{a5}{b1}
\nccurve[angleA=342,angleB= 54]{b2}{b3}
\nccurve[angleA=258,angleB=126]{a5}{b4}
\rput([nodesep=.2,angle= 90]b1){$+$}
\rput([nodesep=.2,angle=162]b2){$-$}
\rput([nodesep=.2,angle=234]b3){$+$}
\rput([nodesep=.2,angle=306]b4){$+$}
\rput([nodesep=.2,angle= 18]b5){$+$}
\endpspicture
\psgoesto
\pspicture[.5](-1.2,-1.2)(1,1.2)\pentanode\ncline{b5}{a5}
\nccurve[angleA=138,angleB=270]{a5}{b1}
\nccurve[angleA=342,angleB= 54]{b2}{b3}
\nccurve[angleA=258,angleB=126]{a5}{b4}
\rput([nodesep=.2,angle= 90]b1){$+$}
\rput([nodesep=.2,angle=162]b2){$+$}
\rput([nodesep=.2,angle=234]b3){$-$}
\rput([nodesep=.2,angle=306]b4){$+$}
\rput([nodesep=.2,angle= 18]b5){$+$}
\endpspicture
$$
$$
\psgoesto
\pspicture[.5](-1.2,-1.2)(1,1.2)\pentanode
\ncline{a3}{b3}\ncline{a4}{b4}\ncline{a5}{b5}
\ncarc[arcangle=6]{a5}{a4}\ncarc[arcangle=6]{a4}{a3}
\nccurve[angleA=168,angleB=270]{a5}{b1}
\nccurve[angleA=114,angleB=342]{a3}{b2}
\rput([nodesep=.2,angle= 90]b1){$+$}
\rput([nodesep=.2,angle=162]b2){$+$}
\rput([nodesep=.2,angle=234]b3){$+$}
\rput([nodesep=.2,angle=306]b4){$-$}
\rput([nodesep=.2,angle= 18]b5){$+$}
\endpspicture
\psgoesto
\pspicture[.5](-1.2,-1.2)(1,1.2)\pentanode\ncline{b3}{a3}
\nccurve[angleA=354,angleB=126]{a3}{b4}
\nccurve[angleA=198,angleB=270]{b5}{b1}
\nccurve[angleA=114,angleB=342]{a3}{b2}
\rput([nodesep=.2,angle= 90]b1){$+$}
\rput([nodesep=.2,angle=162]b2){$+$}
\rput([nodesep=.2,angle=234]b3){$+$}
\rput([nodesep=.2,angle=306]b4){$+$}
\rput([nodesep=.2,angle= 18]b5){$-$}
\endpspicture
\psgoesto
\pspicture[.5](-1.2,-1.2)(1,1.2)\pentanode\ncline{b3}{a3}
\nccurve[angleA=354,angleB=126]{a3}{b4}
\nccurve[angleA=198,angleB=270]{b5}{b1}
\nccurve[angleA=114,angleB=342]{a3}{b2}
\rput([nodesep=.2,angle= 90]b1){$-$}
\rput([nodesep=.2,angle=162]b2){$+$}
\rput([nodesep=.2,angle=234]b3){$+$}
\rput([nodesep=.2,angle=306]b4){$+$}
\rput([nodesep=.2,angle= 18]b5){$+$}
\endpspicture
$$
It would be interesting to compute the cycle structure or establish
properties of this holonomy, which is related to a linear action of the
braid group on
$$\Inv((V_+ \oplus V_-)^{\tensor n}) \cong \bigoplus_s W(s),$$
where the direct sum on the right is
taken over all sign strings of length $n$.

%\bibliographystyle{plain}
%\bibliography{spider}

\end{document}